\documentclass[10pt]{article}
\usepackage[utf8x]{inputenc}
\usepackage[T1]{fontenc}
\usepackage[english]{babel}  
\usepackage{color}
\usepackage[usenames,dvipsnames]{xcolor}
\usepackage{mathptmx} 
\usepackage[scaled=0.90]{helvet}  
\usepackage{courier}
\usepackage[all]{xy}
\usepackage{enumitem}
\usepackage[a4paper,twoside]{geometry}
\usepackage{amsmath,amssymb,mathrsfs}
\usepackage[colorlinks=true]{hyperref}

\usepackage{fancyhdr} 

\usepackage{variations}
\usepackage{dsfont} 
\usepackage{float}
\usepackage{graphicx}
\usepackage{stmaryrd}
\usepackage{wrapfig} 
\usepackage{booktabs}
\usepackage[nottoc, notlof, notlot]{tocbibind}
\usepackage{bbm} 

\usepackage{amsthm}
\usepackage[capitalize,nameinlink]{cleveref}

\newtheorem{prop}{Proposition –}[section]
\newtheorem{lem}{Lemma –}[section]
\newtheorem{coro}{Corollary –}[section]

\newtheorem{theo}{Theorem –}[section]

\theoremstyle{definition}
\newtheorem{defi}{Definition –}[section]

\theoremstyle{remark}
\newtheorem{rema}{Remark –}[section]

\newenvironment{assumption}[1][\unskip]{%
\begin{enumerate}[label=(\textbf{#1}),leftmargin=60pt] 
\item\label{#1}}
{\end{enumerate}}

\definecolor{section}{RGB}{0,51,102}
\definecolor{subsection}{RGB}{64,64,64}

\usepackage{titlesec}

\titleformat*{\section}{\Large\scshape\center}
\titlespacing{\section}{\parindent}{1.5cm}{0.5cm}
\titleformat{\subsection}{\normalfont\bfseries}{\S\ \thesubsection.}{1.5em}{}
\titlespacing{\subsection}{0pc}{0.6cm}{\wordsep}
\titleformat*{\paragraph}{\normalfont\bfseries}

\newenvironment{Proof}{\vspace{1cm}\noindent\underline{Proof.~}}{\hfill$\square$\bigbreak}

\everymath{\displaystyle}
\usepackage{fancybox}

\numberwithin{equation}{section}

\newcommand{\minsep}{\noindent\begin{center}\rule{5cm}{1.2pt}\end{center}}
\newcommand{\sep}{\noindent\rule{\linewidth}{1.2pt}}

\newcommand{\e}{\varepsilon}
\newcommand{\parenth}[1]{\left( #1 \right)}
\newcommand{\dsum}{\mathop{}\displaystyle\sum}
\newcommand{\dint}{\mathop{}\displaystyle\int}

\newcommand{\suiteN}[2]{\parenth{#1}_{#2 \in \N}}

\newcommand{\beal}{\begin{align*}}
\newcommand{\eal}{\end{align*}}	
\newcommand{\beq}{\begin{equation}}
\newcommand{\eeq}{\end{equation}}
\newcommand{\bigmid}[2]{#1~\middle|~#2}

\newcommand{\ensemble}[1]{\left\lbrace #1\right\rbrace}

\newcommand{\abs}[1]{\left\vert #1 \right\vert}

\newcommand{\norme}[1]{\left\Vert #1\right\Vert}

\newcommand{\transp}[1]{{}^\textbf{T}{#1}}

\newcommand{\tr}{\mbox{Tr}}

\newcommand{\scalaire}[2]{\left\langle #1,#2 \right\rangle}

\newcommand{\quadra}[1]{\left\langle #1\right\rangle}

\newcommand{\bsys}{\left\lbrace\begin{matrix}}
\newcommand{\esys}{\end{matrix}\right.}

\renewcommand{\d}{\mathop{}\mathopen{}\mathrm{d}}

\def\Hess{\mathop{\rm Hess}\nolimits}

\newcommand{\supess}{\operatorname{esssup}}

\newcommand{\indicat}{\mathds{1}}

\newcommand{\esp}[1]{\mathbf{E}\left( #1\right)}
\newcommand{\espcond}[2]{\mathbf{E}\mathopen{}\left(#1\middle|#2\right)}

\newcommand{\espcondQi}[3]{\mathbf{E}^{\Q^{#3}}\mathopen{}\left(#1\middle|#2\right)}

\newcommand{\ps}{\text{ a.s}}
\newcommand{\quadravar}[1]{\left\langle #1 \right\rangle} 

\newcommand{\startstop}[3]{{}^{#1}#2^{\rfloor #3}}

\newcommand{\tit}[1]{\textit{#1}}

\newcommand{\tbf}[1]{\textbf{#1}}
\newcommand{\R}{\ensuremath{\mathbf{R}}}

\newcommand{\Q}{\ensuremath{\mathbf{Q}}}

\newcommand{\KK}{\ensuremath{\mathbf{K}}}

\newcommand{\N}{\ensuremath{\mathbf{N}}}
\newcommand{\PP}{\ensuremath{\mathbf{P}}}

\newcommand{\D}{\ensuremath{\mathbf{D}}}

\newcommand{\mcA}{{\mathcal A}} 
\newcommand{\mcB}{{\mathcal B}} 
\newcommand{\mcC}{{\mathcal C}}
 
\newcommand{\mcE}{{\mathcal E}} 
\newcommand{\mcF}{{\mathcal F}}
\newcommand{\mcG}{{\mathcal G}} 
 
\newcommand{\mcH}{{\mathcal H}}

\newcommand{\mcK}{{\mathcal K}} 
\newcommand{\mcL}{{\mathcal L}}
\newcommand{\mcM}{{\mathcal M}}

\newcommand{\mcP}{{\mathcal P}}

\newcommand{\mcS}{{\mathcal S}}
\newcommand{\mcT}{{\mathcal T}}

\newcommand{\mcZ}{{\mathcal Z}} 

\newcommand{\mbbA}{{\mathbb A}} 
\newcommand{\mbbB}{{\mathbb B}} 

\newcommand{\mbbD}{{\mathbb D}}

\newcommand{\mbbK}{{\mathbb K}}

\usepackage[final]{showlabels}
\usepackage{setspace}
\begin{document}
\hypersetup{colorlinks,citecolor=section,linkcolor=section,urlcolor=section}

\title{\textsc{A stability approach for solving multidimensional quadratic BSDEs}}
\maketitle
\begin{center}
\author{Jonathan Harter,\footnote{Université de Bordeaux, Institut de Mathématiques de Bordeaux UMR 5251, 351 Cours de la libération, 33 405 TALENCE, France, jonathan.harter@u-bordeaux.fr}} 
\author{Adrien Richou\footnote{Université de Bordeaux, Institut de Mathématiques de Bordeaux UMR 5251, 351 Cours de la libération, 33 405 TALENCE, France, adrien.richou@u-bordeaux.fr}}
\end{center}
\minsep

\begin{abstract}
We establish an existence and uniqueness result for a class of multidimensional quadratic backward stochastic differential equations (BSDE). This class is characterized by constraints on some uniform \tit{a priori} estimate on solutions of a sequence of approximated BSDEs. We also present effective examples of applications. Our approach relies on the strategy developed by Briand and Elie in [Stochastic Process. Appl. \tbf{123} 2921--2939] concerning scalar quadratic BSDEs.
\end{abstract}
\sep
\section{Introduction}
\paragraph{Backward Stochastic Differential Equations}
Backward stochastic differential equations (BSDEs) have been first introduced in a linear version by Bismut \cite{Bismut:1973ji}, but since the early nineties and the seminal work of Pardoux and Peng \cite{Pardoux:1990ju}, there has been an increasing interest for these equations due to their wide range of applications in stochastic control, in finance or in the theory of partial differential equations. Let us recall that, solving a BSDE consists in finding an adapted pair of processes $(Y,Z)$, where $Y$ is a $\R^d$-valued continuous process and $Z$ is a $\R^{d\times k}$-valued progressively measurable process, satisfying the equation
\begin{equation}
\label{EDSR intro}
 Y_t = \xi + \int_t^T f(s,Y_s,Z_s) \d s-\int_t^T Z_s \d W_s, \quad 0 \leq t \leq T, \quad a.s.
\end{equation}
where $W$ is a $k$-dimensional Brownian motion with filtration $(\mathcal{F}_t)_{t \in \R^+}$, $\xi$ is a $\mathcal{F}_T$-measurable random variable called the terminal condition, and $f$ is a (possibly random) function called the generator. Since the seminal paper of Pardoux and Peng \cite{Pardoux:1990ju} that gives an existence and uniqueness result for BSDEs with a Lipschitz generator, a huge amount of paper deal with extensions and applications. In particular, the class of BSDEs with generators of quadratic growth with respect to the variable $z$, has received a lot of attention in recent years. 
Concerning the scalar case,\tit{i.e.} $d=1$, existence and uniqueness of solutions for quadratic BSDEs has been first proved by Kobylanski in \cite{Kobylanski:2000hl}. Since then, many authors have worked on this question and the theory is now well understood: we refer to \cite{Kobylanski:2000hl,Tevzadze:2008co,Briand:2013cu} when the terminal condition is bounded and to \cite{Briand:2006he,Barrieu:2013gt,Delbaen:2011ji} for the unbounded case. We refer also to \cite{2014arXiv1409.5322G} for a study of BMO properties of $Z$.

In this paper we will focus on existence and uniqueness results for quadratic BSDEs in the multidimensional setting, \tit{i.e.} $d>1$. Let us remark that, in addition to its intrinsic mathematical interest, this question is important due to many applications of such equations. We can mention for example following applications: nonzero-sum risk-sensitive stochastic differential games  in \cite{ELKaroui:2003hh,Hu:2016is}, financial market equilibrium problems for several interacting agents in \cite{Espinosa:2015gsa,Frei:2011kk,Frei:2014hb,Bielagk:2015wf}, financial price-impact models in \cite{Kramkov:2016hc,Kramkov:2016ck}, principal agent contracting problems with competitive interacting agents in \cite{2016arXiv160508099E}, stochastic equilibria problems in incomplete financial markets \cite{2015arXiv150507224K,2016arXiv160300217X} or existence of martingales on curved spaces with a prescribed terminal condition \cite{Darling:1995wd}.

Let us note that moving from the scalar framework to the multidimensional one is quite challenging since tools usually used when $d=1$, like monotone convergence or Girsanov transform, can no longer be used when $d>1$. Moreover, Frei and dos Reis provide in \cite{Frei:2011kk} an example of multidimensional quadratic BSDE with a bounded terminal condition and a very simple generator such that there is no solution to the equation. This informative counterexample show that it is hopeless to try to obtain a direct generalization of the Kobylanski existence and uniqueness theorem in the multidimensional framework or a direct extension of the Pardoux and Peng existence and uniqueness theorem for locally-Lipschitz generators. Nevertheless, we can find in the literature several papers that deal with special cases of multidimensional quadratic BSDEs and we give now a really brief summary of them.

First of all, a quite general result was obtain by Tevzadze in \cite{Tevzadze:2008co}, when the bounded terminal condition is small enough, by using a fixed-point argument and the theory of BMO martingales. Some generalizations with somewhat more general terminal conditions are considered in \cite{Frei:2014hb,Kramkov:2016ck}. In \cite{Cheridito:2015es}, Cheridito and Nam treat some quadratic BSDEs with very specific generators. Before these papers, Darling was already able to construct a martingale on a manifold with a prescribed terminal condition by solving a multidimensional quadratic BSDE (see \cite{Darling:1995wd}). Its proof relies on a stability result obtained by coupling arguments. Recently, the so-called quadratic diagonal case has been considered by Hu and Tang in \cite{Hu:2016is}. To be more precise, they assume that the nth line of the generator has only a quadratic growth with respect to the nth line of $Z$. This type of assumption allows authors to use Girsanov transforms in their \tit{a 
priori} estimates calculations. Some 
little bit more general assumptions are treated by Jamneshan, Kupper and Luo in \cite{Jamneshan:2014uia} (see also \cite{2015arXiv150501796L}). Finally, in the very recent paper \cite{2016arXiv160300217X}, Xing and {\v Z}itkovi{\'c} obtained a general result in a Markovian setting with weak regularity assumptions on the generator and the terminal condition. Instead of assuming some specific hypotheses on the generator, they suppose the existence of a so called Liapounov function which allows to obtain a uniform \tit{a priori} estimate on some sequence $(Y^n,Z^n)$ of approximations of $(Y,Z)$. Their approach relies on analytic methods. We refer to this paper for references on analytic and PDE methods for solving systems of quadratic semilinear parabolic PDEs.

\paragraph{Our approach}
Our approach for solving multidimensional quadratic BSDEs relies on the theory of BMO martingales and stability results as in \cite{Briand:2013cu}. To get more into the details about our strategy, let us recall the sketch of the proof used by Briand and Elie in \cite{Briand:2013cu}. The generator $f$ is assumed to be locally Lipschitz and, to simplify, we assume that it depends only on $z$. First of all, they consider the following approximated BSDE
\begin{equation*}
 Y_t^M = \xi + \int_t^T f(\rho^M(Z_s^M)) ds-\int_t^T Z_s^M dW_s, \quad 0 \leq t \leq T, \quad a.s.
\end{equation*}
where $\rho^M$ is a projection on the centered Euclidean ball of radius $M$. Then existence and uniqueness of $(Y^M,Z^M)$ is obvious since this new BSDE has a Lipschitz generator.
Now, if we assume that $\xi$ is Malliavin differentiable with a bounded Malliavin derivative, they show that $Z^M$ is bounded uniformly with respect to $M$. Thus, $(Y^M,Z^M)=(Y,Z)$ for $M$ large enough. Importantly, the uniform bound on $Z^M$ is obtained thanks to a uniform (with respect to $M$) \tit{a priori} estimate on the BMO norm of the martingale $\int_0^. Z^M_s dW_s$. Subsequently, they extend their existence and uniqueness result for a general bounded terminal condition: $\xi$ is approximated by a sequence $(\xi^n)_{n \in \N}$ of bounded terminal conditions with bounded Malliavin derivatives and they consider $(Y^n,Z^n)$ the solution of the following BSDE
\begin{equation*}
 Y_t^n = \xi^n + \int_t^T f(Z_s^n) ds-\int_t^T Z_s^n dW_s, \quad 0 \leq t \leq T, \quad a.s. 
\end{equation*}
By using a stability result for quadratic BSDEs, they show that $(Y^n,Z^n)$ is a Cauchy sequence that converges to the solution of the initial BSDE \eqref{EDSR intro}. Once again, the stability result used by Briand and Elie relies on a uniform (with respect to $n$) \tit{a priori} estimate on the BMO norm of the martingale $\int_0^. Z^n dW_s$.

The aim of this paper is to adapt this approach in our multidimensional setting. In the first approximation step, we are able to show that $Z^M$ is bounded uniformly with respect to $M$ if we have a small enough uniform (with respect to $M$) \tit{a priori} estimate on the BMO norm of the martingale $\int_0^. Z^M_s dW_s$. But, contrarily to the scalar case, it is not possible to show that we have an \tit{a priori} estimate on the BMO norm of the martingale $\int_0^. Z^M_s dW_s$ under general quadratic assumptions on the generator (let us recall the counterexample provides by Frei and dos Reis in \cite{Frei:2011kk}). So, this \tit{a priori} estimate on the BMO norm of the martingale $\int_0^. Z^M_s dW_s$ becomes in our paper an \tit{a priori} assumption and this assumption has to be verified on a case-by-case basis according to the BSDE structure. In the second approximation step, we are facing the same issue: we are able to show the existence and uniqueness of a solution to \eqref{EDSR intro} by using a 
stability result if we have a small 
enough uniform (with respect to $n$) \tit{a priori} estimate on the BMO norm of the martingale $\int_0^. Z^n_s dW_s$, and this \tit{a priori} estimate becomes, once again, an assumption that has to be verified on a case-by-case basis according to the BSDE structure. Let us emphasize that the estimate on the boundedness of $Z^M$ and the stability result used in the second step come from an adaptation of results obtained by Delbaen and Tang in \cite{Delbaen:2008bh}.
The fact that our results are true only when we have a small enough uniform estimate on the BMO norm of the martingale $\int_0^. Z^M_s dW_s$ or $\int_0^. Z^n_s dW_s$ is the main limitation of our results. Nevertheless, we emphasize that this limitation is related to a crucial open question that could be independently investigated. To be precise, we would like to know if the classical reverse H\"older inequality for exponential of BMO martingales (see Theorem 3.1 in \cite{Kazamaki:2007tr}) stays true in a multidimensional setting, i.e. when we have a matrix valued BMO martingale. For further details we refer the reader to Remark \ref{lem question ouverte}.

To show the interest of these theoretical results, we have to find now some frameworks for which we are able to obtain estimates on the BMO norm of martingales $\int_0^. Z^M_s dW_s$ and $\int_0^. Z^n_s dW_s$. This is the purpose of Section \ref{applications} where results of \cite{Tevzadze:2008co,Darling:1995wd,Hu:2016is} are revisited. Let us note that one interest of our strategy comes from the fact that we obtain these estimates by very simple calculations that allow to easily get new results: for example, we are able to extend the result of Tevzadze when the generator satisfies a kind of monotone assumption with respect to $y$ (see subsection \ref{subsection monotone}). Moreover, we can remark that obtaining such estimates is strongly related to finding a so-called Liapounov function in \cite{2016arXiv160300217X}. Result on the boundedness of $Z$ is also interesting in itself since it allows to consider the initial quadratic BSDE \eqref{EDSR intro} as a simple Lipschitz one which gives access to numerous 
results 
on Lipschitz BSDEs: numerical 
approximation schemes, differentiability, stability, and so on.

\paragraph{Structure of the paper}
In the remaining of the introduction, we introduce notations, the framework and general assumptions. We have collected in Section 2 all our main results in order to improve the readability of the paper. Section 3 contains some general results about SDEs and linear BSDEs adapted from \cite{Delbaen:2008bh}. Section 4 is devoted to the proof of stability properties, existence and uniqueness theorems for multidimensional quadratic BSDEs. Finally, proofs of the applications of previous theoretical results are given in Section 5.

\subsection{Notations}

\begin{itemize}[label=$\diamond$, leftmargin=*, noitemsep]
	\item Let $T >0$. We consider $\left(\Omega, \mcF,(\mcF_t)_{t \in [0,T]},\PP\right)$ a complete probability space where $(\mcF_t)_{t \in [0,T]}$ is a Brownian filtration satisfying the usual conditions. In particular every càdlàg process has a continuous version. Every Brownian motion will be considered relatively to this filtered probability space. A $k$-dimensional Brownian motion $W=\parenth{W^i}_{1 \leqslant i \leqslant k}$ is a process with values in $\R^k$ and with independent Brownian components. Almost every process will be defined on a finite horizon $[0,T]$, either we will precise it explicitly. The stochastic integral of an adapted process $H$ will be denoted by $H \star W$, and the Euclidean   quadratic variation by $\quadravar{.,.}$. The Dolean-Dade exponential of a continuous real local martingale $M$ is denoted by $$\mcE(M) := \exp\left( M-\frac{1}{2} \langle M,M \rangle \right).$$ 
	\item \textit{Linear notions –} 
	On each $\R^p$, the scalar product will be simply denoted by a dot, including the canonical scalar product on $\mcM_{dk}(\R)$: 
	 $$ M . N = \sum_{1\leqslant i \leqslant d,1\leqslant j \leqslant k} M_{i,j} N_{i,j}.$$
	For $A \in \mcM_{dk}(\R)$, $A^{(:,p)}$ will be the column $p \in \ensemble{1,...,k}$ of $A$, and $A^{(l,:)}$ the line $l \in \ensemble{1,...,d}$.
	If $B \in \mcL(\R^{d\times k}, \R^d)$, we write for $i \in \left\lbrace 1, ..., k\right\rbrace$, $B^{(:,i,:)} \in \mcL(\R^d, \R^d)$ the linear map such that $Bx=\sum_{i=1}^k B^{(:,i,:)}x^{(:,i)}$ for all $x \in \R^{d\times k}$.
	If $A$ and $B$ are two processes with values in $\mcM_{dk}(\R)$ and $\R^k$, the quadratic variation $\scalaire{A}{B}$ is the $\R^d$ vector process
	$$ \parenth{\sum_{l=1}^d \scalaire{A^{il}}{B^l}}_{i=1}^d$$
	 and we have the integration by part formula $ \d \parenth{AB}=\d A . B + A . \d B + \d \scalaire{A}{B}$. We can also define the covariation of $(A,B) \in \mcM_{dk}(\R)\times \mcM_{kd'}(\R)$ by 
	 $$ \parenth{\sum_{l=1}^d \scalaire{A^{il}}{B^{lj}}}_{i,j=1}^{d,d'}.$$
	\item \textit{Functional spaces –}  In a general way, Euclidean norms will be denoted by $\abs{.}$ while norms relatively to $\omega$ and $t$ will be denoted by $\norme{.}$.\\ 
For a $\mcF$-adapted continuous process $Y$ with values in $\R^d$ and $1 \leqslant p \leqslant \infty$ , let us define 
	 $$ \norme{Y}_{\mcS^p} = \esp{\sup_{0 \leqslant s \leqslant T} \abs{Y_s}^p}^{1/p}, \quad \text{and} \quad \norme{Y}_{\mcS^\infty} = \supess \sup_{0 \leqslant s \leqslant T} \abs{Y_s}.$$
	If $Z$ is a random variable with values in $\R^d$, we define
	$$ \norme{Z}_{L^p}=\esp{\abs{Z}^p}^{1/p}.$$
	A  continuous martingale $M$ with values in $\R$ is in $\mcH^p(\R)$, or only $\mcH^p$ when it is not necessary to specify the state space, if $\sqrt{\langle M , M \rangle_T}\in L^p$. And we define the $\mcH^p$ norm by
	$$ \norme{M}_{\mcH^p} := \esp{\scalaire{M}{M}^{p/2}_T}^{1/p} < \infty.$$
	If $M$ is a martingale with values in $\R^d$, $M$ is in $\mcH^p(\R^d)$ if $\abs{M}$ is in $\mcH^p(\R)$.
	A real martingale $M=(M_t)_{0 \leqslant t \leqslant T}$ is said to be BMO (bounded in mean oscillation) if there exists a constant $C \geqslant 0$ such that for every stopping time $0 \leqslant \tau \leqslant T$:
	$$\espcond{(M_T-M_{\tau})^2}{\mcF_\tau}\leqslant C^2 \ps.$$
	The best constant $C$ is called the BMO norm of $M$, denoted by $\norme{M}_{BMO(\PP)}$ or sometimes only $\norme{M}_{BMO}$. 	 
   In particular, the one dimensional local martingale $Z \star W^1=\displaystyle\int_0^. Z_s \d W^1_s$ is BMO if there exists a constant $C \geqslant 0$ such that, for all stopping time $\tau$ with values in $[0,T]$, we have
   $$ \espcond{\int_\tau^T \abs{Z_s}^2 \d s}{\mcF_\tau}\leqslant C^2 \ps.$$
   In the sequel, to simplify notations  we will skip the superscript ${.}^1$ on the Brownian motion after a star.
	 For more details about BMO martingales, we can refer to \cite{Kazamaki:2007tr}. 
	 
	 For $k \geqslant 1$, $\mcC^\infty_b(\R^k)$ is the set of all $\mcC^\infty$ functions with values in $\R$ defined on $\R^k$, which have bounded derivatives. 
	 
	 Given  $b_0 \in \R^d$ and a sequence $(\alpha_n) \in (0,1]^\N$, a function  $g : [0,T] \times \R^d \rightarrow \R^d$ is said to be in $\mcC^{(\alpha_n),loc}_{b_0}$ if there exists a sequence $(c_n)$ of positive constants, such that, for all $n \in \N$, 
	 $$\sup_{(t,x) \in [0,T] \times B_n(b_0)} \abs{g(t,x)} +\sup_{(t,x) \neq (t',x') \in [0,T] \times B_n(b_0)} \frac{\abs{g(t,x)-g(t',x')}}{\abs{t-t'}^{\alpha_n/2} +\abs{x-x'}^{\alpha_n}} \leqslant c_n,$$
	 where $B_n(b_0)$ states for the Euclidean ball on $\R^d$ of center $b_0$ and radius $n$. If the last term does not depend on $b_0$, we shall say that $g$ is in $\mcC^{(\alpha_n),loc}$. Finally, for a given $\alpha \in (0,1]$, a function $g : [0,T] \times \R^d \rightarrow \R^d$ is said to be in $\mcC^{\alpha}$ if
	 $$\sup_{(t,x) \neq (t',x') \in [0,T] \times \R^d} \frac{\abs{g(t,x)-g(t',x')}}{\abs{t-t'}^{\alpha/2} +\abs{x-x'}^\alpha}<+\infty.$$
	 
	 \begin{rema}\label{loca Hold et bornée implique glob Hold}
	We can plainly show that if there exists $(\alpha_n)\in (0,1]^\N$ such that a bounded solution $v$ is in $\mcC^{(\alpha_n),loc}$, then $v \in \mcC^{\alpha_1}$.
\end{rema}
	 

	 
	\item \textit{Inequalities –} 
 	BDG inequalities claim that $\norme{.}_{\mcS^p}$ and $\norme{.}_{\mcH^p}$ are equivalent on martingale spaces with two universal constants denoted $C'_p,C_p$. It means that for all continuous local martingales $M$ vanishing at $0$,
	$$ \norme{M}_{\mcH^p} \leqslant C_p \norme{M}_{\mcS^p} $$
	and 
	$$  \norme{M}_{\mcS^p} \leqslant C'_p\norme{M}_{\mcH^p}.$$
	In \cite{Marinelli:2016du}, Marinelli and Röckner deal with martingales taking values in a separable Hilbert space. In particular, the upper constant $C'$ (Proposition 2.1 and Proposition 3.1) defined below is valid for all dimensions:
	$$ C'_p=\begin{cases}
\parenth{\frac{p}{p-1}}^{\frac{p}{2}} \parenth{\frac{p(p-1)}{2}}^2  &\text{if } p>2,\\
4 \sqrt{\frac{2}{p}} &\text{if } p<2,\\
4 &\text{if } p=2.
\end{cases}$$
	We remark that in the case $p=2$, the scalar BDG constant is valid.
	In the following every BDG inequality should be understood with this choice of $C'$. The Doob maximal inequality claims that for every $\R^d$-valued martingale $M$ and $p >1$, 
	$$ \norme{M}_{\mcS^p} \leqslant \frac{p}{p-1} \norme{M_T}_{L^p},$$
	and for $p=\infty$,
	$$ \norme{M}_{\mcS^\infty} \leqslant \norme{M_T}_{L^\infty}.$$
	If $p \in ]1,\infty[$, we will denote by $p^*$ the conjugated exponent of $p$ such that $ \frac{1}{p}+\frac{1}{p^*}=1$. 
	Finally, we say that a process $L=(L_t)_{0 \leqslant t \leqslant T}$ with values in $\R^d$ satisfies a reverse Hölder inequality for some integer $1 \leqslant p < \infty$ if there exists some constant $K_p$ such that for every stopping time $0 \leqslant \tau \leqslant T \ps$, 
	 $$ \espcond{\abs{L_T}^p}{\mcF_\tau} \leqslant K_p \abs{L_\tau}^p \ps.$$ 
	\item \textit{BMO martingales properties –} 
	We recall here several results on BMO martingales that will be useful in the sequel.
	The energy inequality (see \cite{Kazamaki:2007tr}) tells us that for every BMO martingale $M$ and every integer $n \geqslant 1$, we have 
	\beq \esp{\langle M \rangle_T^n} \leqslant n! \norme{M}^{2n}_{BMO},\label{energy}\eeq
	and a conditional version of this inequality is also true: for all $t \in [0,T]$, 
	\beq \espcond{\parenth{\langle M \rangle_T-\langle M \rangle_t}^n}{\mcF_t} \leqslant n! \norme{M}^{2n}_{BMO}.\label{condenergy} \eeq
	Consequently the space of BMO martingales is a subset of $\bigcap_{p \geqslant 1}\mcH^p$. We recall also the so-called Fefferman inequality: for $X \in \mcH^1$ and $Y \in BMO$,
	$$ \esp{\int_0^T \abs{\d \scalaire{X}{Y}_s}}\leqslant \norme{X}_{\mcH^1} \norme{Y}_{BMO}.$$
	This inequality yields the following technical lemma (see \cite{BANUELS:2010vt} and \cite{Delbaen:2008bh} for more details). 
	
	\begin{lem}\label{BMOineq}
	Let $m \geqslant 1$. We consider $X$ an adapted process and $M$ a local martingale. 
	\begin{itemize}
	\item[\textbf{(i)}] If $X \in \mcS^m$ and $M \in BMO$, then $X \star M \in \mcH^m$ and 
	$$ \norme{X \star M}_{\mcH^m} \leqslant \sqrt{2} \norme{X}_{\mcS^m} \norme{M}_{BMO}.$$
	\item[\textbf{(ii)}] If $X \in \mcH^m$ and $M \in BMO$, then $\scalaire{X}{M}_T\in L^m$ and 
	$$ \norme{\scalaire{X}{M}_T}_{L^m} \leqslant \sqrt{2} m \norme{X}_{\mcH^m} \norme{M}_{BMO}.$$
	\end{itemize}
	\end{lem}
	The John-Nirenberg inequality gives a useful estimation on exponential moments of $\langle Z \star W \rangle_T$: if $\norme{Z \star W}_{BMO} <1$, for any stopping time $\tau \in [0,T]$ we have 
	\beq \espcond{e^{\int_\tau^T \abs{Z_s}^2 \d s}}{\mcF_\tau} \leqslant \frac{1}{1-\norme{Z \star W}_{BMO}^2}. \label{JohnNir}\eeq
	We have also a result about changes of probability law and equivalence of BMO norms on a BMO ball (see Lemma A.4 in \cite{Hu:2016is} and Theorem 3.6 in \cite{Kazamaki:2007tr}). 
	
	\begin{prop}\label{BMOequiv}
		Let $B >0$. There are two constants $c_1>0$ and $c_2 >0$ depending only on $B$, such that for any BMO martingale $M$, we have for any  BMO martingale $N$ such that $\norme{N}_{BMO(\PP)} \leqslant B$, \\
		$c_1 \norme{M}_{BMO(\PP)} \leqslant  \norme{M-\langle M,N \rangle}_{BMO(\Q)} \leqslant  c_2 \norme{M}_{BMO(\PP)},$ where $\d \Q = \mcE(N)_T \d \PP$. 
	\end{prop}

	To conclude this paragraph, let us show a technical proposition that will be useful in this paper. 
	\begin{prop}\label{BMOconv}
	Let $m \geqslant 1$ and a sequence of BMO-uniformly bounded local martingales $\suiteN{Z^n \star W}{n}$.  We denote $K = \sup_{n \in \N} \norme{Z^n\star W}_{BMO}<\infty$ and assume that $Z^n\star W$ converge in $\mcH^m$ to a martingale $Z \star W$. Then $Z \star W$ is BMO too and satisfies the same inequality $\norme{Z \star W}_{BMO} \leqslant K$.  
	\end{prop}

	\begin{Proof} 
Let us define by $\mcM$ the measure $\d \mcM = \d \PP \otimes \d x$. Firstly we show that convergence in $\mcH^m$ implies the convergence for the measure $\mcM$. 
Indeed, if $m\geqslant 2$, the Jensen inequality gives us
	$$ \esp{\int_0^T \abs{Z^n_s-Z_s}^2 \d s}\leqslant \norme{Z^n \star W-Z \star W}_{\mcH^m}^2,$$
	and thus we get the convergence in measure, since for all $\varepsilon >0$,
	$$ \mcM (\abs{Z^n-Z} > \varepsilon) \leqslant \frac{1}{\varepsilon^2} \norme{Z^n \star W-Z \star W}_{\mcH^2}^2 \leqslant \frac{1}{\varepsilon^2}  \norme{Z^n \star W-Z \star W}_{\mcH^m}^2. $$ 
Moreover, if $m<2$ we also have
	\begin{align*}
		\mcM (\abs{Z^n-Z} > \varepsilon) \leqslant & \frac{1}{\varepsilon^m} \esp{\int_0^T \abs{Z^n_s-Z_s}^m \d s}\leqslant \frac{T^{1-m/2}}{\varepsilon^m} \esp{\parenth{\int_0^T \abs{Z^n_s-Z_s}^2 \d s}^{m/2}}.
		\end{align*}
	For the both cases, we get convergence in measure. Hence there exists a subsequence $\suiteN{n_k}{k}$ such that 
	$$ \abs{Z^{n_k}}^2 \underset{k \rightarrow \infty}{\longrightarrow}  \abs{Z}^2 \quad  \mcM- \text{a.e.}$$ 
	The Fatou lemma gives us for all stopping time $\tau \in [0,T]$, 
	$$ \int_\tau^T \abs{Z_s}^2 \d s\leqslant \liminf_{k \rightarrow \infty} \int_\tau^T \abs{Z^{n_k}_s}^2 \d s \quad \ps, $$
	and taking the conditional expectation 
	$$ \espcond{\int_\tau^T \abs{Z_s}^2 \d s}{\mcF_\tau} \leqslant \espcond{\liminf_{k \rightarrow \infty} \int_\tau^T \abs{Z^{n_k}_s}^2 \d s}{\mcF_\tau}\leqslant \liminf_{k \rightarrow \infty}  \espcond{\int_\tau^T \abs{Z^{n_k}_s}^2 \d s}{\mcF_\tau}\leqslant K \quad \ps .$$
	Finally $Z \star W$ is BMO and $ \norme{Z \star W}_{BMO} \leqslant K.$ 
	\end{Proof}
	
	\item \textit{Sliceability –} For a process $X$ and a stopping time $\tau$ we denote by ${}^\tau X$ the process started at time $\tau$, that is ${}^\tau X = X_{\max(., \tau)}-X^{\rfloor \tau}$ where $X^{\rfloor \tau}$ is the process stopped at $\tau$. For two stopping times $\tau \leqslant \sigma \ps$, we denote by ${}^{\tau} X^{\rfloor \sigma}$ the process started at $\tau$ and stopped at $\sigma$: 
$$ {}^\tau X^{\rfloor \sigma} = \left({}^\tau X\right)^{\rfloor \sigma}.$$
Associativity property of the stochastic integral can be rewritten with this notation:
$$ \startstop{\tau}{\parenth{H\star W}}{\sigma}=H\star \startstop{\tau}{W}{\sigma}.$$
Between $\tau$ and $\sigma$, the started and stopped process is simply a translation of the stopped process: for all $u$ such that  $\tau \leqslant u \leqslant \sigma \ps$, 
$$\startstop{\tau}{X}{\sigma}_u=X_u-X_\tau.$$

This process is constant after $\sigma$ and vanishes before $\tau$. 
Let us suppose that $X$ is a BMO martingale. We say that $X$ is $\varepsilon$-sliceable if there exists a subsequence of stopping times $0=T_0 \leqslant T_1 \leqslant ... \leqslant T_N=T$, where $N \in \N$ is deterministic, such that 	
$$\norme{{}^{T_n}(X)^{\rfloor T_{n+1}}}_{BMO} \leqslant \varepsilon. $$
The set of all $\varepsilon$-sliceable processes will be denoted by $BMO_\varepsilon$. Schachermayer proved in \cite{Schachermayer:1996tj} that 
$$ \bigcap_{\varepsilon >0} BMO_{\varepsilon}= \overline{\mcH^\infty}^{BMO}.$$

Moreover the BMO norm of a started and stopped stochastic integral process $\startstop{\tau}{Z \star W}{\sigma}$ has a simple expression:
\begin{prop}\label{BMOstartstop}
	$$ \norme{\startstop{\tau}{Z \star W}{\sigma}}_{BMO}=\supess \sup_{\tau ' \in \mcT_{\tau,\sigma}} \espcond{\int_{\tau'}^\sigma \abs{Z_s}^2 \d s}{\mcF_{\tau'}},$$
$$ \text{ where } \mcT_{\tau, \sigma}= \left\lbrace \tau' \text{ stopping time } :  \tau \leqslant \tau' \leqslant \sigma \ps \right\rbrace.$$
\end{prop}
A proof of this proposition is given in the appendix part. 

\item \textit{Malliavin calculus –}
	We denote by
	$$ \mcP = \ensemble{f((g_1 \star W)_{T}, ..., (g_n \star W)_{T}): f \in \mcC^\infty_b(\R^n), g_i \text{ adapted }, n \geqslant 1},$$
	the set of all Wiener functions. For $F \in \mcP$, the Malliavin derivative of $F$ is a progressively measurable process $\D F \in L^2 \big([0,T] \times \Omega, \mcB ([0,T]) \otimes \mcF, \d x \otimes \d \PP\big)$ defined by
	$$ \D_t F = \sum_{i=1}^n \partial_i f((g_1 \star W)_{T}, ..., (g_n \star W)_{T}) g_i(t).$$
	In particular $\D \left((h \star W)_{T}\right)=h$ for all adapted process $h$. We define a kind of Sobolev norm on $\mcP$ with the following definition
	$$ \norme{F}_{1,2}:=\left[\esp{\abs{F}^2}+\esp{\norme{\D F}^2_{L^2(\d x)}} \right]^{1/2}.$$
	We can show that $\D$ is closable, consequently it is possible to extend the definition of $\D$ to $\mbbD^{1,2}=\overline{\mcP}^{1,2}$. Besides, $\mbbD^{1,2}$ is dense in $L^2(\Omega)$. For further considerations on Malliavin calculus  we can refer to \cite{SpringerVerlagGmbHandCoKG:2005we}. 
	We finish this paragraph by the following useful result proved in \cite{SpringerVerlagGmbHandCoKG:2005we} (Proposition 1.2.4). 
	
	\begin{prop}\label{propMaliav}
		Let $\varphi : \R^d \rightarrow \R$. We assume that there exists a constant $K$ such that for all $x,y \in \R^d$,		
		$$ \abs{\varphi(x)-\varphi(y)}\leqslant K\abs{x-y}.$$
		Let $(F^1,...,F^d)$ a vector in $\mbbD^{1,2}(\R^d)\cap L^\infty(\Omega)$. Then $\varphi(F) \in \mbbD^{1,2}(\R^d)$ and there exists a random vector $(G^1,...,G^d)$ such that
		$$ \D \varphi(F)=\sum_{i=1}^d G^i \D F^i, \quad \text{and} \quad \abs{G} \leqslant K.$$
	\end{prop}
\end{itemize}

\subsection{Framework and first assumptions}
In this paper we consider the following quadratic BSDE on $\R^d$: 
\beq Y_t = \xi+\int_t^T f(s,Y_s,Z_s) \mathrm{d}s - \int_t^T Z_s \mathrm{d} W_s, \quad 0 \leqslant t \leqslant T, \quad a.s. \label{quadra}\eeq
where $f$ is a random function $\Omega \times [0,T]\times \R^d \times \R^{d\times k} \rightarrow \R^d$ called the \textit{generator} of the BSDE such that for all $(y,z) \in \R^d \times \R^{d \times k}$ and $t \in [0,T]$, $\parenth{f(t,y,z)}_{0 \leqslant t \leqslant T}$ is progressively measurable, $(Y,Z)$ is a process with values in $\R^{d} \times \R^{d \times k}$ and $\xi \in L^2 \left( \mcF_T, \R^d \right)$. 
	
\begin{defi}
	A solution of BSDE \eqref{quadra} is a process $(Y,Z) \in\mcS^2(\R^{d}) \times \mcH^2(\R^{d\times k})$ satisfying usual integrability conditions and solving initial BSDE:
		\begin{itemize}
		\item[\tbf{(i)}] $\displaystyle\int_0^T \parenth{\abs{f(s,Y_s,Z_s)}^2+\abs{Z_s}^2} \d s < \infty$ \ps.,
		\item[\tbf{(ii)}] $Y_t = \xi+\displaystyle\int_t^T f(s,Y_s,Z_s) \mathrm{d}s - \displaystyle\int_t^T Z_s \mathrm{d} W_s, \quad 0 \leqslant t \leqslant T,$ \ps .
	\end{itemize}
\end{defi}

Some locally Lipschitz assumptions on $f$ and integrability assumptions on $\xi$ and $f$ will be assumed all along this paper. 

\begin{assumption}[H]
\begin{enumerate}[label=(\roman*)]
\item[\tbf{(i)}] For all $(y,y',z,z') \in \left(\R^d\right)^2 \times \left(\R^{d \times k}\right)^2$, we assume that there exists $(K_y,L_y,K_z,L_z) \in (\R^+)^4$ such that $\PP-\ps$ for all $t \in [0,T]$: 
	$$ \vert f(t,y,z)-f(t,y',z) \vert \leqslant (K_y+L_y\vert z \vert^2) \vert y - y' \vert,$$
	$$\vert f(t,y,z)-f(t,y,z') \vert \leqslant \left(K_z+L_z(\vert z \vert +\vert z'\vert )\right) \vert z - z' \vert , $$
\item[\tbf{(ii)}] $\esp{\abs{\xi}^2+ \int_0^T \abs{f(s,0,0)}^2\d s } <+\infty$.
\end{enumerate}
\end{assumption}

We denote by $\mbbB^m(L_y,L_z)$ the following quantity depending on $L_y$ and $L_z$:
\begin{equation}
\label{definition Bm}
\mbbB^m(L_y,L_z):=\begin{cases}
		\frac{-L_zC'_m+\sqrt{mL_y+(L_z C'_m)^2}}{\sqrt{2} mL_y} &\text{ if } L_y \neq 0,\\ \frac{1}{2\sqrt{2} L_z C'_m}&\text{ if } L_y = 0.
\end{cases}
\end{equation}

For all $m >1$, let us denote by $\mcZ_{BMO}^m$ the set 
$$ \mcZ_{BMO}^m = \ensemble{Z, \quad \R^{d \times k}-\text{valued }\text{ process } \Bigg\slash 2mL_y \norme{\abs{Z} \star W}^2_{BMO}+2\sqrt{2}L_z C'_m \norme{\abs{Z} \star W}_{BMO} <1},$$
which can be rewritten as
$$ \mcZ_{BMO}^m = \ensemble{Z, \quad \R^{d \times k}-\text{valued }\text{ process } \Bigg\slash \norme{\abs{Z} \star W}_{BMO} < \mbbB^m(L_y,L_z) },$$
where $\mbbB^m(L_y,L_z)$ is defined in \eqref{definition Bm}.
%
%
We also denote by $\mcZ_{BMO}^{slic,m}$ the set of all $ \R^{d \times k}$-valued processes $Z$ for which there exists a sequence of stopping times $0=T_0\leqslant T_1 \leqslant ...\leqslant T_N=T$ such that $ \startstop{T_{i}}{Z}{T_{i+1}} \in \mcZ_{BMO}^m$ for all $i \in \{0,...,N\}$. ~
\vspace{1cm}

To conclude this introduction, we finally consider an approximation of the BSDE \eqref{quadra}. To this purpose
let us introduce a localisation of $f$ defined by $f^M(t,y,z)=f(t,y, \rho^M(z))$ where $\rho^M : \R^{d \times k} \rightarrow \R^{d \times k}$ satisfies the following properties :
\begin{itemize}
	\item $\rho^M$ is the identity on $\mcB_{\R^{d\times k}}(0,M)$,
	\item $\rho^M$ is the projection on $\mcB_{\R^{d\times k}}(0,M+1)$ outside $\mcB_{\R^{d\times k}}(0,M+1)$ ,
	\item $\rho^M$ is a $\mcC^\infty$ function with $\vert \nabla \rho^M (z) \vert\leqslant 1$ for all $z \in \R^{d \times k}$. 
\end{itemize}  
Thus $f^M$ is a globally Lipschitz function with constants depending on $M$. Indeed we have for all $(t,y,y',z,z') \in [0,T] \times (\R^d)^2 \times \parenth{\R^{d \times k}}^2$, 
\begin{align*}
	\abs{f^M(t,y,z)-f^M(t,y',z')} \leqslant & \abs{f\left(t,y,\rho^M(z) \right)-f\left(t,y',\rho^M(z) \right)}+\abs{f\left(t,y',\rho^M(z) \right)-f\left(t,y',\rho^M(z') \right)}\\\leqslant & \left(K_y+L_y\vert\rho^M(z)\vert^2\right) \abs{y-y'} + \left(K_z+L_z\left(\vert\rho^M(z)\vert+\vert\rho^M(z')\vert\right)\right)\abs{z-z'}
	\\\leqslant & (K_y+L_y (M+1)^2) \abs{y-y'}+(K_z+2L_z(M+1)) \abs{z-z'} .\label{fMlip}
\end{align*}

Then, according to the classical result of Pardoux and Peng in \cite{Pardoux:1990ju}, there exists a unique solution $(Y^M,Z^M)\in\mcS^2(\R^{d}) \times \mcH^2(\R^{d\times k})$ of the localized BSDE 
\beq Y^M_t = \xi+\int_t^T f^M\parenth{s,Y^M_s,Z^M_s} \mathrm{d}s - \int_t^T Z^M_s \mathrm{d} W_s, \quad 0 \leqslant t \leqslant T.\label{quadraloc}\eeq

\section{Main results}\label{main}
We have collect in this section principal results proved in our article. All proofs are postponed to sections \ref{section stabilite} and \ref{section applications}. The following subsection gives some existence and uniqueness results while subsection \ref{subsection enonce applications} is dedicated to particular frameworks where these existence and uniqueness results apply.

\subsection{Some general existence and uniqueness results}
\label{subsection enonce existence et unicite gene}

\subsubsection{Existence and uniqueness results when the terminal condition and the generator have bounded Malliavin derivatives}
\label{subsection enonce existence et unicite Malliavin borne}

We consider here a particular framework where the terminal condition and the random part of the generator have bounded Malliavin derivatives. More precisely, let us consider the following assumptions. 

\begin{assumption}[Dxi,b]	
The Malliavin derivative of $\xi$ is bounded:
	$$\Vert \D \xi \Vert_{\mcS^\infty} = \sup_{0 \leqslant t \leqslant T} \Vert \D_t \xi \Vert_{L^\infty}< \infty.$$
\end{assumption}

\begin{assumption}[Df,b] 
\begin{itemize}
	\item[\tbf{(i)}] For all $(t,y,z) \in [0,T] \times \R^d \times \R^{d\times k}$, we have
$$ f(t,y,z) \in \mbbD^{1,2}(\R^d), \quad  \text{ and } \quad \esp{\int_0^T \int_0^T \abs{\D_u f(s,y,z)} \d u \d s}< \infty.$$ 
	\item[\tbf{(ii)}] There exists $C \geqslant 0$ such that for all  $(u,t,y,z) \in [0,T]^2 \times \R^d \times \R^{d \times k}$,
	$$ \abs{\D_u f(t,y,z)} \leqslant C \parenth{1+\abs{z}^2} \ps.$$
	
	\item[\tbf{(iii)}] For all $(u,t) \in [0,T]^2$, there exists a random variable $C_u(t)$ such that for all $(y^1,z^1,y^2,z^2) \in \parenth{\R^d \times \R^{d\times k}}^2$,  
	$$ \abs{\D_u f(t,y^1,z^1)-\D_u f(t,y^2,z^2)} \leqslant C_u(t) \parenth{\parenth{1+\abs{z^1}^2+\abs{z^2}^2}\abs{y^1-y^2}+ \parenth{1+\abs{z^1}+\abs{z^2})} \abs{z^1-z^2}} \quad \ps.$$
\end{itemize}
\end{assumption}
By recalling that $(Y^M, Z^M)$ is the unique solution of \eqref{quadraloc}, we will also assume that we have an \tit{a priori} estimate on $|Z^M| \star W$ uniform in $M$ and small enough. For a given $m>1$ we consider the following assumption:
\begin{assumption}[BMO,m]
there exists a constant $\mbbK <\mbbB^m(L_y,L_z)$ such that 
$$\displaystyle\sup_{M \in \R^+}\norme{\abs{Z^{M}}\star W}_{BMO} \leqslant \mbbK. $$
\end{assumption}

\begin{theo}[Existence and uniqueness (1)]\label{existuniq1}
Let $m >1$. Under the main assumption \ref{H}, the BMO \tit{a priori} estimate \ref{BMO,m},  and the boundedness of the Malliavin derivatives of $\xi$ and $f$, \ref{Dxi,b}---\ref{Df,b}, the quadratic BSDE \eqref{quadra} has a unique solution  $(Y,Z) \in \mcS^2(\R^d) \times \mcZ^m_{BMO}$ such that
$$\supess_{\Omega \times [0,T]} |Z| <+\infty.$$
\end{theo}


A result similar to \cref{existuniq1} can be obtained when the quadratic growth of $z$ has essentially a diagonal structure. Thus, we replace assumption \ref{H} by the following one:
\begin{assumption}[Hdiag]
$ $ 
\begin{itemize}
 \item There exist $f_{\tbf{diag}} : \Omega \times [0,T] \times \R^{d \times k} \rightarrow \R^d$ and $g:\Omega \times [0,T] \times  \R^d \times \R^{d \times k} \rightarrow \R^d$ such that
 for all $i \in \{1,...,d\}$ we have
 $$ f^i(t,y,z)=f^i_{\tbf{diag}}(t,z)+g^i(t,y,z).$$
 \item There exist five nonnegative constants $L_d,K_{d,y},L_{d,y},K_{d,z},L_{d,z}$ such that for all $(t,y,y',z,z') \in [0,T] \times (\R^d)^2 \times (\R^{d \times k})^2$ and $i \in \ensemble{1,...,d}$: 
$$ \abs{f^i_{\tbf{diag}}(t,z)-f^i_{\tbf{diag}}(t,z')} \leqslant L_{d}\parenth{\abs{z^{(i,:)}}+\abs{(z')^{(i,:)}}}\abs{(z-z')^{(i,:)}},$$ 
$$\abs{g(t,y,z)-g(t,y',z)}\leqslant \parenth{ K_{d,y}+L_{d,y}\abs{z}^2}\abs{y-y'},$$
$$ \abs{g(t,y,z)-g(t,y,z')}\leqslant \parenth{K_{d,z}+L_{d,z}\parenth{\abs{z}+\abs{z'}}} \abs{z-z'}.$$
\end{itemize}

\end{assumption}
This kind of framework has been introduced by Hu and Tang in \cite{Hu:2016is} (see also \cite{Jamneshan:2014uia}). The following result of existence and uniqueness is specific, and do not follows directly from \cref{existuniq1}. Indeed, if an uniform upper bound is assumed (assumption \tbf{(i)} below), we can use specific tools in the diagonal case to obtain an upper bound small enough.

\begin{theo}[Existence and uniqueness (1) - Diagonal Case]\label{existuniq1diag}
We assume that \ref{Hdiag}, \ref{Dxi,b}, \ref{Df,b} hold true and that there exists a constant $\mbbB$ such that
\begin{enumerate}[label=(\roman*)]
 \item[\tbf{(i)}]  $\displaystyle\sup_{M \in \R^+}\norme{\abs{Z^{M}}\star W}_{BMO} \leqslant \mbbB, $
 \item[\tbf{(ii)}] $c_2^2 d L_{d,y} \mbbB^2 <1, \quad \parenth{\frac{c_2}{c_1}\sqrt{L_{d,y}} +\frac{2c_2^2\sqrt{d}}{c_1^2}L_{d,z}}\frac{4\sqrt{d} c_2^2L_{d,z} \mbbB^2}{1-c_2^2 d L_{d,y} \mbbB^2} <1,$
 where $c_1$ and $c_2$ are given by \cref{BMOequiv} with $B=2L_d\mbbB$.
\end{enumerate}
We also assume that $\xi \in L^{\infty}(\Omega, \mcF_T)$ and $f(.,0,0) \in \mcS^\infty(\R^d)$. Then, the quadratic BSDE \eqref{quadra} has a unique solution $(Y,Z) \in \mcS^\infty(\R^d) \times BMO_{\mbbB}$ such that
$$\supess_{\Omega \times [0,T]} |Z| <+\infty.$$
\end{theo}

The main difference between assumptions in Theorem \ref{existuniq1} and Theorem \ref{existuniq1diag} comes from the form of constants used in the bound of the BMO norm. In particular, for any $L_d>0$, there exists $\varepsilon>0$ such that \tbf{(ii)} in Theorem \ref{existuniq1diag} is fulfilled as soon as $L_{d,y}<\varepsilon$ and $L_{d,z} <\varepsilon$ while we cannot take $L_z$ as large as we want in Theorem \ref{existuniq1}.


\subsubsection{Extension to general terminal values and generators}

Now we are able to relax assumptions \ref{Dxi,b} and \ref{Df,b} with some density arguments. 
To do so, we assume that we can write $f$ as a deterministic function $\tbf{f}$ of a progressively measurable continuous process: the randomness of the generator will be contained into this process. 

\begin{assumption}[H']
\begin{itemize} 
\item[\tbf{(i)}] There exists a progressively measurable continuous process $\alpha \in \bigcap_{p \in \N^*} \mcS^p$ with values in $\R^{d'}$, $d' \geqslant 1$, and a function $\tbf{f}: \R^{d'} \times \R^d \times \R^{d \times k} \longrightarrow \R^d$ such that for all $(t,y,z) \in [0,T] \times \R^d \times \R^{d \times k}$: 
$$ f(t,y,z)=\tbf{f}(\alpha_t,y,z).$$
Besides, we assume that \ref{H} holds true for {f}. 

\item[\tbf{(ii)}] There exists $D \in \R^+$ and $\delta \in (0,1]$ such that for all $(y,z) \in \R^d \times \R^{d \times k}$, $(\beta,\beta') \in (\R^{d'})^2$: 
\begin{align}\abs{\tbf{f}(\beta,y,z)-\tbf{f}(\beta',y,z)}\leqslant D \parenth{1+\abs{z}^2} \abs{\beta-\beta'}^{\delta}\label{hypbeta}.\end{align}
\end{itemize}

\end{assumption}

For $\eta \in L^2(\Omega, \mcF_T)$, $\beta \in \mcS^\infty$ and $M \in \R^+$, we denote by $\parenth{Y^{(M,\eta,\beta)},Z^{(M,\eta,\beta)}}$ the unique solution of the BSDE
\beq Y^{(M,\eta,\beta)}_t = \eta +\int_t^T \tbf{f}^M\parenth{\beta_s,Y^{(M,\eta,\beta)}_s,Z^{(M,\eta,\beta)}_s} \d s-\int_t^T Z^{(M,\eta,\beta)}_s \d W_s, \quad 0 \leqslant t \leqslant T, \eeq
where for all $(t,y,z) \in [0,T] \times\R^d \times \R^{d \times k}$ and $\R^{d'}$-valued processes $\alpha$, 
we have $\tbf{f}^M(\alpha_t,y,z)=f^M(t,y,z)$. Finally, assumption \ref{BMO,m} will be replaced by the following one.
\begin{assumption}[BMO2,m]
We assume that $\xi \in L^{2m^*}(\Omega, \mcF_T)$ and that there exists a constant $\mbbK<\mbbB^m(L_y,L_z)$ such that
$$\displaystyle\sup_{M \in \R^+}\sup_{\substack{\norme{\eta}_{L^{2m^*}(\Omega,\mcF_T)}\leqslant \norme{\xi}_{L^{2m^*}(\Omega, \mcF_T)}\\\norme{\beta}_{L^2(\Omega\times [0,T])}\leqslant \norme{\alpha}_{L^2(\Omega\times [0,T])}}}\norme{\abs{Z^{(M,\eta,\beta)}}\star W}_{BMO} \leqslant \mbbK. $$
\end{assumption}

\begin{theo}[Existence and uniqueness (2)]\label{withoutxib}
Let $m>1$. Under the main assumption \ref{H'} and the BMO estimation \ref{BMO2,m}, the quadratic BSDE \eqref{quadra} has a unique solution in $\mcS^{2m^*}(\R^d) \times (\mcH^{m^*}(\R^{d\times k}) \cap \mcZ^m_{BMO})$.
\end{theo}

\begin{rema}
$ $
\begin{itemize}
 \item Let us emphasize that the uniqueness result in Theorem \ref{withoutxib} lies in a different space than the space used in Theorem \ref{existuniq1}.
 \item It is also possible to extend the result of \cref{existuniq1diag} (diagonal case) to more general terminal conditions and generators. Nevertheless, the result obtained would be less general than \cref{withoutxib}. See \cref{extension existence unicite cas diagonal} for more details.
\end{itemize}
\end{rema}

\subsection{Applications to multidimensional quadratic BSDEs with special structures}
\label{subsection enonce applications}
In this subsection we give some explicit frameworks where assumptions \ref{BMO,m} and \ref{BMO2,m} or assumptions \tbf{(i)} and \tbf{(ii)} of \cref{existuniq1diag} are fulfilled. The aim is to show that numerous results on multidimensional quadratic BSDEs already proved in the literature can be obtained with similar assumptions by our approach. We want to underline the simplicity of this approach since we just have to obtain some \tit{a priori} estimates on the BMO norm of $\abs{Z} \star W$ by using classical tools as explained in section \ref{section applications}. Moreover, it is quite easy to construct some << new >> frameworks where \ref{BMO,m} and \ref{BMO2,m} or assumptions \tbf{(i)} and \tbf{(ii)} of \cref{existuniq1diag} are also fulfilled.

\subsubsection{An existence and uniqueness result for BSDEs with a small terminal condition}

In \cite{Tevzadze:2008co}, Tevzadze obtains a result of existence and uniqueness for multidimensional quadratic BSDEs when the terminal condition is small enough by using a contraction argument in $\mcS^\infty \times BMO$. We are able to deal with this kind of assumption with our approach. We consider the following hypothesis.
\begin{assumption}[HQ]
\begin{itemize}
	\item[\tbf{(i)}] There exists $\gamma \in \R^+$ such that for all $(t,y,z) \in [0,T] \times \R^d \times \R^{d \times k}$, we have $\abs{f(t,y,z)}\leqslant \gamma \abs{z}^2,$
	\item[\tbf{(ii)}] $32 \gamma^2 \norme{\xi}_{L^\infty}^2 \leqslant 1.$
\end{itemize} 
\end{assumption}
\begin{prop}
\label{prop Tevzadze revisite}
	Let $m >1$. Under \ref{H'}---\ref{HQ}, and the following condition on $\gamma$:
	$$ \frac{1}{2\sqrt{2}\gamma}\parenth{1-\sqrt{1-32\gamma^2 \norme{\xi}_{L^\infty}^2}}^{1/2}<\mbbB^m(L_y,L_z),$$
	the BSDE \eqref{quadra} has a unique solution in $\mcS^{2m^*}(\R^d) \times (\mcH^{m^*}(\R^{d\times k}) \cap \mcZ^m_{BMO})$.  
	If in addition \ref{Dxi,b} and \ref{Df,b} hold true, there exists an unique solution $(Y,Z) \in \mcS^\infty(\R^d) \times \mcZ^m_{BMO}$ such that
$$\supess_{\Omega \times [0,T]} |Z| <+\infty.$$
\end{prop}

\subsubsection{An existence and uniqueness result for BSDEs with a monotone generator} 
\label{subsection monotone}
In this part we investigate the case where we have for $f$ a kind of monotonicity assumption with respect to $y$. 

\begin{assumption}[HMon]
	\begin{itemize}
		\item[\tbf{(i)}] There exists $\mu >0$ and $\alpha,\gamma \geqslant 0$ such that for all $(s,y,z) \in [0,T] \times \R^d \times \R^{d \times k}$
$$y . f(s,y,z) \leqslant \alpha\abs{y}-\mu \abs{y}^2+\gamma \abs{y} \abs{z}^2,$$
	\item[\tbf{(ii)}] $32 \gamma^2 A^2 \leqslant 1$, where $A=\max \parenth{\norme{\xi}_{L^\infty}, \frac{\alpha}{\mu}}$. 
	\end{itemize} 
\end{assumption}

\begin{prop}
\label{prop Tevzadze monotone}
	Let $m >1$. Under \ref{H'}---\ref{HMon} and the following estimate on $\gamma$:
	$$ \frac{1}{2\sqrt{2}\gamma}\parenth{1-\sqrt{1-32\gamma^2 A^2}}^{1/2}<\mbbB^m(L_y,L_z),$$
	the quadratic BSDE \eqref{quadra} has a solution in $\mcS^{2m^*}(\R^d) \times (\mcH^{m^*}(\R^{d \times k}) \cap \mcZ^m_{BMO})$. If in addition \ref{Dxi,b} and \ref{Df,b} hold true, there exists a unique solution $(Y,Z) \in \mcS^\infty(\R^d) \times  \mcZ^m_{BMO}$ such that
$$\supess_{\Omega \times [0,T]} |Z| <+\infty.$$ 
\end{prop}

\subsubsection{An existence and uniqueness result for diagonal quadratic BSDEs} 

Now we consider the diagonal framework introduced in section \ref{subsection enonce existence et unicite Malliavin borne}. We assume that the generator satisfies \ref{Hdiag}, \tit{i.e.} the generator $f$ can be written as $f(t,y,z) =  f_{\tbf{diag}}(t,z)+g(t,y,z)$ where $f_{\tbf{diag}}$ has a diagonal structure with respect to $z$.

\begin{prop}
\label{proposition existence unicite cas diagonal}
 We assume that
 \begin{enumerate}[label=(\roman*)]
  \item[\tbf{(i)}]  \ref{Hdiag}, \ref{Dxi,b} and \ref{Df,b} hold true,
  \item[\tbf{(ii)}]  there exist nonnegative constants $G_d$ and $G$ such that, for all $(t,y,z)\in [0,T] \times \R^d \times \R^{d \times k}$, we have
  \begin{equation} \label{croissance structure diagonale} \abs{f_{\tbf{diag}}(t,z)} \leqslant G_d \abs{z}^2, \quad \abs{g(t,y,z)} \leqslant G \abs{z}^2,\end{equation}
  \item[\tbf{(iii)}]  $$ \frac{4\sum_{i=1}^d e^{2G_d \norme{\xi^i}_{L^{\infty}}}}{G_d}G \leqslant 1,$$
  \item[\tbf{(iv)}]    $c_2^2 d L_{d,y} (4G_dG)^{-1} <1, \quad \parenth{\frac{c_2}{c_1}\sqrt{L_{d,y}} +\frac{2c_2^2\sqrt{d}}{c_1^2}L_{d,z}}\frac{4 \sqrt{d} c_2^2L_{d,z} (4G_dG)^{-1}}{1-c_2^2 d L_{d,y} (4G_dG)^{-1}} <1,$
 where $c_1$ and $c_2$ are given by \cref{BMOequiv} with $B=2L_d(4G_dG)^{-1/2}$. 
 \end{enumerate}
 Then, the quadratic BSDE \eqref{quadra} has a unique solution  $(Y,Z) \in \mcS^\infty(\R^d) \times BMO_{(4G_dG)^{-1/2}}$ such that
$$\supess_{\Omega \times [0,T]} |Z| <+\infty.$$
\end{prop}

\begin{rema}
\label{remarque autre hypothese possibles pour cas diag}
 The growing assumption \eqref{croissance structure diagonale} is only one example of hypothesis that can be tackled by our approach. It is also possible to obtain the same kind of result by replacing \eqref{croissance structure diagonale} by one of the following assumption:
 \begin{itemize}
  \item We assume that for all $(t,y,z)\in [0,T] \times \R^d \times \R^{d \times k}$,
 $$\abs{g(t,y,z)} \leqslant C(1+\abs{y}) + \varepsilon\abs{z}^2$$
 and $T,\varepsilon$ are supposed to be small enough. This framework is studied in \cite{Hu:2016is,Jamneshan:2014uia}.
  \item We assume that for all $(t,y,z)\in [0,T] \times \R^d \times \R^{d \times k}$,
 $$\abs{g(t,y,z)} \leqslant C(1+\abs{y}).$$
  This situation is already studied in \cite{Hu:2016is}.
\end{itemize}
\end{rema}

\subsubsection{Existence and uniqueness of martingales in manifolds with prescribed terminal condition} 

The problem of finding martingales on a manifold with prescribed terminal value has generated a huge amount of literature. On the one hand with geometrical methods, Kendall in \cite{Kendall:1990bp} treats the case where the terminal value lies in a geodesic ball and is expressed as a functional of the Brownian motion. Kendall gives also a characterisation of the uniqueness in terms of existence of a convex separative function, i.e. a convex function on the product space which vanishes exactly on the diagonal. Besides, in \cite{Kendall:1992eu}, Kendall proved that the property \tit{every couple of points are connected by a unique geodesic} is not sufficient to ensure existence of a separative convex function, which was conjectured by Émery. An approach by barycenters, of the martingale notion on a manifold, is used by Picard in \cite{PicardJean:1994wv} for Brownian filtrations. Arnaudon in \cite{Arnaudon:1997tm} solved the problem in a complex analytic manifold having a convex geometry property for 
continuous filtrations: the main idea is to consider a differentiable family of martingales. For all these results, a convex geometry property is assumed. The first approach using the tool of BSDEs is proposed by Darling in \cite{Darling:1995wd}. 

Let us now define more precisely the problem. A so-called linear connection structure is required to define martingales on a manifold $\mcM$ in a intrinsic way. \textit{A contrario}, for semimartingales, a differential structure is enough. The definition of a martingale can be rewritten with a system of coupled BSDEs having a quadratic growth, so we begin to recall it. We can refer to \cite{Emery:1989cy} for more details about stochastic calculus on manifolds. 
	
Let us consider $(\mcM, \nabla)$ a differential manifold equipped with a linear connection $\nabla$. This is equivalent to give ourselves a Hessian notion or a covariant derivative. We say that a continuous process $X$ is a semimartingale on $M$ if for all $F \in \mcC^2(\mcM)$, $F \circ X$ is a real semimartingale. Consistence of the definition is simply due to the Itô formula. We say that a continuous process $Y$ is a (local) $\nabla$-martingale if for all $F \in \mcC^2(\mcM)$, 
$$ F(Y)_t-\frac{1}{2} \int_0^t \nabla \d F (\d Y, \d Y)_s$$
is a real local martingale on $[0,T]$. Again it is not very hard to see with the Itô formula that this definition is equivalent to the Euclidean  one in the flat case. Let us remember that $\dint_0^{.} \nabla \d F (\d Y, \d Y)_s$ is a notation for the quadratic variation of $Y$ with respect to the $(0,2)$-tensor field $ \nabla \d F$. This notion is defined by considering a proper embedding 
	$(x_i)_{1 \leqslant i \leqslant d}$ into $\R^d$ such that every bilinear form $b$ can be written as $b =b_{ij} \d x^i \otimes \d x^j$ (implicit summation). On the other hand it can be proved that the quantity 
	$$ \int_0^. b (\d Y,\d Y)_s:=\int_0^. b_{ij}(Y_s)\d \scalaire{Y^i}{Y^j}_s$$
	does not depend on $(x^i)_{1 \leqslant i \leqslant d}$ and so the quantity $\dint_0^{.} \nabla \d F (\d Y, \d Y)_s$ is intrinsic.  It is well-known that for all $m \in \mcM$,
	$$ \parenth{\nabla \d F}_{ij}(m)=D_{ij}F(m)-\Gamma^{k}_{ij}(m) D_{k}F(m),$$
	where $\Gamma^{k}_{ij}(m)$ denotes a $\nabla$-Christoffel symbol at the point $m$. The coefficients are symmetric with respect to $i,j$. 
	 Hence martingale property in the domain of a local chart is equivalent to the existence of a process $Z$ such that $(Y,Z)$ solves the following BSDE
		$$ Y_t=\xi + \int_t^T f(s,Y_s,Z_s) \d s-\int_t^T Z_s \d W_s, \quad 0 \leqslant t \leqslant T,$$
	with $f : [0,T] \times \R^d \times \R^{d \times k}\rightarrow \R^d$ defined by $f(s,y,z)=\frac{1}{2}\parenth{\Gamma^{k}_{ij}(y) z^{(i,:)}. z^{(j,:)}}_{1 \leqslant k \leqslant d}$. 
	It is an easy consequence of the representation theorem for Brownian martingales and the definition applied to $F=x^i$. 	We consider in addition the following assumption
	\begin{assumption}[HGam]
		there exists two constants $L_y$ and $L_z$ such that for all $i,j,k \in \{1,...,d\}$  
$$\abs{\Gamma^{k}_{ij}(y)-\Gamma^{k}_{ij}(y')}\leqslant 2 L_y\abs{y-y'}, \quad \abs{\Gamma^{k}_{ij}(y)} \leqslant 2L_z. $$
	\end{assumption}
	
		For example \ref{HGam} is in force if the domain of the chart is a compact set. It is also true if we choose an exponential chart. Without loss of generality we can suppose that $\mcM$ has a global system of coordinates: all the Christoffel symbols will be computed in this system. 	
	
	Under \ref{HGam}, assumption \ref{H} is in force: for all $(y,y',z,z') \in (\R^d)^2 \times \parenth{\R^{d \times k}}^2$, 	
 $$ \abs{f(t,y,z)-f(t,y',z)}\leqslant L_y \abs{z}^2 \abs{y-y'},$$
	and with the symmetric property of the Christoffel symbols, we have
	$$ f(t,y,z')-f(t,y,z)=-\frac{1}{2} \sum_{i,j}\Gamma^{\bullet}_{ij}(y)\parenth{z^{(i,:)}. z^{(j,:)}-(z')^{(i,:)}.(z')^{(j,:)}}=-\frac{1}{2} \sum_{i,j}\Gamma^{\bullet}_{ij}(y)\parenth{z^{(i,:)}-(z')^{(i,:)}}\parenth{z^{(j,:)}+(z')^{(j,:)}},$$
	which implies that
	$$ \abs{f(t,y,z)-f(t,y,z')} \leqslant L_z \parenth{\abs{z}+\abs{z'}}\abs{z-z'}.$$
	
To obtain some important \tit{a priori} estimate for the BMO norm of $Z \star W$, Darling introduce in \cite{Darling:1995wd} a convex geometry assumption.

\begin{defi}
	We say that a function $F \in \mcC^2 (\mcM,\R)$ (seen as a function on $\R^d$) is doubly convex on a set $G \subset \R^d$ if for all $y \in G$ and $z \in \R^d$, 
	$$ \min\ensemble{\Hess F(y)(z,z),\nabla\d F(y)(z,z)} \geqslant 0,$$
	and, for $\alpha >0$, $F$ is $\alpha$-strictly doubly convex on $G$ if for all $y \in G$ and $z \in \R^d$, 
	$$ \min\ensemble{\Hess F(y)(z,z),\nabla\d F(y)(z,z)} \geqslant \alpha \abs{z}^2.$$
\end{defi}
This property means that $F$ is convex with respect to the flat connection, and, with respect to the connection $\nabla$. 

\begin{theo}
\label{theo existence unicite martingale variete}
	Let $m >1$ and assume that:
	\begin{itemize}
		\item[\tbf{(i)}] there exists a function $F^{dc} \in \mcC^2 (\mcM,\R)$, such that $G=\parenth{F^{dc}}^{-1}(]-\infty,0])$ is compact and $\xi \in G$, 		
		\item[\tbf{(ii)}] $F^{dc}$ is doubly convex on $\mcM$, and there exists $\alpha >0$ and $m \geqslant 1$ such that $F^{dc}$ is $\alpha$-strictly doubly convex on $G$ and satisfies
		$$ \parenth{\sup_{(x,y) \in G^2} \ensemble{F^{dc} (x)-F^{dc} (y)}}^{1/2} \leqslant \sqrt{\frac{\alpha}{2}}\times\mbbB^m(L_y,L_z),$$				
		\item[\tbf{(iii)}] \ref{HGam} holds true.
	\end{itemize}
	Then there exists a unique $\nabla$-martingale $Y$ in $\mcS^\infty(\R^d)$ with terminal value $\xi$ such that $\sqrt{|\scalaire{Y}{Y}|} \star W\in  \mcZ^m_{BMO}$. Moreover, if \ref{Dxi,b} holds true we have
$$\supess_{\Omega \times [0,T]} |\scalaire{Y}{Y}| <+\infty.$$ 
\end{theo}

\begin{rema}
 By using the same approach, it should be possible to extend the previous result to $\nabla$-Christoffel symbols that depend on time or even that are progressively measurable random processes. 
\end{rema}

\subsubsection{The Markovian setting} 
The aim of this subsection is to refine some results of Xing and {\v Z}itkovi{\'c} obtained in \cite{2016arXiv160300217X}: in this paper, authors establish existence and uniqueness results for a general class of Markovian multidimensional quadratic BSDEs. Let us start by introducing the Markovian framework. For all $t\in[0,T]$ and $x \in {\R}^{k}$ we denote $X^{t,x}$ a diffusion process satisfying the following SDE
\beq \begin{cases}
      \d X^{t,x}_s = b(s,X^{t,x}_s) \d t + \sigma(s,X^{t,x}_s) \d W_s, \quad s \in [t,T],\\
      X^{t,x}_s = x, \quad s \in [0,t]. 
     \end{cases}
\label{sde}\eeq
In all this part, we assume following assumptions that ensure, in particular, that for all $(t,x) \in [0,T] \times {\R}^{k}$, there exists a unique strong solution of \eqref{sde}.

\begin{assumption}[HX]
 \begin{itemize}
  \item The drift vector $b : [0,T] \times {\R}^{k} \rightarrow {\R}^{k}$ is measurable and uniformly bounded,
  \item The dispersion matrix $\sigma : [0,T] \times {\R}^{k} \rightarrow {\R}^{k \times k}$ is symmetric, measurable and there exists a constant $\Lambda>0$ such that 
  $ \Lambda \abs{u}^2 \geqslant \abs{\sigma(t,x)u}^2 \geqslant \frac{1}{\Lambda} \abs{u}^2$ for all $(t,x) \in [0,T] \times \R^{k}$ and $u \in \R^{k}$,
  \item $b$ and $\sigma$ are Lipschitz functions with respect to $x$.
 \end{itemize}
\end{assumption}

The aim of this subsection is to study the following Markovian BSDE
\begin{equation} 
Y^{t,x}_{u} = \mcG (X^{t,x}_T)+\int_u^T \tbf{f}\parenth{s,X_s^{t,x},Y_s^{t,x},Z_s^{t,x}} \mathrm{d}s - \int_u^T Z_s^{t,x} \mathrm{d} W_s, \quad t \leqslant u \leqslant T, \label{XZquadra} 
\end{equation}
for which we assume following assumptions:
\begin{assumption}[HMark]
 \begin{itemize}
  \item $(s,y,z) \in [0,T] \times \R^d \times \R^{d \times k} \longmapsto \tbf{f}(s,X_s^{t,x},y,z)$ satisfies \ref{H},
  \item There exists $D \in \R^+$ and $\kappa \in (0,1]$ such that for all $(x,x') \in (\R^{k})^2$,
   $$\abs{\mcG(x)} \leqslant D, \quad \abs{\mcG(x)-\mcG(x')} \leqslant D\abs{x-x'}^\kappa.$$
 \end{itemize}
\end{assumption}
As in \cite{2016arXiv160300217X} we say that a pair $(v,w)$ of functions is a continuous Markovian solution of \eqref{XZquadra} if
\begin{itemize}
 \item $v : [0,T] \times \R^k \rightarrow \R^d$ is a continuous function and $w : [0,T] \times \R^k \rightarrow \R^{d\times k}$ is a measurable function,
 \item for all $(t,x) \in [0,T] \times \R^k$, $(Y^{t,x},Z^{t,x}):=(v(.,X^{t,x}),w(.,X^{t,x}))$ is a solution of \eqref{XZquadra}.
\end{itemize}
Two Markovian solutions, $(v,w)$ and $(v',w')$, are considered equal if $v(t,x) = v'(t,x)$ for all $(t,x) \in [0,T] \times \R^k$ and $w=w'$ a.s. with respect to the Lebesgue measure on $[0,T] \times \R^k$.


Some existence and uniqueness results about continuous Markovian solutions of \eqref{XZquadra} are obtained in \cite{2016arXiv160300217X} by assuming the existence of a so-called \tit{Lyapunov function}. We recall here the definition of these functions given in \cite{2016arXiv160300217X}.

\begin{defi}[Lyapunov function associated to $\normalfont{\tbf{g}}$]
	Let $\tbf{g} : [0,T] \times \R^k \times \R^d \times \R^{d \times k} \rightarrow \R^d$ a Borel function. A nonnegative function $F \in \mcC^2(\R^d, \R)$  is said to be a Lyapunov function associated to $\tbf{g}$ if for all $(t,x,y,z) \in [0,T] \times \R^{k} \times \R^d \times \R^{d \times k}$:
	$$ \frac12 \sum_{l=1}^d \Hess F(y) \parenth{z^{(:,l)},z^{(:,l)}}-\d F(y) \tbf{g} \normalfont{\tbf{}}(t,x,y,z) \geqslant \abs{z}^2.$$
\end{defi}

We are now able to give a uniqueness result that partially refine the result given by \cite{2016arXiv160300217X}.

\begin{theo}[Uniqueness for the Markovian case]
\label{th unicite cas markovien}
	We assume that 
\begin{itemize}
		\item[\tbf{(i)}]\ref{HX} and \ref{HMark} are in force.\item[\tbf{(ii)}] there exists a Lyapunov function $F$ associated to $\tbf{f}$. 
\end{itemize}		
Then \eqref{XZquadra} admits at most one continuous Markovian solution $(v,w)$ such that $v$ is bounded.
\end{theo}

Moreover, we are also able to precise the regularity of the solution when it exists.

\begin{theo}[Regularity of the Markovian solution]
\label{th regularite markov}
We assume that:
	\begin{itemize}
		\item[\tbf{(i)}] \ref{HX} and \ref{HMark} are in force,
		\item[\tbf{(ii)}] there exists $D \in \R^+$ and $\kappa \in (0,1]$ (same constant $\kappa$ as in \ref{HMark}) such that for all $(s,x,x',y,z)\in  [0,T] \times(\R^{k})^2\times \R^d \times \R^{d \times k}$,
$$\abs{\tbf{f}(s,x,y,z)} \leqslant D(1+\abs{y}+\abs{z}^2),\quad \abs{\tbf{f}(s,x,y,z)-\tbf{f}(s,x',y,z)}\leqslant D(1+|z|^2) \abs{x-x'}^{\kappa},$$
 \item[\tbf{(iii)}] there exists a Lyapunov function $F$ associated to $\tbf{f}$. 
 \end{itemize}
 
If $(v,w)$ is a continuous Markovian solution of \eqref{XZquadra} such that $v$ is bounded, then $v \in \mcC^{\kappa}$. 
Particularly, if $\kappa =1$ then $w$ is essentially bounded: the multidimensional quadratic BSDE \eqref{XZquadra} becomes a standard multidimensional Lipschitz BSDE by a localisation argument.
\end{theo}

\begin{rema}
 An existence result is given by Theorem 2.7 in \cite{2016arXiv160300217X}. A less general existence result can be obtained thanks to our approach by combining estimates obtained by Xing and {\v Z}itkovi{\'c} in Theorem 2.5 of \cite{2016arXiv160300217X}, small BMO estimates obtained in the proof of \cref{th unicite cas markovien} and \cref{remexistence-uniqueness2-avecepsBMO} but the approach is less direct than in \cite{2016arXiv160300217X}. Concerning the uniqueness, Xing and {\v Z}itkovi{\'c} have proved a uniqueness result for generators that do not depend on $y$: our result allows to fill this small gap. Finally,  Xing and {\v Z}itkovi{\'c} prove that there exists a Markovian solution that satisfies $v \in \mcC^{\kappa',loc}$ with  $\kappa' \in (0,\kappa]$. Thus, our regularity result gives a better estimation of the solution regularity since the regularity of the terminal condition and the generator is retained. In particular, we obtain that $Z$ is bounded when $\kappa=1$ which can have important 
applications, as pointed out in the introduction.
\end{rema}

\begin{rema}
The existence of a Lyapunov function seems to be an \textit{ad hoc} theoretical assumption at first sight but Xing and {\v Z}itkovi{\'c} provide in \cite{2016arXiv160300217X} a lot of examples and concrete criteria to obtain such kind of functions. Moreover we can note that the Lyapunov function can be used to obtain \tit{a priori} estimates on $\norme{\abs{Z} \star W}_{BMO}$ (see the proof of \cref{th unicite cas markovien} and \cref{th regularite markov}). 
\end{rema}

\section{Generalities about SDEs and linear BSDEs}\label{tools}

We collect in this section some technical results that will be useful for section \ref{section stabilite} and section \ref{section applications}.

\subsection{The linear case: representation of the solutions} 

We investigate here the following linear BSDE
\beq
U_t = \zeta+\int_t^T \left(A_s U_s +B_s V_s+f_s\right) \mathrm{d}s - \int_t^T V_s \mathrm{d} W_s, \quad 0 \leqslant t \leqslant T \label{BSDElingen}, \eeq 
where $\zeta \in L^2(\mcF_T, \R^d), f \in L^2(\Omega \times [0,T])$ and  $A,B,f$ are three processes with values in  $\mathcal{L}(\R^{d}, \R^{d})$, $\mcL \left(\R^{d\times k}, \R^d \right)$ and $\R^{d}$. For the linear case we have an explicit formulation of the solution. Let us begin to recall the classical scalar formula which can be obtained using the Girsanov transform. 

\begin{rema}[One-dimensional case ($d=1$)] It is well-known that the solution of \eqref{BSDElingen} is given by the formula 
$$ \quad U_t = \espcond{ S_t^{-1} S_T  \zeta+\int_{t}^T S_t^{-1} S_s  f_s \mathrm{d}s }{\mathcal{F}_t}, \quad 0 \leqslant t \leqslant T, $$
where
$$S_{t}=\exp \left( \int_0^t B_s \mathrm{d} W_s - \frac{1}{2} \int_0^t \vert B_s \vert^2 \mathrm{d}s + \int_0^t A_s \mathrm{d}s \right)=\mathcal{E}\left( B \star W \right)_t \exp\left(\int_0^t A_s \mathrm{d}s \right).$$
\end{rema} 
To extend this last formula in the general case we define, as in \cite{Delbaen:2008bh}, a process $S$ as the unique strong solution of  
$$ \mathrm{d} S_t = \sum_{p=1}^k S_t B_t^{(:,p,:)} \mathrm{d} W^p_t + S_t A_t\mathrm{d}t, \quad S_0=I_{d \times d}.$$

\begin{prop}[Formula for $U$]\label{expressiongeneralesolulin}
$ $
\begin{itemize}
	\item[\tbf{(i)}] The process $S$ is almost surely invertible for all $t \in [0,T]$ and $S^{-1}$ is the solution of 
	$$ \d S_t^{-1}=\left[ \parenth{ \sum_{p=1}^k \parenth{B_t^{(:,p,:)}}^2 -A_t} \d t - \sum_{p=1}^k B_t^{(:,p,:)} \d W_t^p \right] S_t^{-1}, \quad S_0^{-1}=I_{d \times d}.$$
	\item[\tbf{(ii)}] The BSDE \eqref{BSDElingen} has a unique solution $(U,V)$ in $\mathcal{S}^2 \left( \R^d \right) \times \mathcal{H}^2 \left( \R^{d\times k} \right)$, and $U$ is given by: 
\beq U_t=\espcond{ S_t^{-1} S_T  \zeta+\int_{t}^T S_t^{-1} S_s  f_s \mathrm{d}s }{\mathcal{F}_t}.\label{formula}\eeq
\end{itemize} 
\end{prop}

\begin{Proof}	
	Existence and uniqueness of a solution $(U,V)$ in $\mcS^2(\R^{d}) \times \mcH^2(\R^{d\times k})$ is guaranteed by the Pardoux and Peng result in \cite{Pardoux:1990ju}.
	The solution $(U,V)$ satisfies 
	$$ U_t=\zeta+\int_t^T \parenth{ A_s U_s+\sum_{p=1}^k B_s^{(:,p,:)} V_s^{(:,p)} +f_s}\d s-\sum_{p=1}^k \int_t^T V_s^{(:,p)} \d W_s^p.$$
	The Itô formula gives the invertibility of $S$ and the formula for $S^{-1}$ on the one hand. On the other hand:
	$$ \d \parenth{S_t U_t}=-S_t f_t \d t+\sum_{p=1}^k \parenth{S_t B_t^{(:,p,:)}U_t+S_t V_t^{(:,p)}}\d W_t^p,$$
	and thus we get, for all $t \in [0,T]$,
   $$ S_t U_t=S_T \zeta + \int_t^T S_s f_s \d s- \int_t^T \sum_{p=1}^k \parenth{S_s B_s^{(:,p,:)}U_s+S_s V_s^{(:,p)}}\d W_s^p.$$
	By taking the conditional expectation $ S_t U_t = \espcond{S_T \zeta +\displaystyle\int_t^T S_s f_s \d s }{\mcF_t}.$ Adaptability and invertibility of $S$ give the result. 
\end{Proof}

\subsection{A result about SDEs}
We consider a SDE on $\R^{d\times d}$ of the form 
\begin{equation}\label{SDE}
	X_t = X_0 + \int_0^t F(s,X_s) \mathrm{d}s + \sum_{p=1}^k \int_0^t G^p(s,X_s) \mathrm{d}W^p_s,
\end{equation}
where $ F : \Omega \times [0,T] \times \R^{d\times d} \rightarrow \R^{d \times d}$ and for all $p \in \ensemble{1,...,k}$, $G^p : \Omega \times [0,T] \times \R^{d \times d} \rightarrow \R^{d \times d}$ are progressively measurable functions. We start by recalling a result of Delbaen and Tang (see \cite{Delbaen:2008bh}, Theorem 2.1) about existence and uniqueness of a solution to the equation \eqref{SDE}, under BMO assumptions. 

\begin{prop}\label{SDEprop}
	Let $m \geqslant 1$. We suppose that there are two non-negative adapted processes $\alpha$ and $\beta$ such that 
	\begin{enumerate}
		\item[\textbf{(i)}] (Regularity) $F(t,0)=0$, $G(t,0)=0$ and for all $(x_1,x_2,t)\in (\R^{d\times d})^2 \times [0,T]$, $$\abs{F(t,x_1)-F(t,x_2)}\leqslant \alpha_t \abs{x_1-x_2} \quad \ps,$$ 
$$\sum_{p=1}^k \abs{G^p(t,x_1)-G^p(t,x_2)}^2\leqslant \beta_t^2 \abs{x_1-x_2}^2 \quad\ps.$$
 
		\item[\textbf{(ii)}](Sliceability) $(\sqrt{\alpha} \star W$, $\beta \star W) \in BMO_{\e_1}\times BMO_{\e_2}$ with the condition 
	$$2m\varepsilon_1^2+\sqrt{2}\varepsilon_2 C'_m < 1.$$
	\end{enumerate}
	Then there exists a solution $X \in \mcS^m(\R^d)$ to the equation \eqref{SDE} and a constant $K_{m,\varepsilon_1,\varepsilon_2}$ such that
	$$ \norme{X}_{\mcS^m} \leqslant K_{m,\varepsilon_1,\varepsilon_2} \norme{X_0}_{L^m}.$$
\end{prop}

 For the reader convenience a proof of this result can be found in the appendix. From this last proposition we can deduce the following corollary (see \cite{Delbaen:2008bh}, Corollary 2.1)

\begin{coro}\label{ReverseSDEprop}
	Let $m \geqslant 1$. We suppose that there are two non-negative adapted processes $\alpha$ and $\beta$ such that 
	\begin{enumerate}
		\item[\textbf{(i)}] (Regularity) $F(t,0)=0$, $G(t,0)=0$ and for all $(x_1,x_2,t)\in (\R^{d\times d})^2 \times [0,T]$, $$\abs{F(t,x_1)-F(t,x_2)}\leqslant \alpha_t \abs{x_1-x_2} \ps,$$ 
$$\sum_{p=1}^k \abs{G^p(t,x_1)-G^p(t,x_2)}^2\leqslant \beta_t^2 \abs{x_1-x_2}^2 \ps.$$
 
		\item[\textbf{(ii)}](Sliceability) $(\sqrt{\alpha} \star W$, $\beta \star W) \in BMO_{\e_1}\times BMO_{\e_2}$ with the condition 
	$$2m\varepsilon_1^2+\sqrt{2}\varepsilon_2 C'_m < 1.$$
	\end{enumerate}
	For $t \in [0,T]$, let $X^{t,I_d}$ the unique solution defined on $[t,T]$ of the SDE \eqref{SDE} such that $X^{t,I_d}_t=I_d$. Then $X^{t,I_d}$ is in $\mcS^m(\R^d)$ and satisfies for a constant $K_{m}$ depending only on $C'_m$, $m$, $k$ and $\varepsilon_1,\varepsilon_2$:
	\beq \espcond{\sup_{t \leqslant s \leqslant T} \abs{X^{t,I_d}_s}^m}{\mathcal{F}_t} \leqslant K_{m,\varepsilon_1,\varepsilon_2}^m.\label{reverse Holder}\eeq
	In particular, if $X$ is an invertible solution to the equation \eqref{SDE} and if $F$ and $G$ are linear with respect to $x$, we get the reverse Hölder inequality	
	$$ \espcond{\sup_{t \leqslant s \leqslant T} \abs{X_t^{-1} X_s}^m}{\mathcal{F}_t} \leqslant K_{m,\varepsilon_1,\varepsilon_2}^m.$$
\end{coro}

\begin{Proof}
	We can use \cref{SDEprop}. For all $t \in [0,T]$ and all event $\mbbA \in \mathcal{F}_t$,
	$$ \norme{X^{t,I_d} \times \indicat_\mbbA}_{\mcS^m([t,T])} \leqslant K_{m,\varepsilon_1,\varepsilon_2} \norme{I_d \times \indicat_\mbbA}_{L^m}.$$
	Then we get, for all $t \in [0,T]$,
	$$ \esp{\sup_{t \leqslant s \leqslant T} \abs{X^{t,I_d}_s \times \indicat_{\mbbA}}^m} \leqslant K_{m,\varepsilon_1,\varepsilon_2}^m \esp{\abs{\indicat_{\mbbA}}^m},$$
	and we have
	$$ \esp{\sup_{t \leqslant s \leqslant T} \abs{X^{t,I_d}_s}^m \times \indicat_{\mbbA}} \leqslant K_{m,\varepsilon_1,\varepsilon_2}^m \esp{\indicat_{\mbbA}}.$$
	Finally, the definition of conditional expectation gives us the result. If $X$ is invertible and if $F$ and $G$ are linear with respect to $x$, the process $X_t^{-1}X$ is for all $t$ a solution taking the value $I_d$ at $s=t$. The particular case is shown using \eqref{reverse Holder}. 
\end{Proof}

\begin{rema}
\label{lem question ouverte}
 The main limitation of \cref{ReverseSDEprop} comes from assumption \tbf{(ii)}: we need to have a small BMO norm estimate on processes $(\sqrt{\alpha} \star W$, $\beta \star W)$ to get a reverse Hölder inequality. It is well known that we have a more general result when $d=1$: if $\alpha=0$ and $\beta \star W \in BMO$ then there exists $m>1$ (that depends on the BMO norm of $\beta \star W \in BMO$) such that $X$ satisfies a reverse Hölder inequality with the exponent $m$ (see Theorem 3.1 in \cite{Kazamaki:2007tr} and references inside, or \cite{Chikvinidze:2013gj} for a new recent proof). We do not know if this result stays true in the multidimensional framework but we emphasize that this is a crucial open question. Indeed, if such a result is true, then we whould be able to prove that \cref{existuniq1} and \cref{withoutxib} stay true without assuming $\mbbK <\mbbB^m(L_y,L_z)$ in hypothesis \ref{BMO,m} (at least when the generator is Lipschitz with respect to $y$, i.e. $L_y=0$).
\end{rema}

\subsection{Estimates for the solution to BSDE \eqref{BSDElingen}} 
We come back to the linear BSDE \eqref{BSDElingen}, and we want to obtain some $\mcS^q$-estimations for $U$ with $q$ large enough, including $q=\infty$, under BMO assumptions.

\begin{prop}\label{condfinalegene}
Let $m \geqslant 1$. We assume that $B$ and $A$ are adapted, bounded respectively by two non negative processes $\beta$ and $\alpha$ such that: $(\sqrt{\alpha} \star W$, $\beta \star W) \in BMO_{\e_1}\times BMO_{\e_2}$ with the condition 
	$$2m\varepsilon_1^2+\sqrt{2}\varepsilon_2 C'_m < 1.$$
Then
\begin{itemize} 
\item[\textbf{(i)}] If $\zeta \in L^\infty (\Omega, \mcF_T)$ and $f \in \mcS^\infty$, then $U \in \mcS^\infty(\R^d)$ and 
$$ \norme{U}_{\mathcal{S}^\infty} \leqslant K_{m,\varepsilon_1,\varepsilon_2} \parenth{\norme{\zeta}_{L^\infty}+T \norme{f}_{\mcS^\infty}},$$
\item[\textbf{(ii)}] Let us assume that $m>1$. If $\zeta \in L^\infty(\Omega, \mcF_T)$, and $\sqrt{\abs{f}}\star W\in BMO$, then $U \in \mcS^\infty$ and  
		$$ \norme{U}_{\mcS^\infty} \leqslant (m^*)! K_{m,\varepsilon_1,\varepsilon_2}  \parenth{\norme{\xi}_{L^\infty}+ \norme{\sqrt{\abs{f}}\star W}_{BMO}^{2}}.$$ 
\item[\textbf{(iii)}]  Let us assume that $m>1$. For all $q > m^*=\frac{m}{m-1}$, if $\parenth{\zeta,\displaystyle\int_0^T \abs{f_s} \d s} \in L^{q} \times L^{q}$, then $U \in \mcS^q(\R^d)$ and 
$$\norme{U}_{\mcS^q}^q \leqslant 2^{q-1} K_{m,\varepsilon_1,\varepsilon_2}^q \parenth{\frac{q}{q-m^*}}^{q/m^*} \parenth{\norme{\zeta}_{L^q}^q+\norme{\displaystyle\int_0^T \abs{f_s}\mathrm{d}s}_{L^{q}}^{q}}.$$
\end{itemize}
In the following we will denote simply $\mcK_{q,m,\varepsilon_1,\varepsilon_2}=2^{q-1} K_{m,\varepsilon_1,\varepsilon_2}^q \parenth{\frac{q}{q-m^*}}^{q/m^*}$.
\end{prop}

\begin{Proof}
	The formula \eqref{formula} gives us, for all $t \in [0,T]$: 
$$ \abs{ U_t } \leqslant  \espcond{\abs{S_t^{-1} S_T}\abs{\zeta}}{\mathcal{F}_t} +\espcond{\int_t^T \abs{S_t^{-1}S_s} \abs{f_s}\mathrm{d}s}{\mathcal{F}_t},$$
with 
$$ \mathrm{d} S_t = \sum_{p=1}^k S_t B_t^{(:,p,:)} \mathrm{d} W^p_t + S_t A_t\mathrm{d}t, \quad S_0=I_{d \times d}.$$
	$S$ is the solution of an SDE on $\R^{d\times d}$ for which we can use \cref{ReverseSDEprop}
	by taking, for all $1 \leqslant p \leqslant k$ and $(x,y) \in (\R^{d \times d})^2$, $G^p(s,x)=x B^{(:,p,:)}_s$ and $F(s,y)=y A_s$. Let us note that $\abs{B^{(:,p,:)}} \leqslant \abs{B}$ for all $p \in \ensemble{1,...,k}$. Thus there exists a constant $K_{m,\varepsilon_1,\varepsilon_2}$ such that:	$$ \esp{\sup_{t \leqslant s \leqslant T} \abs{S_t^{-1} S_s}^m}\leqslant K_{m,\varepsilon_1,\varepsilon_2}^m.$$
	
	\begin{itemize}[label=$\diamond$, leftmargin=*, noitemsep]
		\item If $\zeta \in L^\infty$ and $f \in \mcS^\infty$, by using the Hölder inequality we have 
			\begin{align*} \abs{U_t} \leqslant & \norme{\xi}_{L^\infty}K_{m,\varepsilon_1,\varepsilon_2}+\norme{f}_{\mcS^\infty}\espcond{(T-t)\sup_{t \leqslant s \leqslant T} \abs{S_t^{-1}S_s}^{m}}{\mcF_t}^{1/m}\\\leqslant & K_{m,\varepsilon_1,\varepsilon_2}\parenth{\norme{\zeta}_{L^\infty}+T\norme{f}_{\mcS^\infty}}.
		\end{align*}
		
		\item Let us consider $m>1$ and assume that $\zeta \in L^\infty$, $\sqrt{\abs{f}}\star W$ is $BMO$. Then, by using Hölder and energy inequalities
			\begin{align*} \abs{U_t} \leqslant & K_{m,\varepsilon_1,\varepsilon_2}\norme{\xi}_{L^\infty}+K_{m,\varepsilon_1,\varepsilon_2}\espcond{\parenth{\int_t^T \abs{f_s}\d s}^{m^*}}{\mcF_t}^{1/m^*}\\\leqslant & (m^*)! K_{m,\varepsilon_1,\varepsilon_2}\parenth{\norme{\xi}_{L^\infty}+ \norme{\sqrt{\abs{f}}\star W}_{BMO}^{2}}.
		\end{align*}
		\item Let us consider $m>1$ and $q >m^*$. We get, for all $t \in [0,T]$, 
		\begin{align*}
			\abs{U_t}^q \leqslant & 2^{q-1} \parenth{\espcond{\abs{S_t^{-1}S_T}\abs{\zeta}}{\mathcal{F}_t}^q+\espcond{\int_t^T \abs{S_t^{-1}S_s} \abs{f_s}\mathrm{d}s}{\mathcal{F}_t}^q} \\\leqslant & 2^{q-1}\espcond{\abs{S_t^{-1} S_T}^m}{\mcF_t}^{q/m} \espcond{\abs{\zeta}^{m^*}}{\mathcal{F}_t}^{q/m^*}\\&+2^{q-1} \espcond{\sup_{t \leqslant s\leqslant T} \abs{S_t^{-1}S_s}^m}{\mathcal{F}_t}^{q/m} \espcond{\parenth{\int_t^T \abs{f_s}\mathrm{d}s}^{m^*}}{\mathcal{F}_t}^{q/m^*} \\\leqslant & 2^{q-1} K_{m,\varepsilon_1,\varepsilon_2}^q\parenth{\espcond{\abs{\zeta}^{m^*}}{\mathcal{F}_t}^{q/m^*}+\espcond{\parenth{\int_0^T \abs{f_s}\mathrm{d}s}^{m^*}}{\mathcal{F}_t}^{q/m^*}}.
		\end{align*}
The processes $M_t=\espcond{\abs{\zeta}^{m^*}}{\mathcal{F}_t}$ and $N_t=\espcond{\parenth{\displaystyle\int_0^T \abs{f_s}\mathrm{d}s}^{m^*}}{\mcF_t}$ are two martingales with terminal values, respectively given by $\abs{\zeta}^{m^*}$ and $\parenth{\displaystyle\int_0^T \abs{f_s}\mathrm{d}s}^{m^*}$. Hence the Doob maximal inequality gives us, if $q > m^*$,  
	$$ \esp{\sup_{0 \leqslant t \leqslant T} \abs{M_t}^{q/m^*}}=\norme{M}_{\mcS^{q/m^*}}^{q/m^*}\leqslant \parenth{\frac{q}{q-m^*}}^{q/m^*} \norme{M_T}^{q/m^*}_{L^{q/m^*}}=\parenth{\frac{q}{q-m^*}}^{q/m^*} \norme{\zeta}_{L^{q}}^{q},$$
	and
	$$ \esp{\sup_{0 \leqslant t \leqslant T} \abs{N_t}^{q/m^*}}\leqslant\parenth{\frac{q}{q-m^*}}^{q/m^*} \norme{\displaystyle\int_0^T \abs{f_s}\mathrm{d}s}_{L^{q}}^{q}.$$
	So we obtain the announced result:
	$$\norme{U}_{\mcS^q}^q \leqslant 2^{q-1} K_{m,\varepsilon_1,\varepsilon_2}^q \parenth{\frac{q}{q-m^*}}^{q/m^*} \parenth{\norme{\zeta}_{L^q}^q+\norme{\displaystyle\int_0^T \abs{f_s}\mathrm{d}s}_{L^{q}}^{q}}.$$
	\end{itemize}
\end{Proof}

\begin{coro}[Affine upper bound]\label{condfinalegenerale}
Let $m \geqslant 1$. Let us consider $A$ and $B$ adapted, bounded respectively by two real processes $\alpha$ and $\beta$ of the form 
$$ \alpha_s=K+L \mcA_s, \quad  \beta_s=K'+L' \mcB_s,$$
with $(K,L,K',L') \in (\R^+)^4$, $\mcA$, $\mcB$ two non negative real processes such that $\sqrt{\mcA}\star W$ and $\mcB \star W$ are BMO with the condition 
$$2 m L \norme{\sqrt{\mcA}\star W}^2_{BMO}+\sqrt{2} L' \norme{\mcB\star W}_{BMO} C'_m<1.$$
We have the following estimates, with constants $K_m$, $\mcK_{q,m}$ depending only on $m,q,K_y,K_z,L_y,L_z$ and the BMO norms $\norme{\sqrt{\mcA}\star W}_{BMO}$, $\norme{\mcB \star W}_{BMO}$:
\begin{itemize} 
\item[\textbf{(i)}] If $\zeta \in L^\infty (\Omega, \mcF_T)$ and $f \in \mcS^\infty$, then $U \in \mcS^\infty(\R^d)$ and 
$$ \norme{U}_{\mathcal{S}^\infty} \leqslant K_{m} \parenth{\norme{\zeta}_{L^\infty}+T \norme{f}_{\mcS^\infty}},$$
\item[\textbf{(ii)}] Let us assume that $m>1$. If $\zeta \in L^\infty (\Omega, \mcF_T)$, $\sqrt{\abs{f}}\star W\in BMO$, then $U \in \mcS^\infty$ and  
		$$ \norme{U}_{\mcS^\infty} \leqslant (m^*)! K_{m}  \parenth{\norme{\xi}_{L^\infty}+ \norme{\sqrt{\abs{f}}\star W}_{BMO}^{2}}.$$ 
\item[\textbf{(iii)}]  Let us assume that $m>1$. For all $q > m^*=\frac{m}{m-1}$, if $\parenth{\zeta,\displaystyle\int_0^T \abs{f_s} \d s} \in L^{q} \times L^{q}$, then $U \in \mcS^q(\R^d)$ and 
$$\norme{U}_{\mcS^q}^q \leqslant 2^{q-1} K_{m}^q \parenth{\frac{q}{q-m^*}}^{q/m^*} \parenth{\norme{\zeta}_{L^q}^q+\norme{\displaystyle\int_0^T \abs{f_s}\mathrm{d}s}_{L^{q}}^{q}}.$$
\end{itemize}
In the following we will denote simply $\mcK_{q,m}=2^{q-1} K_{m}^q \parenth{\frac{q}{q-m^*}}^{q/m^*}$.\end{coro}

\begin{Proof}
We obtain easily estimates about BMO-norms of $\sqrt{\alpha} \star W$ and $\beta \star W$ by using the triangle inequality, 
$$\Vert \sqrt{\alpha}\star W \Vert_{BMO}\leqslant  \sqrt{K T}+\sqrt{L}\Vert \sqrt{\mcA}\star W \Vert_{BMO}, \quad \Vert \beta \star W \Vert_{BMO}\leqslant K' \sqrt{T} +L'\Vert \mcB \star W\Vert_{BMO} , $$
and it follows that $\sqrt{\alpha} \star W,\beta \star W$ are $BMO$. 
To use \cref{condfinalegene} we just have to show that $\sqrt{\alpha}\star W$ and $\beta \star W$ are respectively $\varepsilon_1$ and $\varepsilon_2$ sliceable with $2m \varepsilon_1^2 +\sqrt{2}\varepsilon_2 C'_m <1$. To this end, we consider the following uniform sequence of deterministic stopping times $$ T_j=j \frac{T}{N}, \quad j \in \ensemble{0,...,N}, $$
and a parameter $\eta >0$. With \cref{BMOstartstop} and defining $\eta=\frac{T}{N}$, previous inequalities become on $[T_i,T_{i+1}]$
\beq \norme{\startstop{T_{i}}{\sqrt{\alpha}\star W}{T_{i+1}}}_{BMO}\leqslant  \sqrt{K \eta}+\sqrt{L} \Vert \sqrt{\mcA}\star W \Vert_{BMO}, \label{sliceA} \eeq
\beq \Vert \startstop{T_{i}}{\beta \star W}{T_{i+1}} \Vert_{BMO}\leqslant K' \sqrt{\eta} +L'\Vert \mcB \star W \Vert_{BMO}. \label{sliceB}\eeq
By taking $\eta$ small enough, we get $2m \varepsilon_1^2 +\sqrt{2}\varepsilon_2 C'_m <1$ since the following upper bound holds true
$$ 2 m L \norme{\sqrt{\mcA}\star W}^2_{BMO}+\sqrt{2} L' \norme{\mcB\star W}_{BMO} C'_m<1.$$
\end{Proof}

\begin{rema}\label{amel}
	In inequalities \eqref{sliceA} and \eqref{sliceB}, we have used that $\norme{\startstop{T_{i}}{\mcB\star W}{T_{i+1}}}_{BMO} \leqslant \norme{\mcB \star W}_{BMO}$ and \linebreak $\norme{\startstop{T_{i}}{\sqrt{\mcA}\star W}{T_{i+1}}}_{BMO} \leqslant \norme{\sqrt{\mcA} \star W}_{BMO}$. We can easily obtain a more general result by replacing the following assumption:
	$\mcA$, $\mcB$ are two positive real processes such that $\sqrt{\mcA}\star W,\mcB \star W$ are BMO with the condition 
$$ 2 m L \norme{\sqrt{\mcA}\star W}^2_{BMO}+\sqrt{2} L' \norme{\mcB\star W}_{BMO} C'_m<1,$$
by the new one: $\mcA$, $\mcB$ are two positive real processes such that $\sqrt{\mcA}\star W,\mcB \star W$ are in $BMO_{\varepsilon_1}$ and $BMO_{\varepsilon_2}$ with the condition 
$$ 2 m L \varepsilon_1^2+\sqrt{2} L'\varepsilon_2 C'_m<1.$$
\end{rema}

\begin{rema}
We have not mentioned the dependence of the constants with respect to $\norme{\sqrt{\mcA} \star W}_{BMO}$ and $\norme{\mcB\star W}_{BMO}$ in notations but we will precise it explicitly when it is important. 
\end{rema}


\section{Stability, existence and uniqueness results for general multidimensional quadratic BSDEs}
\label{section stabilite}

\subsection{Proofs of \cref{existuniq1} and \cref{existuniq1diag}}
$ $
We recall first that $(Y^M, Z^M)$ is the unique solution of the localized BSDE \eqref{quadraloc}. 
To show \cref{existuniq1} we begin to prove the following proposition which gives an uniform estimates for $Z^M$. This is the keystone of our procedure. 

\begin{prop}\label{unifbound}
Let $m>1$. If assumptions \ref{H}---\ref{BMO,m}---\ref{Dxi,b}---\ref{Df,b} hold true then 
$$\sup_M \supess_{\Omega \times [0,T]} |Z^M| <+\infty.$$ 
\end{prop}

\begin{Proof}
\paragraph{Step 1 --- Malliavin differentiation.} We assume that $f$ is continuously differentiable with respect to $(y,z)$. This assumption is not restrictive by considering a smooth regularization of $f$. 

Recalling assumptions \ref{Dxi,b} and \ref{Df,b}, Proposition 5.3 in \cite{ELKaroui:1997dn} gives us that for all $0 \leqslant u \leqslant t \leqslant T$, $Y_t^M$ and $Z_t^M$ are respectively in $\mbbD^{1,2}(\R^d)$ and  $\mbbD^{1,2}(\R^{d\times k})$. Moreover the process $(\D_u Y^M, \D_u Z^M)=(\D_u Y^M_t, \D_u Z^M_t)_{0\leqslant t \leqslant T}$ solves for all $u$ the following linear BSDE in $\R^{d \times k}$:
\begin{align}\label{BSDEderivee}
	 \D_u Y^M_t &= \D_u \xi + \int_t^T \Big(\nabla_y f^M\parenth{s,Y^M_s,Z^M_s} \D_u Y^M_s + \nabla_z f^M\parenth{s,Y^M_s,Z^M_s} \D_u Z^M_s \nonumber\\&+(\D_u f^M)\parenth{s,Y^M_s,Z^M_s}\Big) \mathrm{d} s -\int_t^T \D_u Z^M_s \mathrm{d} W_s,
\end{align}

and $(\D_t Y_t)_{0 \leqslant t \leqslant T}$ is a version of $(Z_t)_{0 \leqslant t \leqslant T}$. In particular, there exists a continuous version of $Z$.
Let us emphasize that BSDE \eqref{BSDEderivee} means that for each $p \in \ensemble{1,...,k}$, 
\begin{align}\label{BSDEderiveep}
	 \D^p_u Y^M_t &= \D^p_u \xi + \int_t^T \Big(\nabla_y f^M\parenth{s,Y^M_s,Z^M_s} \D^p_u Y^M_s + \nabla_z f^M\parenth{s,Y^M_s,Z^M_s} \D^p_u Z^M_s \nonumber\\&+(\D^p_u f^M)\parenth{s,Y^M_s,Z^M_s}\Big) \mathrm{d} s -\int_t^T \D^p_u Z^M_s \mathrm{d} W_s, 
\end{align}

besides $\D^p  Y^M$ is a process with values in $\R^{d}$ for each $p \in \{1,...,k\}$. 

\paragraph{Step 2 --- $\mcS^\infty$-Estimation.}
We are looking for an $\mcS^\infty$-estimate of $\D_u Y^M$ for all $u \in [0,T]$ applying results of section   \ref{tools}. Since $\vert \nabla_z \rho^M(z) \vert \leqslant 1$, we obtain the following inequalities by recalling the main assumption \ref{H}, 
$$ \abs{\nabla_y f^M\left(s,Y^M_s,Z^M_s\right)} = \abs{\nabla_y f\left(s,Y^M_s,\rho^M (Z^M_s)\right)} \leqslant  K_y+L_y \abs{Z^M_s}^2 ,$$
$$ \abs{\nabla_z f^M\left(s,Y^M_s,Z^M_s\right)} = \abs{\nabla_z f\left(s,Y^M_s,\rho^M(Z^M_s)\right)} \leqslant  K_z +2 L_z \abs{Z^M_s}.$$

\noindent Let us consider the two positive processes $\alpha^M$ and $\beta^M$ defined below, 
$$\alpha^M=K_y+L_y\abs{Z^M}^2, \quad \beta^M=K_z+2L_z\abs{Z^M}.$$

For all $p \in \ensemble{1,...,k}$, by recalling \ref{BMO,m}, we can apply \cref{condfinalegenerale} \tbf{(iii)}, to the BSDE \eqref{BSDEderiveep} with the following constants and processes: 
$$L=L_y, \quad K=K_y, \quad K'=K_z, \quad L'=2L_z, \quad \mcA=\abs{Z^M}^2, \quad \mcB=\abs{Z^M}.$$
Thus, we obtain, for all $u \in [0,T]$,
	\begin{align*} \norme{\D_u Y^M}_{\mathcal{S}^\infty} \leqslant & \sum_{p=1}^k   \norme{\D_u^p Y^{M}}_{\mathcal{S}^\infty} \leqslant (m^*)! K_{m} \sum_{p=1}^k \parenth{\norme{\D^p_u \xi}_{L^\infty}+\norme{\sqrt{\abs{\D^p_u f^M(.,Y^M,Z^M)}}\star W}^2_{BMO}}\\\leqslant & C_{m} \parenth{\norme{\D_u \xi}_{L^\infty}+\norme{\sqrt{\abs{\D_u f^M(.,Y^M,Z^M)}}\star W}^2_{BMO}},
	\end{align*}
where $C_m$ does not depend on $M$. Indeed, it is important to remark that the constant $K_m$ given by \cref{condfinalegenerale} depends on $\norme{\abs{Z^M} \star W}_{BMO}$ and so, could depend on $M$. But, by checking the proof of \cref{SDEprop} in the Appendix it is easy to see that the constant $K_m$ given by \cref{condfinalegenerale} is equal to 
$$\displaystyle\sum_{i=0}^{N-1} \parenth{\frac{1}{1-2m L_y \norme{Z^M\star W}^2_{BMO}-2L_z\norme{Z^M\star W}_{BMO} C'_m}}^{i} \leqslant \displaystyle\sum_{i=0}^{N-1} \parenth{\frac{1}{1-2m L_y \mbbK^2-2L_z\mbbK C'_m}}^{i},$$
where $N$ is an integer large enough and the uniform bound with respect to $M$ follows.

Under the assumption \ref{Df,b} together with \ref{BMO,m}, the last term has a $\mcS^\infty$-upper bound uniform with respect to $M$. Indeed we have, for all $(u,t) \in [0,T]^2$,
$$ \espcond{\int_t^T  \abs{\D_u f^M(s,Y^M_s,Z^M_s)}\d s}{\mcF_t} \leqslant C \parenth{T+\norme{\abs{Z^M} \star W}^2_{BMO}},$$
hence we deduce
$$ \sup_M \norme{\sqrt{\abs{\D_u f^M(.,Y^M,Z^M)}}\star W}_{BMO} \leqslant \sqrt{C}\parenth{\sqrt{T}+\sup_M \norme{\abs{Z^M} \star W}_{BMO}}.$$
The last supremum is finite under assumption \ref{BMO,m} and we obtain the announced result since
$$\sup_M\norme{Z^M}_{\mathcal{S}^{\infty}} = \sup_M\norme{(\D_tY_t^M)_{t \in [0,T]}}_{\mathcal{S}^{\infty}} \leqslant \sup_M\sup_u \norme{\D_u Y^M}_{\mathcal{S}^\infty} <+\infty.$$
When  $f$ is not continuously differentiable with respect to $(y,z)$ we consider a smooth regularization of $f$ and we obtain by this classical approximation that
$$\sup_M \supess_{\Omega \times [0,T]} |Z^M| <+\infty.$$ 
\end{Proof}

We are now able to prove \cref{existuniq1}.

\begin{Proof}[of \cref{existuniq1}]
	For the existence result, we can fix $M^{\star} > \sup_{M} \supess_{\Omega \times [0,T]} |Z^M| $ according to \cref{unifbound}. Thanks to assumptions on $f$ and $f^M$, we get
	$$ f^{M^{\star}} \left(s,Y^{M^{\star}}_s, Z^{M^{\star}}_s \right)=f\left(s,Y^{M^{\star}}_s, \rho^{M*} \parenth{Z^{M^{\star}}_s} \right)=f\left(s,Y^{M^{\star}}_s, Z^{M^{\star}}_s \right) \quad \mathbb{P} \otimes [0,T] \quad \text{a.e.}$$
Then, $\left( Y^{M^{\star}}, Z^{M^{\star}}\right)$ becomes a solution of the quadratic BSDE \eqref{quadra} in $\mcS^2(\R^d) \times (\mcS^\infty(\R^{d\times k}) \cap \mcZ^m_{BMO})$.
The uniqueness comes from the classical uniqueness result of Pardoux and Peng \cite{Pardoux:1990ju}: indeed, if we have two solutions $(Y^1,Z^1)$ and $(Y^1,Z^1)$ with $\supess_{\Omega \times [0,T]} |Z^1|+|Z^2|<+\infty$ then they are solution of the Lipschitz localized BSDE \eqref{quadraloc} where $M = \supess_{\Omega \times [0,T]} |Z^1|+|Z^2|$.

\end{Proof}

\begin{rema}
 \label{remexistence-uniqueness-avecepsBMO}
 
 By using \cref{amel}, \cref{existuniq1} can be adapted if we replace \ref{BMO,m} by the following one: there exist a  constant $\mbbK$ and a sequence $0=T_0\leqslant T_1 \leqslant ...\leqslant T_N=T$ of stopping times (that does not depend on $M$) such that 
\begin{itemize}\item[\tbf{(i)}] $2m L_y \mbbK^2+\sqrt{2} L_z\mbbK C'_m<1,$ \text{ or equivalently } $\mbbK < \mbbB^m(L_y,L_z),$
\item[\tbf{(ii)}] for all $i\in\{0,...N-1\}$, $\displaystyle\sup_{M \in \R^+}\norme{ \startstop{T_i}{\abs{Z^M} \star W}{T_{i+1}}}_{BMO} \leqslant \mbbK. $
\end{itemize}
In this case, if all the other assumptions of \cref{existuniq1} are fulfilled, then the quadratic BSDE \eqref{quadra} has a unique solution $(Y,Z) \in \mcS^\infty(\R^d) \times \mcZ^{slic,m}_{BMO}$ such that
$$\supess_{\Omega \times [0,T]} |Z| <+\infty.$$
\end{rema}

We do not give the proof of \cref{existuniq1diag} since it is quite similar to the proof of \cref{existuniq1}. Indeed, the main point is to show that \cref{unifbound} stays true. To do that we just have to mimic its proof and replace the application of \cref{condfinalegenerale} \tbf{(iii)} by a new tailored one adapted to the diagonal framework and proved by using the same strategy as in the proof of \cref{diagstab}.

\subsection{Stability result}

With the classical linearisation tool we can prove a stability theorem for the BSDE \eqref{quadra} by using results of section \ref{tools}. Let us consider two solutions of \eqref{quadra} in $\R^d \times \R^{d \times k}$, denoted $(Y^1,Z^1)$ and $(Y^2,Z^2)$, with terminal conditions $\xi^1$ and $\xi^2$ and generators respectively $f_1$ and $f_2$:
$$Y^1_t = \xi^1+\int_t^T f_1 \parenth{s,Y^1_s,Z^1_s} \mathrm{d}s - \int_t^T Z^1_s \mathrm{d} W_s, \quad 0 \leqslant t \leqslant T,$$
$$Y^2_t = \xi^2+\int_t^T f_2\parenth{s,Y^2_s,Z^2_s} \mathrm{d}s - \int_t^T Z^2_s \mathrm{d} W_s, \quad 0 \leqslant t \leqslant T.$$
We assume that $f_1$,$f_2$ satisfies the usual conditions \ref{H}. Let us denote 
$$\delta Y_s=Y^1_s-Y^2_s, \quad \delta Z_s=Z^1_s-Z^2_s,\quad \delta F_s=f_1(s,Y^1_s,Z^1_s)-f_2(s,Y^2_s,Z^2_s),$$
$$\delta f_s=f_1(s,Y^2_s,Z^2_s)-f_2(s,Y^2_s,Z^2_s) \quad\text{ and } \quad \delta \xi=\xi^1-\xi^2.$$ 

The process $(\delta Y, \delta Z)$ solves the BSDE 
\beq \delta Y_t = \delta \xi + \int_t^T \delta F_s \mathrm{d}s-\int_t^T \delta Z_s \mathrm{d}W_s, \quad \quad 0 \leqslant t \leqslant T. \label{eqdifference}\eeq

\begin{theo}[Stability result] \label{Stability (1)}
	Let $m >1$, $p > \frac{m^*}{2}$ and let us suppose that 
	\begin{itemize}
		\item[\textbf{(i)}] $ 2 m L_y \norme{\abs{Z^1} \star W}_{BMO}^2 + \sqrt{2} L_z \parenth{\norme{\abs{Z^1} \star W}_{BMO}+\norme{\abs{Z^2} \star W}_{BMO}} C'_m <1$
		or $(Z^1,Z^2) \in \mcZ_{BMO}^m\times \mcZ_{BMO}^m$,
		\item[\textbf{(ii)}] $(\xi_1,\xi_2)\in (L^{2p})^2,$
		\item[\textbf{(iii)}]  $\dint_0^T\abs{\delta f_s}\d s \in L^{2p}$.
	\end{itemize}
	Then, there exists a constant $\widetilde{K_{p}}\parenth{\norme{\abs{Z^1} \star W}_{BMO},\norme{\abs{Z^2} \star W}_{BMO}}$ (depending only on $p, K_y, L_y, K_z, L_z,T$ and the BMO norms of $\abs{Z^1} \star W$ and $\abs{Z^2} \star W$) such that 
$$ \norme{\delta Y}_{\mcS^{2p}}^p+\norme{\abs{\delta Z} \star W}_{\mcH^{p}}^p \leqslant \widetilde{K_{p}}\parenth{\norme{\abs{Z^1} \star W}_{BMO},\norme{\abs{Z^2} \star W}_{BMO}} \parenth{ \norme{\delta \xi}_{L^{2p}}^p +\norme{\int_0^T\abs{\delta f_s}\d s}_{L^{2p}}^p}.$$
\end{theo}

\begin{Proof}
We firstly assume that
$$ 2 m L_y \norme{\abs{Z^1} \star W}_{BMO}^2 + \sqrt{2} L_z \parenth{\norme{\abs{Z^1} \star W}_{BMO}+\norme{\abs{Z^2} \star W}_{BMO}} C'_m <1.$$
By using the classical linearisation tool, we can rewrite \eqref{eqdifference} as
	$$ \delta Y = \delta \xi + \int_t^T \left( A_s \delta Y_s+B_s (\delta Z_s) +\delta f_s\right) \mathrm{d}s-\int_t^T \delta Z_s \mathrm{d}W_s,$$
	where 
\begin{itemize}[label=$\diamond$, leftmargin=*, noitemsep]
\item $B$ is a $\mcL(\R^{d \times k}, \R^d)$ process defined by blocks by, for all $i \in \{1,...,k\}$,
$$ B^i_s= \begin{cases} \frac{f_1(s,Y^2_s,Z^1_s)-f^1_1(s,Y^2_s,Z^2_s)}{\abs{\delta Z_s}^2} (\transp \delta Z_s^{(:,i)}) &\text{ if } \delta Z_s \neq 0,\\0&\text{ otherwise,}\end{cases}$$
and $B_s \delta Z_s = \sum_{i=1}^k B_s^i \delta Z_s^{(:,i)}$,	

\item $A$ is a $\mcL(\R^{d},\R^d)$-process defined by 
	 $$ A_s = \begin{cases}\frac{f_1(s,Y^1_s,Z^1_s)-f_1(s,Y^2_s,Z^1_s)}{\abs{\delta Y_s}^2}\left( \transp \delta Y_s \right) &\text{ if } \delta Y_s \neq 0, \\ 0 &\text{ otherwise,}\end{cases}$$
	\end{itemize}
Assumption \ref{H} on $f_1$ and $f_2$ gives the following inequalities: 
	$$ \abs{B_s}\leqslant K_z+L_z\left(\abs{Z^1_s}+\abs{Z^2_s}\right),$$
	$$\abs{A _s} \leqslant K_y+L_y\abs{Z^1_s}^2.$$
	
\paragraph{Step 1 – Control of $\delta Y$.} $A$ and $B$ are bounded respectively by two real processes $\alpha$ and $\beta$ defined by
$$\alpha=K_y+L_y\vert Z^1 \vert^2, \quad \beta= K_z+L_z\left(\vert Z^1 \vert+\vert Z^2 \vert\right),$$
and $(\delta Y, \delta Z)$ solves a linear BSDE of the form \eqref{BSDElingen} with $\delta f$ instead of $f$. 
We can apply \cref{condfinalegenerale}, (\tbf{iii}) with 
$$\mcB=\abs{Z^1}+\abs{Z^2}, \quad \mcA=\abs{Z^1}^2,\quad L'=L_z, \quad K=K_y, \quad K'=K_z, \quad \text{ and } \quad L=L_y,$$
which gives, for all $q >1$ such that $q > m^*$,
\beq \norme{\delta Y}_{\mcS^q}^q \leqslant \mcK_{q,m} \parenth{ \norme{\delta \xi}^q_{L^{q}} +\norme{\int_0^T \abs{\delta f_s}\mathrm{d}s}_{L^{q}}^q}.\label{deltaY}\eeq

\paragraph{Step 2 – Control of $\delta Z$.}  The Itô formula applied to $\abs{\delta Y}^2$ gives us
\begin{align}
	 \int_0^T \abs{\delta Z_s}^2 \d s &= \abs{\delta \xi}^2 -\abs{\delta Y_0}^2-2 \int_0^T \delta Y_s . (\delta Z_s \d W_s)+2 \int_0^T (\delta Y . \delta F)_s \d s \nonumber\\ &\leqslant \abs{\delta \xi}^2-2\int_0^T \delta Y_s . (\delta Z_s \d W_s)+2\int_0^T (\delta Y . \delta F)_s \d s.\label{itocarrez}
\end{align}
Recalling assumption \ref{H} we have
$$ \abs{\delta F_s}=\abs{f_1(s,Y^1_s,Z^1_s)-f_2(s,Y^2_s,Z^2_s)}\leqslant \left(K_y+L_y\abs{Z^1_s}^2 \right) \abs{\delta Y_s}+\left(K_z+L_z\left(\abs{Z^1_s}+\abs{Z^2_s}\right)\right) \abs{\delta Z_s}+\abs{\delta f_s}.$$
With the Cauchy-Schwarz and Young inequalities, we get 
\begin{align*}
 {}&2\int_0^T (\delta Y . \delta F)_s \d s \leqslant  2\int_0^T \abs{\delta Y_s} \abs{\delta F_s} \d s\\\leqslant & 2\int_0^T \left[\left(K_y+L_y\abs{Z^1_s}^2 \right) \abs{\delta Y_s}^2+\left(K_z+L_z\left(\abs{Z^1_s}+\abs{Z^2_s}\right)\right) \abs{\delta Y_s}\abs{\delta Z_s}+\abs{\delta f_s} \abs{\delta Y_s} \right] \d s\\\leqslant &  2\parenth{\sup_{0 \leqslant s \leqslant T} \abs{\delta Y_s}^2}\int_0^T\left[K_y+L_y\abs{Z^1_s}^2 +\left(K_z+L_z\left(\abs{Z^1_s}+\abs{Z^2_s}\right)\right)^2+\frac{1}{2}\right]\d s+\frac{1}{2} \int_0^T \abs{\delta Z_s}^2 \d s+\parenth{\int_0^T  \abs{\delta f_s} \d s}^2.
\end{align*}

\noindent By using this last inequality in \eqref{itocarrez} we obtain
\begin{align*} \frac{1}{2}\int_0^T \abs{\delta Z_s}^2 \d s \leqslant &\abs{\delta \xi}^2-2\int_0^T \delta Y_s . (\delta Z_s \d W_s)\\&+  2\parenth{\sup_{0 \leqslant s \leqslant T} \abs{\delta Y_s}^2}\int_0^T\left[K_y+L_y\abs{Z^1_s}^2 +\left(K_z+L_z\left(\abs{Z^1_s}+\abs{Z^2_s}\right)\right)^2+\frac{1}{2}\right]\d s+\parenth{\int_0^T  \abs{\delta f_s} \d s}^2.\end{align*}

\noindent Thus, for all $p \geqslant 1$, there exists a constant $K$ depending only on $p$ such that 
\begin{align*}
	{}&\norme{\abs{\delta Z} \star W}_{\mcH^p}^p \leqslant K \Bigg[ \norme{\delta \xi}^p_{L^p}+\esp{\left(\sup_{0\leqslant t \leqslant T} \abs{\int_0^t \delta Y_s .(\delta Z_s \d W_s)}\right)^{p/2}}\\&+\esp{\left(\sup_{0 \leqslant s \leqslant T} \abs{\delta Y_s}^2  \int_0^T\left[K_y+L_y\abs{Z^1_s}^2+\frac{1}{2} \left(K_z+L_z\left(\abs{Z^1_s}+\abs{Z^2_s}\right)\right)^2+1\right]\d s\right)^{p/2}}+\esp{\parenth{\int_0^T \abs{\delta f_s} \d s}^{p}} \Bigg].
\end{align*}

\noindent In the following we keep the notation $K$ for all constants appearing in the upper bounds.
Then, according to the $BDG$ inequalities, we get for all $p \geqslant 1$:
$$ \esp{\left(\sup_{0\leqslant t \leqslant T}\abs{\int_0^t \delta Y_s.( \delta Z_s  \d W_s)}\right)^{p/2}}=\norme{\int_0^. \delta Y_s.( \delta Z_s  \d W_s)}_{\mcS^{p/2}}^{p/2}\leqslant (C'_{p/2})^{p/2} \norme{\int_0^. \delta Y. (\delta Z_s \d W_s)}_{\mcH^{p/2}}^{p/2}.$$
Since we have
$$\norme{\int_0^. \delta Y. (\delta Z_s \d W_s)}_{\mcH^{p/2}}^{p/2}=\esp{\left(\sum_{i=1}^k \int_0^T \left( \delta Y_s.\delta Z_s^{(:,i)}\right)^2\d s\right)^{p/4}}\leqslant \esp{ \left(\sup_{0 \leqslant s \leqslant T} \abs{\delta Y_s}^2\times \int_0^T \abs{\delta Z_s}^2\d s\right)^{p/4}},$$
then the Cauchy-Schwartz inequality gives us
$$\esp{\left(\sup_{0\leqslant t \leqslant T} \abs{\int_0^t \delta Y_s.( \delta Z_s  \d W_s)}\right)^{p/2}}\leqslant (C'_{p/2})^{p/2} \norme{\abs{\delta Z} \star W}_{\mcH^{p}}^{p/2} \norme{\delta Y}_{\mcS^{p}}^{p/2}.$$

\noindent Moreover we obtain with Cauchy-Schwarz and Young inequalities:
\begin{align*}
	{}&\norme{\abs{\delta Z} \star W}_{\mcH^p}^p \\\leqslant &K\Bigg[ \norme{\delta \xi}^p_{L^p}+ \norme{\abs{\delta Z} \star W}_{\mcH^{p}}^{p/2} \norme{\delta Y}_{\mcS^{p}}^{p/2}+\\&\norme{\delta Y}_{\mcS^{2p}}^{p}\esp{\left(\int_0^T \left[ 1+\abs{Z^1_s}^2+\abs{Z^2_s}^2\right]\d s\right)^{p}}^{1/2}+\esp{ \parenth{\int_0^T \abs{\delta f_s} \d s}^p}\Bigg]\\\leqslant & K\Bigg[ \norme{\delta \xi}^p_{L^p}+\norme{\delta Y}_{\mcS^{p}}^{p}+\norme{\delta Y}_{\mcS^{2p}}^{p}\esp{\left(\int_0^T \left[ 1+\abs{Z^1_s}^2+\abs{Z^2_s}^2\right]\d s\right)^{p}}^{1/2}+\esp{ \parenth{\int_0^T \abs{\delta f_s} \d s}^{p}} \Bigg]\\&+\frac{1}{2} \norme{\abs{\delta Z} \star W}_{\mcH^{p}}^{p}.
\end{align*}

The energy inequality allows us to bound $\esp{\left(\dint_0^T \left[ 1+\abs{Z^1_s}^2+\abs{Z^2_s}^2\right]\d s\right)^{p}}$ by
 $$K\parenth{1+\norme{\abs{Z^1} \star W}_{BMO}^{2p}+\norme{\abs{Z^2} \star W}_{BMO}^{2p}},$$ 
which is finite recalling assumption \tbf{(i)}.
\noindent Finally, for all $p \geqslant 1$, there exists a constant $K$ (which depends only on $p, K_y, L_y, K_z, L_z,T$ and the BMO norms of $\abs{Z^1} \star W$, $\abs{Z^2} \star W$) such that 
\beq\norme{\abs{\delta Z} \star W}_{\mcH^p}^p \leqslant K\parenth{\norme{\delta \xi}_{L^p}^p+\norme{\delta Y}_{\mcS^{2p}}^{p}+\norme{\int_0^T \abs{\delta f_s} \d s}_{L^p}^{p}}. \label{deltaZ}
\eeq

\paragraph{Step 3 – Stability.}
Considering $p > \frac{m^*}{2}$ and combining \eqref{deltaY} where $q=2p$ with \eqref{deltaZ}, we obtain existence of a constant $\widetilde{K_{p}}\parenth{\norme{\abs{Z^1} \star W}_{BMO},\norme{\abs{Z^2} \star W}_{BMO}}$ which depends only on $p, K_y, L_y, K_z, L_z,T,K$ and the BMO norms of $\abs{Z^1} \star W$, $\abs{Z^2} \star W$ such that 
$$ \norme{\delta Y}_{\mcS^{2p}}^p+ \norme{\abs{\delta Z} \star W}_{\mcH^{p}}^p \leqslant \widetilde{K_{p}}\parenth{\norme{\abs{Z^1} \star W}_{BMO},\norme{\abs{Z^2} \star W}_{BMO}} \parenth{ \norme{\delta \xi}_{L^{2p}}^p +\norme{\int_0^T\abs{\delta f_s}\d s}_{L^{2p}}^p}.$$
\paragraph{Step 4 – The alternative assumption \tbf{(i)}} We can deal with assumption $(Z^1,Z^2) \in \mcZ_{BMO}^m\times \mcZ_{BMO}^m$ by changing the linearization step in the proof. We can remark that $\delta F_s=\widetilde{A_s}\delta Y_s+\widetilde{B_s}\delta Z_s+\delta f_s$, where 
	$$ \widetilde{A_s}\delta Y_s = \frac{1}{2} \parenth{f^1(s,Y^1_s,Z^1_s)-f^1(s,Y^2_s,Z^1_s)+f^1(s,Y^1_s,Z^2_s)-f^1(s,Y^2_s,Z^2_s)}, $$
	$$ \widetilde{B_s}\delta Z_s = \frac{1}{2}\parenth{f^1(s,Y^2_s,Z^1_s)-f^1(s,Y^2_s,Z^2_s)+f^1(s,Y^1_s,Z^1_s)-f^1(s,Y^1_s,Z^2_s)},$$
	and we get symmetric bounds for $ \widetilde{A}$ and $ \widetilde{B}$:
		$$ \abs{\widetilde{A_s}} \leqslant K_y+\frac{L_y}{2}\parenth{\abs{Z^1_s}^2+\abs{Z^2_s}^2}, \quad \abs{\widetilde{B_s}}\leqslant K_z+L_z\left(\abs{Z^1_s}+\abs{Z^2_s}\right).$$
	Then \tbf{(i)} becomes 
		$$ mL_y \parenth{\norme{\abs{Z^1} \star W}^2_{BMO}+\norme{\abs{Z^2} \star W}^2_{BMO}}+\sqrt{2} L_z C'_m \parenth{\norme{\abs{Z^1} \star W}_{BMO}+\norme{\abs{Z^2} \star W}}_{BMO} <1,$$
	wich is fulfilled as soon as we have $(Z^1,Z^2) \in \mcZ_{BMO}^m\times \mcZ_{BMO}^m$. 
\end{Proof}

\begin{rema}
	A more restrictive stability result is already obtained in \cite{Kramkov:2016hc} (see Theorem 2.1). 
\end{rema}

\begin{rema}
\label{remarque stabilite sliceable}
 By using \cref{amel}, it is clear that \cref{Stability (1)} stays true when $Z^1$ and $Z^2$ are only in $\mcZ_{BMO}^{slic,m}$. Indeed, if we denote $0=T^j_0 \leqslant T^j_1\leqslant ...\leqslant T^j_{N^j}$ the sequence of stopping times associated to $Z^j\star W$ for $j \in \{1,2\}$, we can define a new common sequence of stopping times:
 $$T_{k N^1+i} = (T^1_i\vee T^2_{k}) \wedge T^2_{k+1}, \quad i \in \{0,...,N^1-1\}, \, k \in \{0,...,N^2-1\}.$$
 Then, by applying the stability result on each interval $[T_i,T_{i+1}]$ for $i \in \{0,N^1N^2-1\}$ we obtain
 $$ \norme{\delta Y}_{\mcS^{2p}}^p+ \norme{\abs{\delta Z} \star W}_{\mcH^{p}}^p \leqslant N^1N^2\prod_{k=0}^{N^1N^2-1}\widetilde{K_{p}}\parenth{\norme{\startstop{T_k}{\abs{Z^{1}}\star W}{T_{k+1}}}_{BMO},\norme{\startstop{T_k}{\abs{Z^{2}}\star W}{T_{k+1}}}_{BMO}} \parenth{ \norme{\delta \xi}_{L^{2p}}^p +\norme{\int_0^T\abs{\delta f_s}\d s}_{L^{2p}}^p}.$$
 Obviously, when sequences of stopping times are the same for $\abs{Z^1} \star W$ and $\abs{Z^2} \star W$, we can use it directly as the common sequence of stopping time.
\end{rema}

\subsection{Stability result for the diagonal quadratic case} \label{Stabdiagsubsect}

We give here a specific result when the quadratic growth of $z$ has essentially a diagonal structure: we assume that assumption \ref{Hdiag} is in force. As explained in section \ref{subsection enonce existence et unicite Malliavin borne}, this kind of framework has been introduced by Hu and Tang in \cite{Hu:2016is} (see also \cite{Jamneshan:2014uia}). 

To simplify notations in this paragraph, the line $i$ of $z$ will be denoted in a simple way by $(z)^i$, or $z^i$ if there is no ambiguity, instead of $z^{(i,:)}$.
Let us consider two solutions $(Y^1,Z^1)$ and $(Y^2,Z^2)$ which correspond to terminal conditions $\xi^1,\xi^2$ and generators $f_1=f_{\tbf{diag},1}+g_1, f_2=f_{\tbf{diag},2}+g_2$. We have for all $i \in \ensemble{1,...,d}$,
$$ \delta Y^i_t=\delta \xi^i+\int_t^T \delta F_s^i \d s -\int_t^T \delta Z^i_s. \d W_s,$$
with
$$\delta Y_s:=Y^1_s-Y^2_s, \quad \delta Z_s:=Z^1_s-Z^2_s,$$
$$\delta F_s:=\parenth{f_{\tbf{diag},1}(s,Z^1_s)-f_{\tbf{diag},2}(s,Z^2_s)}+\parenth{g_1(s,Y^1_s,Z^1_s)-g_2(s,Y^2_s,Z^2_s)}\quad\text{ and } \quad \delta \xi:=\xi^1-\xi^2.$$
We also define
$$ \delta f_s=f_1(s,Y^2_s,Z^2_s)-f_2(s,Y^2_s,Z^2_s).$$ 

\begin{theo}[Stability result for the diagonal quadratic case] \label{diagstab}
	Let us assume that
	\begin{itemize}
	\item[\tbf{(i)}] $f_1$ and $f_2$ satisfy \ref{Hdiag},
	\item[\tbf{(ii)}] there exists $\mbbB > 0$ such that $(Y^1,Z^1)$ and $(Y^2,Z^2)$ are in $\mcS^\infty (\R^d) \times BMO_{\mbbB}$ and	
	\begin{equation} \label{condition structure cas diagonal}
	c_2^2 d L_{d,y} \mbbB^2 <1, \quad \parenth{\frac{c_2}{c_1}\sqrt{L_{d,y}} +\frac{2\sqrt{d}c_2^2}{c_1^2}L_{d,z}}\frac{4 \sqrt{d}c_2^2L_{d,z} \mbbB^2}{1-c_2^2 d L_{d,y} \mbbB^2} <1, 
	\end{equation}
	where $c_1$ and $c_2$ are given by \cref{BMOequiv} with $B=2L_d\mbbB$.
	\end{itemize}
	Then there exists a constant $\widetilde{K^{\tbf{diag}}} \parenth{\norme{\abs{Z^1} \star W}_{BMO},\norme{\abs{Z^2} \star W}_{BMO}} $ depending only on $\mbbB$ and constants in \ref{Hdiag} such that
	$$ \norme{\delta Y}_{\mcS^\infty} + \norme{\abs{\delta Z} \star W}_{BMO} \leqslant \widetilde{K^{\tbf{diag}}}\parenth{\norme{\abs{Z^1} \star W}_{BMO},\norme{\abs{Z^2} \star W}_{BMO}} \parenth{\norme{\delta \xi}_{L^\infty}+\norme{\sqrt{\abs{\delta f}}\star W}^2_{BMO}}.$$
\end{theo}

\begin{rema}
 For a given $L_d$ and a given $\mbbB$, condition \eqref{condition structure cas diagonal} is fulfilled as soon as $L_{d,y}$ and $L_{d,z}$ are small enough. 
\end{rema}

\begin{Proof}
\paragraph{Step 1 – Control of $\delta Y$.}
	We write $\delta F^i$ as 
	\begin{align*} \delta F^i_s &= \parenth{f_{\tbf{diag},1}(s,Z^1_s)-f_{\tbf{diag},1}(s,Z^2_s)}^i+\parenth{g_1(s,Y^1_s,Z^1_s)-g_1(s,Y^2_s,Z^2_s)}^i+\delta f_s^i \\
	&=\beta^i_s \delta Z^i_s+\alpha^i_s\delta Y_s+\tr(\gamma^i_s\delta Z_s)+\delta f_s^i,	
	\end{align*}
	where $\beta^i,\alpha^i$ and $\gamma^i$ are defined by:
	$$\beta^i_s=\begin{cases} \frac{f^i_{\tbf{diag},1}(s,Z^1_s)-f^i_{\tbf{diag},1}(s,Z^2_s)}{\abs{ \delta Z^{i}_s}^2} \parenth{\transp \delta Z^{i}_s} &\text{ if } \delta Z^{i}_s \neq 0, \\ 0 &\text{ otherwise}\end{cases},$$
		$$\alpha^i_s=\begin{cases} \frac{g^i_{1}(s,Y^1_s,Z^1_s)-g^i_{1}(s,Y^2_s,Z^1_s)}{\abs{\delta Y_s}^2}\parenth{\transp \delta Y_s} &\text{ if } \delta Y_s \neq 0 \\ 0 &\text{ otherwise}\end{cases}, \quad \gamma^i_s=\begin{cases} \frac{g^i_{1}(s,Y^2_s,Z^1_s)-g^i_{1}(s,Y^2_s,Z^2_s)}{\abs{\delta Z_s}^2}\parenth{\transp \delta Z_s} &\text{ if } \delta Z_s \neq 0 \\ 0 &\text{ otherwise}\end{cases}.$$  
	Since we have the following estimate on $\beta ^i$, for all $i \in \ensemble{1,...,d}$, 
	$$\abs{\beta^i} \leqslant L_d \parenth{\abs{(Z^{1})^i}+\abs{(Z^{2})^i}},$$
	and that $\parenth{|(Z^{1})^i| \star W,|(Z^{2})^i| \star W} \in BMO \times BMO$, we deduce that $|\beta^i| \star W$ is BMO too and $\mcE \parenth{\beta^i\star W}$ is an uniform integrable martingale. Consequently we can apply the Girsanov theorem:
	\begin{align}
	\delta Y^i_t &=\delta \xi^i + \int_t^T \parenth{\alpha^i_s\delta Y_s+\tr(\gamma^i_s\delta Z_s)+\delta f_s^i }\d s - \int_t^T \delta Z^i_s . \parenth{ \d W_s - \beta^i_s \d s} \nonumber \\&= \delta \xi^i + \int_t^T \parenth{\alpha^i_s\delta Y_s+\tr(\gamma^i_s\delta Z_s)+\delta f_s^i} \d s - \int_t^T \delta Z^i_s . \d  \overline{W_s}^i\label{diaglin},
	\end{align}
	where $\overline{W}^i$ is a Brownian motion with respect to the probability $\Q^i$ defined by $\d \Q^i=\mcE \parenth{\beta^i\star W}_T \d \PP$. Taking the $\Q^i$-conditional expectation we get 
	$$ \delta Y^i_t = \espcondQi{\delta \xi^i + \int_t^T \parenth{\alpha^i_s\delta Y_s+\tr(\gamma^i_s\delta Z_s)+\delta f_s^i}\d s}{\mcF_t}{i}.$$

Following estimates hold true:	
$$ \abs{\alpha}\leqslant K_{d,y}+ L_{d,y}\abs{Z^1}^2,\quad \abs{\gamma} \leqslant K_{d,z}+L_{d,z}\parenth{\abs{Z^1}+\abs{Z^2}},$$
	and consequently, we obtain
	$$ \norme{\sqrt{\abs{\alpha}}\star W}_{BMO}^2 \leqslant  K_{d,y}T +L_{d,y} \norme{\abs{Z^1} \star W}_{BMO}^2, \quad  \norme{\abs{\gamma}\star W}_{BMO} \leqslant  K_{d,z}\sqrt{T}+L_{d,z} \parenth{\norme{\abs{Z^1} \star W}_{BMO}+\norme{\abs{Z^2} \star W}_{BMO}}.$$
	
Since $\sqrt{\left(\abs{\alpha^i}\abs{\delta Y}+\abs{\gamma^i}\abs{\delta Z}+\abs{\delta f^i}\right)}\star W$ is a BMO martingale, we can apply \cref{BMOequiv}: there exists a constant $c_2$ that depend only on $L_d$ and $\mbbB$ such that

\begin{align*}
\abs{\delta Y^i_t}\leqslant & \norme{\delta \xi^i}_{L^\infty} + c_2^2\parenth{ \norme{\delta Y}_{\mcS^\infty}\norme{\sqrt{\abs{\alpha^i}}\star W}_{BMO}^{2}+\norme{ \abs{\delta Z} \star W}_{BMO}\norme{\abs{\gamma^i} \star W}_{BMO}+\norme{\sqrt{\abs{\delta f^i}}\star W}^2_{BMO}},\nonumber
\end{align*}
and consequently we get
\begin{align} 
	\norme{\delta Y}_{\mcS^\infty} &\leqslant \sum_{i=1}^d \norme{\delta Y^i}_{\mcS^\infty} \nonumber\\
	&\leqslant \sqrt{d}\norme{\delta \xi}_{L^\infty} + c_2^2\parenth{d\norme{\delta Y}_{\mcS^\infty}\norme{\sqrt{\abs{\alpha}}\star W}_{BMO}^{2}+d\norme{\abs{\delta Z} \star W}_{BMO}\norme{\abs{\gamma}\star W}_{BMO}+d\norme{\sqrt{\abs{\delta f}}\star W}^2_{BMO}}\nonumber\\
	&\leqslant \sqrt{d}\norme{\delta \xi}_{L^\infty}+c_2^2 d\norme{\delta Y}_{\mcS^\infty}\parenth{K_{d,y}T+L_{d,y} \norme{\abs{Z^1} \star W}_{BMO}^2}\nonumber\\
	&+c_2^2\norme{\abs{\delta Z} \star W}_{BMO} \parenth{ dK_{d,z}\sqrt{T}+dL_{d,z} \parenth{\norme{\abs{Z^1} \star W}_{BMO}+\norme{\abs{Z^2} \star W}_{BMO}}}+c_2^2 d\norme{\sqrt{\abs{\delta f}}\star W}^2_{BMO} \label{quadnonslice}.
\end{align}

As in the proof of \cref{SDEprop}, now we slice $[0,T]$ in small pieces. We consider $\eta = \frac{T}{N}$ with $N\in \N^*$ and we set $T_i=i\eta$ for $i \in \{0,...,N\}$.
The process $\delta Y$ is equal to \beq \delta Y_t=\delta Y_T\indicat_{\{T\}}(t)+\sum_{k=1}^{N-1} \widetilde{\delta Y}^k_t \indicat_{[T_k,T_{k+1}[}(t),\quad \text{with} \quad \widetilde{\delta Y}^k_t =  \delta Y_t\indicat_{[T_k,T_{k+1}]}(t), \, k \in \{0,...,N-1\} \label{sumsliceY}.\eeq

On the interval $[T_k,T_{k+1}]$ the inequality \eqref{quadnonslice} becomes:
\begin{align*}\norme{\widetilde{\delta Y}^k}_{\mcS^\infty} \leqslant & \sqrt{d}\norme{\widetilde{\delta Y}^k_{T_{k+1}}}_{L^\infty}+ dc_2^2 \norme{\sqrt{\abs{\delta f}}\star W}^2_{BMO} \nonumber\\
&+c_2^2 \Bigg(d \norme{\widetilde{\delta Y}^k}_{\mcS^\infty}\parenth{K_{d,y}\eta+L_{d,y} \mbbB^2}+\norme{\startstop{T_k}{\abs{\delta Z}\star W}{T_{k+1}}}_{BMO} \parenth{ dK_{d,z}\sqrt{\eta}+2 d L_{d,z} \mbbB}\Bigg).
\end{align*}
Then, we can choose $N$ large enough to get $1-c_2^2d(K_{d,y}\eta+L_{d,y} \mbbB^2)>0$. Finally we obtain
\begin{align}\label{diagYcontrol}
\norme{\widetilde{\delta Y}^k}_{\mcS^\infty} &\leqslant \frac{\sqrt{d}}{1-c_2^2 d (K_{d,y}\eta+L_{d,y} \mbbB^2)}\norme{\widetilde{\delta Y}^k_{T_{k+1}}}_{L^\infty}\nonumber\\&+ \norme{\startstop{T_k}{\abs{\delta Z} \star W}{T_{k+1}}}_{BMO} \frac{c_2^2 d(K_{d,z}\sqrt{\eta}+2 L_{d,z} \mbbB)}{1-c_2^2 d (K_{d,y}\eta+L_{d,y} \mbbB^2)}+\norme{\sqrt{\abs{\delta f}}\star W}^2_{BMO}\frac{c_2^2 d}{1-c_2^2 d (K_{d,y}\eta+L_{d,y} \mbbB^2)}.
\end{align}

\paragraph{Step 2 – Control of $\delta Z$.} 
	Applying the Itô formula for the process $\abs{\delta Y^i}^2$ and 
	taking the $\Q^i$-conditional expectation, we get for all $t \in [0,T]$,
\begin{align*}\abs{\delta Y_t^i}^2+\espcondQi{\int_t^{T} \abs{\delta Z^i_s}^2 \d s}{\mcF_t}{i}=&\abs{\delta \xi^i}^2+2 \espcondQi{\int_t^T \parenth{\delta Y^i_s. (g_1(s,Y^1_s,Z^1_s)-g_1(s,Y^2_s,Z^2_s))^i} \d s}{\mcF_t}{i}\\&+2\espcondQi{ \int_t^T \delta Y_s^i.(\delta f^i)_s \d s}{\mcF_t}{i}.	
\end{align*}
Martingales $\sqrt{\abs{\delta Y^i} \abs{(g_1(.,Y^1,Z^1)-g_1(.,Y^2,Z^2))^i}}\star W$ and $\sqrt{\abs{\delta Y^i}{\abs{\delta f^i}}}\star W$ are BMO, since under \ref{Hdiag} we have, for all $s \in [0,T]$,
	$$ \abs{g_1(s,Y^1_s,Z^1_s)-g_1(s,Y^2_s,Z^2_s)} \leqslant \parenth{ K_{d,y}+L_{d,y}\abs{Z^1_s}^2}\abs{\delta Y_s}+\parenth{K_{d,z}+L_{d,z}\parenth{\abs{Z^1_s}+\abs{Z^2_s}}}\abs{\delta Z_s}.$$
By using \cref{BMOequiv}, there exist two constants $c_1>0$ and $c_2>0$ that depend only on $L_d$ and $\mbbB$ such that
\begin{align*}c_1^2 \norme{\abs{\delta Z^i} \star W}_{BMO}^2 \leqslant& \norme{\delta \xi^i}^2_{L^\infty}+2 c_2^2 \norme{ \sqrt{\abs{\delta Y^i} \abs{(g_1(.,Y^1,Z^1)-g_1(.,Y^2,Z^2))^i}}\star W}_{BMO}^2\\
&+2c_2^2 \norme{ \sqrt{\abs{\delta Y^i}{\abs{\delta f^i}}}\star W}_{BMO}^2.
\end{align*}
By summing with respect to $i$ and by using assumption \ref{Hdiag} we obtain
\begin{align*}c_1^2 \norme{\abs{\delta Z} \star W}_{BMO}^2 &\leqslant\norme{\delta \xi}^2_{L^\infty}+2 c_2^2 d \norme{\delta Y}_{\mcS^\infty}  \norme{ \sqrt{ \abs{g_1(.,Y^1,Z^1)-g_1(.,Y^2,Z^2)}}\star W}_{BMO}^2 + 2c_2^2 d \norme{\delta Y}_{\mcS^\infty} \norme{ \sqrt{{\abs{\delta f}}}\star W}_{BMO}^2\\
&\leqslant \norme{\delta \xi}^2_{L^\infty}+2 c_2^2vd \left( T K_{d,y}+ L_{d,y} \norme{\abs{Z^1} \star W}_{BMO}^2 \right)  \norme{\delta Y}_{\mcS^\infty}^2\\
&\quad + 2c_2^2d\left(K_{d,z} T+ L_{d,z}\parenth{\norme{\abs{Z^1} \star W}_{BMO}+\norme{\abs{Z^2} \star W}_{BMO}}  \right)\norme{\abs{\delta Z} \star W}_{BMO} \norme{\delta Y}_{\mcS^\infty} \\
&\quad +2 c_2^2 d\norme{\delta Y}_{\mcS^\infty}  \norme{ \sqrt{{\abs{\delta f}}}\star W}_{BMO}^2.
\end{align*}

Once again, for each $k \in \{0,...,N-1\}$ we can write this inequality on $[T_{k},T_{k+1}]$, and with the same notations as in \eqref{sumsliceY} we obtain

\begin{align*} c_1^2 \norme{\startstop{T_k}{\abs{\delta Z} \star W}{T_{k+1}}}_{BMO}^2 \leqslant& \norme{\widetilde{\delta Y}^k_{T_{k+1}}}^2_{L^\infty}+2c_2^2 d\parenth{\eta K_{d,y}+L_{d,y} \mbbB^2} \norme{\widetilde{\delta Y}^k}^2_{\mcS^\infty}\\
&+2c_2^2d \parenth{\eta K_{d,z} +2 \mbbB L_{d,z}}\norme{\widetilde{\delta Y}^k}_{\mcS^\infty}\norme{\startstop{T_k}{\abs{\delta Z} \star W}{T_{k+1}}}_{BMO}+2 c_2^2 d \norme{\widetilde{\delta Y}^k}_{\mcS^\infty}\norme{\sqrt{\abs{\delta f}}\star W}^2_{BMO}.\end{align*}

We apply Young inequality to the terms $\norme{\startstop{T_k}{\abs{\delta Z} \star W}{T_{k+1}}}_{BMO}$ and  $\norme{\sqrt{\abs{\delta f}}\star W}^2_{BMO}$ and we obtain
\begin{align*}
	{}2 c_2^2 d \parenth{\eta K_{d,z} +2 \mbbB L_{d,z}}\norme{\widetilde{\delta Y}^k}_{\mcS^\infty}\norme{\startstop{T_k}{\abs{\delta Z} \star W}{T_{k+1}}}_{BMO}\leqslant 2 \frac{c_2^4 d^2}{c_1^2}  \parenth{\eta K_{d,z} +2 \mbbB L_{d,z}}^2 \norme{\widetilde{\delta Y}^k}_{\mcS^\infty}^2+\frac{c_1^2}{2} \norme{\startstop{T_k}{\abs{\delta Z} \star W}{T_{k+1}}}_{BMO}^2,
\end{align*}
and for all $\varepsilon >0$
$$ 2 c_2^2 d \norme{\widetilde{\delta Y}^k}_{\mcS^\infty}\norme{\sqrt{\abs{\delta f}}\star W}^2_{BMO} \leqslant \varepsilon^2\norme{\widetilde{\delta Y}^k}_{\mcS^\infty}^2+\frac{c_2^4 d^2}{\varepsilon^2}\norme{\sqrt{\abs{\delta f}}\star W}^4_{BMO}.$$
Consequently we get 
\begin{align}\label{deltaZdiag} \frac{c_1}{\sqrt{2}} \norme{\startstop{T_k}{\abs{\delta Z} \star W}{T_{k+1}}}_{BMO} \leqslant & \norme{\widetilde{\delta Y}^k_{T_{k+1}}}_{L^\infty}+\parenth{\sqrt{2}c_2 \sqrt{d}\parenth{\sqrt{\eta K_{d,y}}+\sqrt{L_{d,y}} \mbbB}+\sqrt{2} \frac{c_2^2 d}{c_1}  \parenth{\eta K_{d,z}+2 \mbbB L_{d,z}}+\varepsilon}\norme{\widetilde{\delta Y}^k}_{\mcS^\infty}\nonumber\\&+\frac{c_2^2 d}{\varepsilon}\norme{\sqrt{\abs{\delta f}}\star W}^2_{BMO}.\end{align}

\paragraph{Step 3 – Stability.}
Combining \eqref{diagYcontrol} and \eqref{deltaZdiag}, we can obtain a stability result on $[T_k,T_{k+1}]$ as soon as  $\eta$ and $\varepsilon$ are sufficiently small to get 
$ c_2^2 d (K_{d,y}\eta+L_{d,y} \mbbB^2) <1$
and 
$$ \frac{\sqrt{2}}{c_1} \parenth{\sqrt{2}c_2\sqrt{d} \parenth{\sqrt{\eta K_{d,y}}+\sqrt{L_{d,y}} \mbbB}+\sqrt{2} \frac{c_2^2 d}{c_1}  \parenth{\eta K_{d,z}+2 \mbbB L_{d,z}}+\varepsilon}\frac{c_2^2 d(K_{d,z}\sqrt{\eta}+2L_{d,z} \mbbB)}{1-c_2^2 d (K_{d,y}\eta+L_{d,y} \mbbB^2)} <1.$$
We obtain the existence of a constant $K$ which does not depend on $k$ such that 
$$ \norme{\widetilde{\delta Y}^k}_{\mcS^\infty} + \norme{\startstop{T_k}{\abs{\delta Z} \star W}{T_{k+1}}}_{BMO} \leqslant K \parenth{\norme{\widetilde{\delta Y}^k_{T_{k+1}}}_{L^\infty}+\norme{\sqrt{\abs{\delta f}}\star W}^2_{BMO}}.$$
Since $\norme{\delta Y}_{\mcS^\infty}\leqslant \sum_{k=0}^{N-1} \norme{\widetilde{\delta Y}^k}_{\mcS^\infty}$ and $\norme{\abs{\delta Z} \star W}_{BMO} \leqslant  \sum_{k=0}^{N-1} \norme{\startstop{T_k}{\abs{\delta Z} \star W}{T_{k+1}}}_{BMO},$ by a direct iteration we finally obtain a constant $K$ such that
$$ \norme{\delta Y}_{\mcS^\infty} + \norme{\abs{\delta Z} \star W}_{BMO} \leqslant K \parenth{\norme{\delta \xi}_{L^\infty}+\norme{\sqrt{\abs{\delta f}}\star W}^2_{BMO}},$$
and $K$ depends only on $\mbbB$ and constants in \ref{Hdiag}. 
\end{Proof}

\subsection{Proof of \cref{withoutxib}}
$ $

\cref{withoutxib} is proved by relaxing assumptions \ref{Dxi,b} and \ref{Df,b} of \cref{existuniq1} thanks to some density arguments.
To ensure the convergence, the keystone result will be the stability \cref{Stability (1)}.

\begin{Proof}[of \cref{withoutxib}]
	\paragraph{Step 1– Approximations.}
	We can approach $\xi$ with a sequence of random variables $(\xi^n)_{n \in \N}$ such that for every $n$, $\xi^n$ has a bounded Malliavin derivative:
$$ \norme{\D \xi^n}_{\mcS^\infty} < \infty.$$
More precisely $\xi^n$ can be chosen of the form $\Phi^n(W_{t_1}, ..., W_{t_n})$ where $\Phi^n \in \mcC^\infty_b(\R^{n})$, $(t_1,...,t_n)\in [0,T]^n$ and $\xi^n$ tends to $\xi$ in every $L^p$ for $p \geqslant 1$ (see \cite{SpringerVerlagGmbHandCoKG:2005we}, Exercise 1.1.7). 

Since $\alpha$ is adapted, we can approach this process with a sequence of sample processes $\alpha^n$ of the form 
$$ \alpha_t^n=\sum_{i=0}^{p_n-1} \alpha_{t_i^n} \indicat_{[t_i^n,t_{i+1}^n[}(t),$$
where $(t_i^n)_{i=0}^{p_n}$ is a sequence of subdivisions of $[0,T]$, with $\sup_{0 \leqslant i \leqslant p_n-1}\abs{t_{i+1}^n-t_i^n} \longrightarrow_{n \rightarrow \infty} 0$, and, for all $0 \leqslant i\leqslant p_n-1$, $n\in\N$, $\alpha^{i,n}$ is  a $\mcF_{t_i^n}$-measurable random variable. We have a convergence of this sequence to $\alpha$ in $L^2(\Omega \times [0,T])$:
$$ \esp{\int_0^T \abs{\alpha^n_s-\alpha_s}^2\d s} \underset{n \rightarrow \infty}{\longrightarrow} 0.$$
We can assume in addition that for all $n$ and for all $0 \leqslant i\leqslant p_n$, $\alpha^{i,n}$ has a bounded Malliavin derivative since this set is dense in $L^2(\Omega)$. It is obvious that for all $0 \leqslant u \leqslant T$ and $0 \leqslant t \leqslant T$, $$ \D_u \alpha^n_t =\dsum_{i=0}^{p_n-1} \D_u\alpha^{i,n}_t \indicat_{[t_i^n,t_{i+1}^n[}(u).$$
According to \cref{propMaliav} applied to $\varphi=\tbf{f}(.,y,z)$, there exists for all $n \in \N$ and $t \in [0,T]$ a bounded random variable $\tbf{G}$ such that
$$ \D_t \tbf{f}(\alpha^n_t,y,z)=\tbf{G}.\D_t \alpha^n_t, \quad \text{ and } \quad  \abs{ \tbf{G}} \leqslant D(1+\abs{z}^2).$$
For each $n \in \N$: $\xi^n$ satisfies \ref{Dxi,b}, $\tbf{f}(\alpha^n_{.},.,.)$ satisfies \ref{Df,b} and \ref{BMO,m} is fulfilled. So, we can apply \cref{existuniq1}: there exists an unique solution $(Y^n,Z^n) \in \mcS^2(\R^d) \times  \mcZ^m_{BMO}$ of the equation
$$Y^{n}_t = \xi^n+\int_t^T \tbf{f}\parenth{\alpha_s^n,Y^{n}_s,Z^{n}_s} \mathrm{d}s - \int_t^T Z^{n}_s \mathrm{d} W_s, \quad 0 \leqslant t \leqslant T.$$

\paragraph{Step 2– Application of the stability result.}
	We can assume that for all $n$, $\norme{\xi^n}_{L^{2m^*}} \leqslant \norme{\xi}_{L^{2m^*}}$. If it is not true, we consider the sequence $\widetilde{\xi^n}=\frac{\norme{\xi}_{L^{2m^*}}}{\norme{\xi}_{L^{2m^*}}+\norme{\xi^n-\xi}_{L^{2m^*}}}\xi^n$ instead of $\xi^n$. The same argument allows us to assume that 
	$$\norme{\alpha^n}_{L^{2m^*}(\Omega\times [0,T])}\leqslant \norme{\alpha}_{L^{2m^*}(\Omega\times [0,T])}.$$
	\noindent Under \ref{BMO2,m}, we have the estimate
	$$mL_y \norme{\abs{Z^{n}} \star W}^2_{BMO}+\sqrt{2}L_z\norme{\abs{Z^{n}} \star W}_{BMO} C'_m \leqslant m L_y \mbbK^2+\sqrt{2}L_z\mbbK C'_m< \frac{1}{2}.$$
	Hence, for all $n_1,n_2 \in \N$, we can use \cref{Stability (1)} for $p=m^*$ which gives us:
	\begin{align*}\norme{Y^{n_1}-Y^{n_2}}_{\mcS^{2m^*}}^{m^*}&+ \norme{\abs{Z^{n_1}- Z^{n_2}}\star W}_{\mcH^{m^*}}^{m^*} \\\leqslant & \widetilde{K}_{m^*}\parenth{\norme{\xi^{n_1}-\xi^{n_2}}_{L^{2m^*}}^{m^*}+\esp{\parenth{\int_0^T \abs{\tbf{f}\parenth{\alpha^{n_1}_t,Y_t^{n_2},Z_t^{n_2}}-\tbf{f}\parenth{\alpha^{n_2}_t,Y_t^{n_2},Z_t^{n_2}}}\d t}^{2m^*}}^{1/2}}.
	\end{align*}	
where the constant $\widetilde{K}_{m^*}$ appearing does not depend on $n$ under \ref{BMO2,m}. This fact was already highlighted in the proof of \cref{unifbound} where an explicit formula for $\widetilde{K}_{m^*}$ was given. We recall that $\suiteN{\xi^n}{n}$ is a Cauchy sequence in $L^{2m^{*}}$, so
	$$\norme{\xi^{n_1}-\xi^{n_2}}_{L^{2m^*}} \underset{n_1,n_2 \rightarrow \infty}{\longrightarrow} 0.$$
For the second term, we use the Hölder inequality:
\begin{align}
	{}&\esp{\parenth{\int_0^T \abs{\tbf{f}\parenth{\alpha^{n_1}_t,Y_t^{n_2},Z_t^{n_2}}-\tbf{f}\parenth{\alpha^{n_2}_t,Y_t^{n_2},Z_t^{n_2}}}\d t}^{2m^*}}\nonumber\\
	\leqslant & D^{2m^*} \esp{\parenth{\int_0^T \parenth{1+\abs{Z_t^{n_2}}^2} \abs{\alpha_t^{n_1}-\alpha_t^{n_2}}^{\delta}\d t}^{2m^*}}\nonumber\\
	\leqslant &D^{2m^*} \norme{\alpha^{n_1}-\alpha^{n_2}}_{\mcS^{4m^*\delta}}^{2m^*\delta} \esp{\parenth{\int_0^T \parenth{1+\abs{Z_t^{n_2}}^2}\d t}^{4m^*}}^{1/2}.\label{alphacont}
\end{align}

Since $\abs{Z^n} \star W\in BMO$, and \ref{BMO2,m} holds true, 
we have
$$\sup_{n_2 \in \N} \esp{\parenth{\int_0^T \parenth{1+\abs{Z_t^{n_2}}^2}\d t}^{4m^*}}^{1/2} <+\infty.$$
Moreover, by using the uniform continuity of trajectories of $\alpha$ on $[0,T]$, we get:
$$ \sup_{t \in [0,T]}\abs{\alpha_t^{n_1}-\alpha_t^{n_2}}\underset{n_1,n_2 \rightarrow \infty}{\longrightarrow} 0 \quad \ps.$$
Since we have
$$ \esp{\sup_{t \in [0,T]}\abs{\alpha_t^{n_1}-\alpha_t^{n_2}}^{4m^*+1}} \leqslant 2C\esp{\sup_{t \in [0,T]}\abs{\alpha_t}^{4m^*+1}}<+\infty,$$
then, a uniform integrability argument gives us 
$$ \esp{\sup_{t \in [0,T]}\abs{\alpha_t^{n_1}-\alpha_t^{n_2}}^{4m^*}} \underset{n_1,n_2 \rightarrow \infty}{\longrightarrow} 0.$$
Finally we get
$$\esp{\parenth{\int_0^T \abs{\tbf{f}\parenth{\alpha^{n_1}_t,Y_t^{n_2},Z_t^{n_2}}-\tbf{f}\parenth{\alpha^{n_2}_t,Y_t^{n_2},Z_t^{n_2}}}\d t}^{2m^*}}\underset{n_1,n_2 \rightarrow \infty}{\longrightarrow} 0.$$

Consequently $(Y^n,Z^n \star W)_{n \in \N}$ is a Cauchy sequence in $\mcS^{2m^*}(\R^d)\times \mcH^{m^*}(\R^{d \times k})$, thus it converges in $\mcS^{2m^*}(\R^d)\times \mcH^{m^*}(\R^{d \times k})$ to a process $(Y,Z\star W)$, and $(Y,Z)$ solves the BSDE \eqref{quadra}. 
Finally, according to \cref{BMOconv} the upper bound for $\norme{\abs{Z^n} \star W}_{BMO}$ holds true for $\norme{\abs{Z} \star W}_{BMO}$ and so the uniqueness follows from the stability theorem.
\end{Proof}

\begin{rema}\label{simplifié}
	If $f$ is a deterministic function, then the assumption \ref{Df,b} is not required.
\end{rema}

\begin{rema}\label{simplifié2}
	If we replace the inequality \eqref{hypbeta} by the new one: there exist $\eta >0$, $D >0$ and $\delta \in (0,1]$ such that for all $(\beta,\beta',y,z) \in (\R^{d'})^2\times \R^d \times \R^{d \times k}$ we have  
	$$\abs{\normalfont{\tbf{f}}(\beta,y,z)-\normalfont{\tbf{f}}(\beta',y,z)}\leqslant D \parenth{1+\abs{z}^{2-\eta}}\abs{\beta-\beta'}^{\delta},$$
	then we do not have to assume that $\alpha$ is a continuous process. Indeed, we can change the inequality \eqref{alphacont} by the following one: by using the Hölder and Cauchy-Schwartz inequalities we have, for all $p > 1$,
	\begin{align*} &\esp{\parenth{\int_0^T \abs{\normalfont{\tbf{f}}\parenth{\alpha^{n_1}_t,Y_t^{n_2},Z_t^{n_2}}-\normalfont{\tbf{f}}\parenth{\alpha^{n_2}_t,Y_t^{n_2},Z_t^{n_2}}}\d t}^{2m^*}} \\
	\leqslant &D^{2m^*} \esp{\parenth{\int_0^T \parenth{1+\abs{Z^{n_2}_t}^{2-\eta}}\abs{\alpha_t^{n_1}-\alpha_t^{n_2}}^{\delta}\d t}^{2m^*}}  \\
	\leqslant &D^{2m^*} \esp{\parenth{\int_0^T \parenth{1+\abs{Z^{n_2}_t}^{2-\eta}}^p \d t}^{2m^*/p}\times \parenth{\int_0^T\abs{\alpha_t^{n_1}-\alpha_t^{n_2}}^{\delta p^*}\d t}^{2m^*/p^*}}\\
	\leqslant &D^{2m^*} \esp{\parenth{\int_0^T \parenth{1+\abs{Z^{n_2}_t}^{2-\eta}}^p \d t}^{2m^*}}^{1/p}\esp{\parenth{\int_0^T\abs{\alpha_t^{n_1}-\alpha_t^{n_2}}^{\delta p^*}\d t}^{2m^*}}^{1/p^*}\\
	\leqslant &D^{2m^*}T^{\frac{4m^*}{p^*}-2m^*}  \esp{\parenth{\int_0^T \parenth{1+\abs{Z^{n_2}_t}^{2-\eta}}^p \d t}^{2m^*}}^{1/p}\norme{\abs{\alpha^{n_1}-\alpha^{n_2}}^{\delta p^*}\star W}^{4m^*/p^*}_{\mcH^{4m^*}}.
	\end{align*}
	With the energy inequality, the first term is uniformly bounded with respect to $n_2$ under the assumption \ref{BMO2,m} by choosing $1 < p \leqslant \frac{2}{2-\eta}$. The second one tends to zero when $n_1,n_2$ go to infinity since the convergence in every $\mcH^r$ for $r>1$ holds true.
\end{rema}

\begin{rema}
 \label{remexistence-uniqueness2-avecepsBMO}
 By using \cref{amel} once again, \cref{withoutxib} can be adapted if we replace the assumption \ref{BMO2,m} by the following one: $\xi \in L^{2m^*}$ and there exist a constant $\mbbK$ and a sequence $0=T_0\leqslant T_1 \leqslant ...\leqslant T_N=T$ of stopping times (that does not depend on $M$) such that 
\begin{itemize}\item[\tbf{(i)}] $2 m L_y \mbbK^2+2\sqrt{2} L_z\mbbK C'_m<1,$ 
\item[\tbf{(ii)}] for all $i\in\{0,...N-1\}$, $\displaystyle\sup_{M \in \R^+}\sup_{\substack{\norme{\eta}_{L^{2m^*}}\leqslant \norme{\xi}_{L^{2m^*}}\\\norme{\beta}_{L^2(\Omega\times [0,T])}\leqslant \norme{\alpha}_{L^2(\Omega\times [0,T])}}}\norme{\startstop{T_i}{\abs{Z^{(M,\eta,\beta)}}\star W}{T_{i+1}} }_{BMO} \leqslant \mbbK. $
\end{itemize}
In this case, if all the other assumptions of \cref{withoutxib} are fulfilled, then the quadratic BSDE \eqref{quadra} has a unique solution $(Y,Z) \in \mcS^\infty(\R^d) \times \mcZ^{slic,m}_{BMO}$ such that
$$\supess_{\Omega \times [0,T]} |Z|<+\infty.$$
\end{rema}

\begin{rema}
\label{extension existence unicite cas diagonal}
 It is possible to extend the existence and uniqueness result for the diagonal case given by \cref{existuniq1diag} to more general terminal conditions and generators. More precisely, it is possible to apply the same strategy as for the proof of \cref{withoutxib} by applying the stability result given by \cref{diagstab} instead of \cref{Stability (1)}. Nevertheless we can only obtain an existence and uniqueness result for terminal conditions (resp. generators) that can be approximated in $L^\infty$ (resp. $BMO$) by terminal conditions satisfying assumption \ref{Dxi,b} (resp. generators satisfying assumption \ref{Df,b}).
\end{rema}


\section{Proofs of section \ref{subsection enonce applications} results}

\label{applications}
\label{section applications}

\subsection{Proof of \cref{prop Tevzadze revisite}}

We start by proving some uniform (with respect to $M$) \tit{a priori} estimates on $\parenth{Y^M,Z^M}$.

\begin{prop}\label{propQ}
	Let us assume that \ref{H} and \ref{HQ} are in force.  Then $\abs{Z^M} \star W \in BMO$, $Y^M \in \mcS^\infty$ and we have the following estimates:
	\begin{itemize}
		\item[\textbf{(i)}] $\norme{\abs{Z^M} \star W}_{BMO}^2 \leqslant \frac{1}{8 \gamma^2} \parenth{1-\sqrt{1-32\gamma^2 \norme{\xi}_{L^\infty}^2}}$,
		\item[\textbf{(ii)}]  $\norme{Y^M}_{\mcS^\infty} \leqslant \norme{\xi}_{L^\infty}+\gamma \norme{\abs{Z^M} \star W}_{BMO}^2$.
	\end{itemize}
\end{prop}

\noindent We can note that upper bound do not depend on $M$.  

\begin{Proof}
To simplify notations in the proof, we skip the superscript $M$ on $(Y^M,Z^M)$ and $f^M$.
	The unique solution $(Y,Z) \in \mcS^2 \times \mcH^2$ of \eqref{quadraloc} can be constructed with a Picard principle as in the seminal paper of Pardoux and Peng (see \cite{Pardoux:1990ju}). We consider a sequence $\suiteN{Y^n,Z^n}{n}$ such that $\suiteN{Y^n,Z^n}{n}$ tends to $(Y,Z)$ in $\mcS^2\parenth{\R^d} \times \mcH^2 \parenth{\R^{d\times k}}$. This sequence is given by
	$$ Y^{n+1}_t=\xi +\int_t^T f\parenth{s,Y^n_s,Z^n_s} \d s-\int_t^T Z^{n+1}_s \d W_s, \quad 0 \leqslant t \leqslant T, \quad (Y^0,Z^0)=(0,0).$$ 
	
	\noindent We will prove with an induction that: for all $n \in \N$, $Y^n \in \mcS^\infty$, $\abs{Z^n} \star W \in BMO$ and 
	$$\norme{\abs{Z^n} \star W}_{BMO}^2 \leqslant \frac{1}{8 \gamma^2} \parenth{1-\sqrt{1-32\gamma^2 \norme{\xi}_{L^\infty}^2}}.$$
	The case $n=0$ is obviously satisfied. Let us suppose that $Y^n \in \mcS^\infty$ and $\abs{Z^n} \star W \in BMO$. Then for all $t \in [0,T]$, under \ref{HQ},
	\beq \abs{Y^{n+1}_t} \leqslant \espcond{\abs{\xi}}{\mcF_t}+\gamma \times \espcond{\int_t^T \abs{Z^n_s}^2\d s}{\mcF_t}.\label{estimeeQ}\eeq
We get $Y^{n+1} \in \mcS^\infty$ since $\norme{Y^{n+1}}_{\mcS^\infty} \leqslant \norme{\xi}_{L^\infty}+\gamma \norme{\abs{Z^n}\star W}_{BMO}^2$. Itô formula gives the following equality 
	 $$ \abs{Y^{n+1}_t}^2=\abs{\xi}^2+ 2 \int_t^T Y^{n+1}_s . f(s,Y^n_s,Z^n_s) \d s - 2\int_t^T Y^{n+1}_s . \parenth{Z_s^{n+1}\d W_s} - \int_t^T \abs{Z^{n+1}_s}^2 \d s.$$
	 By taking conditional expectation we get for every stopping time $\tau$:
	 \begin{align*} \abs{Y^{n+1}_\tau}^2+\espcond{\int_\tau^T \abs{Z^{n+1}_s}^2 \d s}{\mcF_\tau}\leqslant & \espcond{\abs{\xi}^2}{\mcF_t}+2 \espcond{\int_\tau^T Y^{n+1}_s . f(s,Y^n_s,Z^n_s) \d s}{\mcF_\tau}\\\leqslant &  \norme{\xi}^2_{L^\infty}+2 \norme{Y^{n+1}}_{\mcS^\infty} \espcond{\int_\tau^T \gamma \abs{Z^n_s}^2 \d s}{\mcF_\tau}\\\leqslant &  \norme{\xi}^2_{L^\infty}+2\gamma \norme{Y^{n+1}}_{\mcS^\infty} \norme{\abs{Z^n} \star W}_{BMO}^2.
	 \end{align*}
	 Taking the essential supremum with respect to $\tau$ in following inequalities, we obtain
	 \begin{align*}
	 	\abs{Y^{n+1}_\tau}^2 &\leqslant \norme{\xi}^2_{L^\infty}+2\gamma \norme{Y^{n+1}}_{\mcS^\infty} \norme{\abs{Z^n} \star W}_{BMO}^2\\ \espcond{\int_\tau^T \abs{Z^{n+1}_s}^2 \d s}{\mcF_\tau} &\leqslant \norme{\xi}^2_{L^\infty}+2\gamma \norme{Y^{n+1}}_{\mcS^\infty} \norme{\abs{Z^n} \star W}_{BMO}^2.
	 \end{align*}
	Thus $|Z^{n+1}| \star W\in BMO$ and we have:
	 \begin{align*} \norme{Y^{n+1}}^2_{\mcS^\infty}+\norme{\abs{Z^{n+1}}\star W}^2_{BMO}  \leqslant & 2\norme{\xi}^2_{L^\infty}+4\gamma \norme{Y^{n+1}}_{\mcS^\infty} \norme{\abs{Z^n} \star W}_{BMO}^2\\\leqslant & 2\norme{\xi}^2_{L^\infty}+\norme{Y^{n+1}}_{\mcS^\infty}^2+4\gamma^2\norme{\abs{Z^n} \star W}_{BMO}^4,
	 \end{align*}
	 which leads to
	$$ \norme{\abs{Z^{n+1}}\star W}^2_{BMO} \leqslant 2\norme{\xi}^2_{L^\infty} + 4\gamma^2 \norme{\abs{Z^n} \star W}^4_{BMO}.$$ 
	 Using the induction assumption we obtain 
	 $$ \norme{\abs{Z^{n+1}}\star W }^2_{BMO} \leqslant \frac{1}{8 \gamma^2} \parenth{1-\sqrt{1-32\gamma^2 \norme{\xi}_{L^\infty}^2}}.$$ 
	 The induction is achieved. 
	Now we can use  \cref{BMOconv} with $K=\frac{1}{8 \gamma^2} \parenth{1-\sqrt{1-32\gamma^2 \norme{\xi}_{L^\infty}^2}}$: since $Z^{n}\star W$ tends to $Z\star W$ in $\mcH^2$, we conclude that $\norme{\abs{Z}\star W}_{BMO}^2 \leqslant \frac{1}{8 \gamma^2} \parenth{1-\sqrt{1-32\gamma^2 \norme{\xi}_{L^\infty}^2}}$. Finally, we use that $Y^n$ tends to $Y$ in $\mcS^2$ to pass to the limit into \eqref{estimeeQ} and to obtain the final upper bound on $\norme{Y}_{\mcS^\infty}$ .
\end{Proof}

\begin{Proof}[of \cref{prop Tevzadze revisite}]
The proof of the proposition is a direct consequence of \cref{withoutxib} together with \cref{propQ}: since the map $$x \in \R^+ \longmapsto \frac{1}{2\sqrt{2} \gamma} \parenth{1-\sqrt{1-32\gamma^2 x^2}}^{\frac{1}{2}}$$ 
is nondecreasing, the assumption \ref{BMO2,m} is satisfied.
\end{Proof}

\subsection{Proof of \cref{prop Tevzadze monotone}}
$ $

Once again, we start by proving some uniform (with respect to $M$) \tit{a priori} estimates on $(Y^M,Z^M)$.

\begin{prop}
	Let us assume that \ref{H} and \ref{HMon} are in force. Then $\abs{Z^M} \star W \in BMO$, $Y^M \in \mcS^\infty$ and we have 
	\begin{itemize}
		\item[\tbf{(i)}] $ \supess \sup_{t \in [0,T]} \espcond{\dint_t^T e^{-\mu(s-t)}\abs{Z^M_s}^2 \d s}{\mcF_t} \leqslant \frac{1}{8\gamma^2}\parenth{1-\sqrt{1-32\gamma^2 A^2}},$
	\item[\tbf{(ii)}] 	$ \norme{Y^M}_{\mcS^\infty} \leqslant \frac{1}{4\gamma}\parenth{1-\sqrt{1-32\gamma^2 A^2}}+\sqrt{2A^2+ \frac{1}{16 \gamma^2} \parenth{1-\sqrt{1-32\gamma^2 A^2}}^2}.$
	\end{itemize}
		with $A=\max \parenth{\norme{\xi}_{L^\infty}, \frac{\alpha}{\mu}}$. 
\end{prop}

\begin{Proof}
To simplify notations in the proof, we skip the superscript $M$ on $(Y^M,Z^M)$ and $f^M$.
	The unique solution $(Y,Z)$ of \eqref{quadraloc} can  be constructed with a Picard principle. We consider a sequence $\suiteN{Y^n,Z^n}{n}$ such that $\suiteN{Y^n,Z^n \star W}{n}$ tends to $(Y,Z\star W)$ in $\mcS^2\parenth{\R^d} \times \mcH^2 \parenth{\R^{d\times k}}$, with
	$$ Y^{n+1}_t=\xi +\int_t^T f\parenth{s,Y^{n+1}_s,Z^n_s} \d s-\int_t^T Z^{n+1}_s \d W_s, \quad 0 \leqslant t \leqslant T, \quad (Y^0,Z^0)=(0,0), .$$ 
	We can easily show that replacing $Y^{n+1}$ by $Y^n$ in the generator does not affect the convergence of the scheme since $f$ is a Lipschitz function. 
	Moreover, applying Itô formula to $e^{Kt} \abs{Y_t^{n+1}}^2$ with $K$ large enough, we justify with classical inequalities that for all $n \in \N$, $Y^{n+1} \in \mcS^\infty$, with a bound that depend on $M$ for the moment.
	Applying Itô formula to the process $e^{-\mu t} \abs{Y_t^{n+1}}^2$, we obtain
	\begin{align*} e^{-\mu t} \abs{Y_t^{n+1}}^2= &e^{-\mu T} \abs{\xi}^2-\int_t^T \parenth{-\mu e^{-\mu s}\abs{Y_s^{n+1}}^2-2e^{-\mu s} Y_s^{n+1}.f(s,Y^{n+1}_s,Z^n_s)+e^{-\mu s}\abs{Z_s^{n+1}}^2}\d s\\&-2\int_t^T e^{-\mu s}Y^{n+1}_s . \parenth{Z_s^{n+1} \d W_s}.\end{align*}
	Taking conditional expectation, and using assumption \ref{HMon}, we get:
	\begin{align*}\abs{Y_t^{n+1}}^2 \leqslant & e^{-\mu(T-t)} \norme{\xi}_{L^\infty}^2+\espcond{\int_t^T 2e^{-\mu(s-t)}\parenth{\alpha\abs{Y_s^{n+1}}-\frac{\mu}{2} \abs{Y_s^{n+1}}^2+\gamma \abs{Y_s^{n+1}} \abs{Z_s^{n}}^2}\d s}{\mcF_t}\\&-\espcond{\int_t^T e^{-\mu(s-t)} \abs{Z_s^{n+1}}^2 \d s}{\mcF_t}.	
	\end{align*}
	With the Young inequality we have the following estimate for all $n$ and $s \in [0,T]$:
	$$ \alpha\abs{Y_s^{n+1}}\leqslant \frac{\mu}{2} \abs{Y_s^{n+1}}^2+\frac{\alpha^2}{2 \mu},$$
	and thus
	\begin{align*}\abs{Y_t^{n+1}}^2 \leqslant & e^{-\mu(T-t)}\norme{\xi}_{L^\infty}^2+\espcond{\int_t^T 2e^{-\mu(s-T)}\parenth{\frac{\alpha^2}{2 \mu}+\gamma \abs{Y_s^{n+1}} \abs{Z_s^{n}}^2}\d s}{\mcF_t}\\&-\espcond{\int_t^T e^{-\mu(s-t)} \abs{Z_s^{n+1}}^2 \d s}{\mcF_t}\\\leqslant & \underbrace{e^{-\mu(T-t)}\norme{\xi}_{L^\infty}^2+\frac{\alpha^2}{\mu^2}\parenth{1-e^{-\mu(T-t)}}}_{\leqslant A^2}+2\gamma \times \espcond{\int_t^T e^{-\mu(s-t)}\abs{Y_s^{n+1}} \abs{Z_s^{n}}^2\d s}{\mcF_t}\\&-\espcond{\int_t^T e^{-\mu(s-t)} \abs{Z_s^{n+1}}^2 \d s}{\mcF_t}.
	\end{align*}
	Finally we obtain
	$$\abs{Y_t^{n+1}}^2 + \espcond{\int_t^T e^{-\mu(s-t)} \abs{Z_s^{n+1}}^2 \d s}{\mcF_t} \leqslant A^2 + 2\gamma \times \espcond{\int_t^T e^{-\mu(s-t)}\abs{Y_s^{n+1}} \abs{Z_s^{n}}^2\d s}{\mcF_t}.$$
	Then
	\begin{align} \norme{Y^{n+1}}_{\mcS^\infty}^2 &+ \supess \sup_{t \in [0,T]}  \espcond{\int_t^T e^{-\mu(s-t)} \abs{Z_s^{n+1}}^2 \d s}{\mcF_t} \nonumber\\\leqslant & 2 A^2+4\gamma \norme{Y^{n+1}}_{\mcS^\infty} \supess \sup_{t \in [0,T]}\espcond{\int_t^T e^{-\mu(s-t)}\abs{Z_s^{n}}^2\d s}{\mcF_t}\label{Ymono}\\ \leqslant & 2 A^2+\norme{Y^{n+1}}_{\mcS^\infty}^2+(2\gamma)^2 \supess \sup_{t \in [0,T]}\espcond{\int_t^T e^{-\mu(s-t)}\abs{Z_s^{n}}^2\d s}{\mcF_t}^2\nonumber .
	\end{align}
	Once again with an induction we show easily that for all $n \in \N$, $\abs{Z^n} \star W \in BMO$, $Y^n \in \mcS^\infty$ and 
	$$ \supess \sup_{t \in [0,T]}\espcond{\int_t^T e^{-\mu(s-t)}\abs{Z_s^{n}}^2\d s}{\mcF_t} \leqslant \frac{1}{8\gamma^2}\parenth{1-\sqrt{1-32\gamma^2 A^2}}.$$
	Moreover the inequality \eqref{Ymono} gives us
	$$ \norme{Y^{n+1}}_{\mcS^\infty} \leqslant 2\gamma \times \supess \sup_{t \in [0,T]}\espcond{\int_t^T e^{-\mu(s-t)}\abs{Z_s^{n}}^2\d s}{\mcF_t} +\sqrt{2A^2+4\gamma^2 \parenth{\supess \sup_{t \in [0,T]}\espcond{\int_t^T e^{-\mu(s-t)}\abs{Z_s^{n}}^2\d s}{\mcF_t} }^2}.$$
	Letting $n$ to infinity, and with the \cref{BMOconv}, we finally get 	
	$$ \supess \sup_{t \in [0,T]} \espcond{\int_t^T e^{-\mu(s-t)}\abs{Z_s}^2 \d s}{\mcF_t} \leqslant \frac{1}{8\gamma^2}\parenth{1-\sqrt{1-32\gamma^2 A^2}}.$$
	and so we deduce that
	$$ \norme{Y}_{\mcS^\infty} \leqslant \frac{1}{4\gamma}\parenth{1-\sqrt{1-32\gamma^2 A^2}}+\sqrt{2A^2+ \frac{1}{16 \gamma^2} \parenth{1-\sqrt{1-32\gamma^2 A^2}}^2}.$$
\end{Proof}

\begin{rema}\label{monBMO}
	What about the $BMO$ norm of $Z^M$ ? --- we can slice $[0,T]$ with a uniform sequence $(T_i)_{i=1}^N$ such that $0 = T_0 \leqslant T_1 \leqslant ... \leqslant T_N=T$ and for all $i$, $h=\abs{T_{i+1}-T_i}=\frac{T}{N}$. The last inequality can be used for each started and stopped process $\startstop{T_i}{\abs{Z^M} \star W}{T_{i+1}}$:
	\begin{align*}
		\supess \sup_{T_i \leqslant t \leqslant T_{i+1}} \espcond{\int_t^{T_{i+1}} \abs{Z^M_s}^2 \d s}{\mcF_t} \leqslant & e^{\mu h} \parenth{\supess \sup_{T_i \leqslant t \leqslant T_{i+1}} \espcond{\int_t^{T_{i+1}} e^{-\mu (s-t)} \abs{Z^M_s}^2 \d s}{\mcF_t}}\\\leqslant & e^{\mu h} \parenth{\supess \sup_{T_i \leqslant t \leqslant T_{i+1}} \espcond{\int_t^{T} e^{-\mu (s-t)} \abs{Z^M_s}^2 \d s}{\mcF_t}} \\\leqslant &  \frac{e^{\mu h}}{8\gamma^2}\parenth{1-\sqrt{1-32\gamma^2 A^2}}.
	\end{align*}
\end{rema}

We are now in position to prove \cref{prop Tevzadze monotone}.

\begin{Proof}[of \cref{prop Tevzadze monotone}]
The previous \cref{monBMO} shows that for all $M \in \R^+$ and all $h>0$, the process $\abs{Z^M} \star W$ is $$\parenth{\frac{e^{\frac12 \mu h}}{2\sqrt{2}\gamma}\parenth{1-\sqrt{1-32\gamma^2 A^2}}^{1/2}}- \text{scliceable}.$$
 We just have to apply an adaptation of \cref{withoutxib} given by \cref{remexistence-uniqueness2-avecepsBMO}. 
\end{Proof}

\subsection{Proof of \cref{proposition existence unicite cas diagonal}}

We consider the diagonal framework introduced in section \ref{main} and subsection \ref{Stabdiagsubsect}. We assume that the generator satisfies \ref{Hdiag}, so the generator $f$ can be written as $f =  f_{\tbf{diag}}(t,z)+g(t,y,z)$ where $f_{\tbf{diag}}$ is diagonal with respect to $z$. If we want to apply \cref{existuniq1diag}, we have to obtain a uniform estimate on $\norme{\abs{Z^{M}}\star W}_{BMO}$ where $(Y^M,Z^M)$ is the unique solution of the Lipschitz localized BSDE \eqref{quadraloc}. This is the purpose of the following lemma.
\begin{prop}
\label{estimee YM ZM cas diag}
 Let us assume that there exist nonnegative constants $G_d$ and $G$ such that
 \begin{enumerate}[label=(\roman*)]
  \item[\tbf{(i)}] for all $(t,y,z) \in [0,T] \times \R^d \times \R^{d \times k}$,  $\quad\abs{f_{\tbf{diag}}(t,z)} \leqslant G_d \abs{z}^2, \quad \abs{g(t,y,z)} \leqslant G \abs{z}^2.$
  \item[\tbf{(ii)}] \beq \label{relation Gd G} \frac{4\sum_{i=1}^d e^{2G_d \norme{\xi^i}_{L^{\infty}}}}{G_d}G \leqslant 1.\eeq \end{enumerate}
Then, $\abs{Z^M} \star W \in BMO$, $Y^M \in \mcS^\infty$ and we have following estimates:
$$\norme{\abs{Z^{M}}\star W}_{BMO} \leqslant (4G_dG)^{-1/2}, \quad \norme{Y^M}_{\mcS^{\infty}} \leqslant \norme{\xi}_{L^\infty} +\frac{\sqrt{d}\log 2}{2G_d}.$$

\end{prop}

\begin{Proof}
To simplify notations in the proof, we skip once again the superscript $M$ on $\parenth{Y^M,Z^M}$ and $f^M$. The unique solution $(Y,Z) \in \mcS^2\parenth{\R^d} \times \mcH^2\parenth{\R^{d\times k}}$ of \eqref{quadraloc} can be constructed with a Picard principle slightly different than the one used in the seminal paper of Pardoux and Peng (see \cite{Pardoux:1990ju}). We consider a sequence $\suiteN{Y^n,Z^n}{n}$ defined by
$$Y^{n+1}_t=\xi +\int_t^T f_{\tbf{diag}}(s,Z^{n+1}_s)+ g\parenth{s,Y^{n}_s,Z^{n}_s} \d s-\int_t^T Z^{n+1}_s \d W_s, \quad 0 \leqslant t \leqslant T, \quad (Y^0,Z^0)=(0,0).$$
Obviously, we can easily show that $\suiteN{Y^n,Z^n}{n}$ tends to $(Y,Z \star W)$ in $\mcS^2\parenth{\R^d} \times \mcH^2 \parenth{\R^{d\times k}}$ since we are in the Lipschitz framework. 
We will prove by induction that: for all $n \in \N$, $Y^n \in \mcS^\infty$, $\abs{Z^n} \star W \in BMO$ and
\begin{equation}
\label{hypothese recurrence cas diagonal}
\norme{\abs{Z^{n}}\star W}_{BMO} \leqslant (4G_dG)^{-1/2}, \quad \norme{Y^n}_{\mcS^{\infty}} \leqslant \norme{\xi}_{L^\infty} +\frac{\sqrt{d}\log 2}{2G_d}.
\end{equation}
The result is obvious for $n=0$. Let us assume that for a given $n \in \N$ we have $Y^n \in \mcS^\infty\parenth{\R^d}$, $\abs{Z^n} \star W \in BMO$ and \eqref{hypothese recurrence cas diagonal} is true. The Lipschitz regularity of $f$ gives us that $(Y^{n+1},\abs{Z^{n+1}}\star W) \in \mcS^{\infty}\parenth{\R^d} \times BMO$. By following the idea of \cite{Briand:2006he}, we introduce the $C^2$ function $\varphi : \R \longmapsto (0,+\infty)$ defined by 
$$\varphi : x \mapsto \frac{e^{2G_d\abs{x}}-1-2G_d\abs{x}}{\abs{2G_d}^2}, \quad \text{so that} \quad \varphi''(.)-2G_d\abs{\varphi'(.)} = 1.$$
We pick a stopping time $\tau$ and applying It\^o's formula to the regular function $\phi$, we compute, for all components $i \in \{1,...,d\}$,
\begin{align*}
 \varphi\parenth{Y_{\tau}^{n+1,i}} =& \varphi\parenth{\xi^i} +\int_{\tau}^T \left( \varphi'\parenth{Y_s^{n+1,i}} \parenth{f_{\tbf{diag}}(s,Z_s^{n+1,i})+g(s,Y^{n,i}_s,Z_s^{n,i})}-\frac{\varphi''\parenth{Y_s^{n+1,i}}\abs{Z_s^{n+1,i}}^2}{2} \right) \d s 
 \\ &- \int_{\tau}^T \varphi'\parenth{Y_s^{n+1,i}} Z_s^{n+1,i} \d W_s\\
 \leqslant &  \varphi\parenth{\norme{\xi^i}_{L^{\infty}}} + \int_\tau^T \parenth{ \frac{2G_d \abs{\varphi'\parenth{Y_s^{n+1,i}}}-\varphi''\parenth{Y_s^{n+1,i}}}{2} \abs{Z_s^{n+1,i}}^2 + G  \abs{\varphi'(Y_s^{n+1,i})} \abs{Z_s^{n}}^2} \d s\\
 & - \int_{\tau}^T \varphi'(Y_s^{n+1,i}) Z_s^{n+1,i} \d W_s.
\end{align*}
Since $\varphi''(.)-2G_d\abs{\varphi'(.)} = 1$, $\varphi \geqslant 0$ and $\abs{\varphi'(x)} \leqslant (2G_d)^{-1} e^{2G_d\norme{Y^{n+1,i}}_{\mcS^{\infty}}}$ whenever $\abs{x} \leqslant \norme{Y^{n+1,i}}_{\mcS^{\infty}}$, taking the conditional expectation with respect to $\mathcal{F}_{\tau}$, we compute
$$\frac{1}{2}  \espcond{ \int_\tau^T \abs{Z_s^{n+1,i}}^2 \d s}{\mcF_\tau} \leqslant \varphi(\norme{\xi^i}_{L^{\infty}}) + G(2G_d)^{-1} e^{2G_d\norme{Y^{n+1,i}}_{\mcS^{\infty}}} \norme{\abs{Z^{n}}\star W}_{BMO}^2.$$
Thus, we get the estimate
\begin{equation}
\label{estimee BMO Zn+1 cas diag}
\frac{1}{2} \norme{\abs{Z^{n+1,i}}\star W}_{BMO}^2 \leqslant \varphi(\norme{\xi^i}_{L^{\infty}}) + G(2G_d)^{-1} e^{2G_d\norme{Y^{n+1,i}}_{\mcS^{\infty}}} \norme{\abs{Z^{n}}\star W}_{BMO}^2.
\end{equation}
By using the \tit{a priori} estimate given by Proposition 1 in \cite{Briand:2006he} we also have, for all stopping time $\tau$ and $i \in \{1,...,d\}$,
$$e^{2 G_d\abs{Y^{n+1,i}_\tau}} \leqslant e^{2G_d \norme{\xi^i}_{L^{\infty}}}\espcond{e^{2G_d G\int_\tau^T \abs{Z^{n}_s}^2 \d s}}{\mcF_\tau}.$$
Then, the John-Nirenberg inequality \eqref{JohnNir} coupled with the induction assumption on $Z^n$ gives us
$$e^{2G_d\norme{Y^{n+1,i}}_{\mcS^{\infty}}} \leqslant \frac{e^{2G_d \norme{\xi^i}_{L^{\infty}}}}{1-2G_dG \norme{\abs{Z^{n}}\star W}_{BMO}^2}.$$
We put this last inequality into \eqref{estimee BMO Zn+1 cas diag} to obtain
$$ \norme{\abs{Z^{n+1}}\star W}_{BMO}^2 \leqslant \frac{\sum_{i=1}^d e^{2G_d \norme{\xi^i}_{L^{\infty}}}}{2G_d^2} \frac{1}{1-2G_dG \norme{\abs{Z^{n}}\star W}_{BMO}^2}.$$
Since we have assumed that $\norme{\abs{Z^{n}}\star W}_{BMO}^2 \leqslant (4G_dG)^{-1}$ and \eqref{relation Gd G} is fulfilled, we get
$$\norme{\abs{Z^{n+1}}\star W}_{BMO}^2 \leqslant (4G_dG)^{-1}$$
and, by using previous calculations,
$$\norme{Y^{n+1,i}}_{\mcS^{\infty}} \leqslant \norme{\xi^i}_{L^\infty} +\frac{\log 2}{2G_d}, \quad \text{and so} \quad \norme{Y^{n+1}}_{\mcS^{\infty}} \leqslant \norme{\xi}_{L^\infty} +\frac{\sqrt{d}\log 2}{2G_d},$$
which concludes the induction. Finally, we just have to use the fact that $\suiteN{Y^n,Z^n\star W}{n}$ tends to $(Y,Z \star W)$ in $\mcS^2\parenth{\R^d} \times \mcH^2 \parenth{\R^{d\times k}}$ and \cref{BMOconv} to obtain the desired result.
\end{Proof}

Then the proof of \cref{proposition existence unicite cas diagonal} is direct: we just have to apply \cref{existuniq1diag} by using \cref{estimee YM ZM cas diag}.
 
\begin{rema}
 As explained in \cref{extension existence unicite cas diagonal}, it is possible to extend \cref{proposition existence unicite cas diagonal} to more general terminal conditions and generators.
\end{rema}

\subsection{Proof of \cref{theo existence unicite martingale variete}}

\begin{Proof}[of \cref{theo existence unicite martingale variete}]
Let us consider for all $M \in \R^+$ a smooth map $h^M : \R^{d\times k} \rightarrow \R$ satisfying:
$$ h^M(z)=\begin{cases}
	0 & \text{if } \abs{z} \leqslant M, \\  \parenth{\frac{\abs{z}}{M+1}}^2-1 &\text{if } \abs{z} \geqslant M
+1,\end{cases}$$
and let us define a localisation $\rho^M : \R^{d \times k} \rightarrow \R^{d \times k}$ given by
$$ \rho^M(z)=\frac{z}{\sqrt{1+h^M(z)}}.$$
%
As usual we denote by $(Y^M,Z^M)$ the solution obtained by replacing $f$ by $f^M$. This choice of $\rho^M$ will be useful in the following computations. 
For $F \in \mcC^2 (\mcM,\R)$, the Itô formula with $F$ seen as a function on $\R^d$ gives for all stopping time $\tau$:
\begin{align*} \espcond{F(\xi)-F(Y^M_\tau)}{\mcF_\tau}&=\espcond{\int_{\tau}^T \left( \frac{1}{2} \sum_{l=1}^k \Hess F(Y^M_s)\left(Z_s^{M, (:,l)},Z_s^{M,(:,l)}\right)-\d F\left(Y^M_s\right) f^M\left(Y^M_s,Z^M_s\right)\right) \d s}{\mcF_\tau},
\end{align*}
since the local martingale part is a martingale because $F$ has bounded first derivative.
By using the definition of $f$, its formulation in the local chart and the link between $Z^M$ and $\rho^M(Z^M)$, we get 

\begin{align} \label{Itoconv}
{}&\espcond{F(\xi)-F(Y^M_\tau)}{\mcF_\tau}\nonumber \\&= \frac{1}{2}\espcond{\int_{\tau}^T \sum_{l=1}^k\left(h^M(Z^M_s) \Hess F(Y^M_s)\left(\rho^M(Z_s^{M})^{(:,l)},\rho^M(Z_s^{M})^{(:,l)}\right)+\nabla \d F\left(Y^M_s\right)\left(\rho^M(Z_s^{M})^{(:,l)},\rho^M(Z_s^{M})^{(:,l)}\right) \right) \d s}{\mcF_\tau}.
\end{align}

By using same arguments as Darling in \cite{Darling:1995wd} we can show that $Y^M$ is in $G$ almost surely.
Indeed, we know with \eqref{Itoconv} applied to $F=F^{dc}$, integrating between $\tau$ and $\sigma$ together with \tbf{(ii)}, that for all $\sigma \geqslant \tau \ps$, 
	$$ \espcond{F^{dc}(Y^M_\sigma)}{\mcF_\tau} \geqslant F^{dc}(Y^M_\tau)\ps.$$
	Let us consider for all $n \in \N$ the following sequence of stopping times: $\sigma^n=\inf\ensemble{\bigmid{u \geqslant \tau}{F^{dc}(Y^M_{u})\leqslant \frac{1}{n} \ps }}.$
	Each $\sigma^n$ is finite almost surely since $\xi \in G$. Continuity of $Y^M$ gives for all $n \in \N$,  $F^{dc}\parenth{Y^M_{\sigma^n}} \leqslant \frac{1}{n} \ps$. So we get for all stopping time $\tau$:
	$$ F^{dc}\parenth{ Y^M_\tau }\leqslant \espcond{ F^{dc}\parenth{Y^M_{\sigma^n}} }{\mcF_\tau} \leqslant \frac{1}{n} \quad  \ps.$$
and, consequently, $\PP(Y_t \in G)=1$ for all $t \in [0,T]$. Moreover the $\alpha$-strictly doubly convexity on $G$ gives us

\begin{align}
	&{}\espcond{F^{dc}(\xi)-F^{dc}(Y^M_\tau)}{\mcF_\tau} \nonumber\\\geqslant &\frac{\alpha}{2}\espcond{\int_{\tau}^T \sum_{l=1}^k\left(h^M(Z^M_s)\abs{\rho^M(Z_s^{M})^{(:,l)}}^2+\abs{\rho^M(Z_s^{M})^{(:,l)}}^2\right) \d s}{\mcF_\tau}=\frac{\alpha}{2}  \espcond{\int_{\tau}^T \abs{Z^M_s}^2 \d s}{\mcF_\tau}.
\end{align} 

And finally, continuity of $F^{dc}$ on $G$ yields
	$$\norme{\abs{Z^M} \star W}_{BMO} \leqslant \sqrt{\frac{2} {\alpha} \times \parenth{\sup_{(x,y) \in G \times G} \ensemble{F^{dc}(x)-F^{dc}(y)}}}.$$
	Thus, assumption \tbf{(ii)} ensures assumption \ref{BMO,m}. Since the terminal value is bounded (in $G$), \cref{withoutxib} together with  \cref{simplifié} gives the result.
\end{Proof}

\subsection{Proofs of \cref{th unicite cas markovien} and \cref{th regularite markov}}

\begin{Proof}[of \cref{th unicite cas markovien} and \cref{th regularite markov}]

\paragraph{Uniqueness} We start by proving \cref{th unicite cas markovien}. Let us consider two continuous Markovian solutions $(v,w)$ and $(\tilde{v},\tilde{w})$ such that $v$ and $\tilde{v}$ are bounded. We set $(t,x) \in [0,T] \times \R^k$ and we denote $(Y^{t,x},Z^{t,x})$ (resp. $(\tilde{Y}^{t,x},\tilde{Z}^{t,x})$) the solution of the BSDE \eqref{XZquadra} associated to  $(v,w)$ (resp. $(\tilde{v},\tilde{w})$). The idea of the proof is to compare the two solutions by using the stability result given by Remark  \ref{remarque stabilite sliceable}. In order to do that, we must show that $\abs{Z^{t,x}} \star W$ and $\abs{\tilde{Z}^{t,x}} \star W$ are $\varepsilon$-sliceable BMO martingales. 
By Remark 2.6 part (2) in \cite{2016arXiv160300217X}, we know that there exist $b_0 \in \R^d$ and $(\alpha_n) \in (0,1]^{\N}$ such that $v,\tilde{v}\in\mcC^{(\alpha_n),loc}_{b_0}$. Moreover, by following same arguments than in the proof of Theorem 2.9 in \cite{2016arXiv160300217X} we can show that $v,\tilde{v}\in\mcC^{(\alpha_n),loc}$. Finally, we just have to apply \cref{loca Hold et bornée implique glob Hold} to conclude that there exists $\kappa' \in (0,1]$ such that $v,\tilde{v}\in\mcC^{\kappa'}$. 

Now, let us apply the Itô formula to $F(Y^{t,x})$: we consider two stopping times $\tau$ and $\sigma$ such that $\tau \leqslant \sigma \ps$ and we take the conditional expectation
	\begin{align*}
		{}&\espcond{F\parenth{v(\sigma, X^{t,x}_\sigma)}-F\parenth{v(\tau, X^{t,x}_\tau)}}{\mcF_\tau}\\&=\espcond{-\int_{\tau}^\sigma \d F(v(u,X^{t,x}_u) \tbf{f}(u,X^{t,x}_u,Y^{t,x}_u,Z^{t,x}_u) \d u+\frac12 \sum_{l=1}^k \int_\tau^\sigma \Hess F\parenth{v(u,X^{t,x}_u)}\parenth{\parenth{Z^{t,x}_u}^{(:,l)}, \parenth{Z^{t,x}_u}^{(:,l)}}\d u }{\mcF_\tau}\\& \geqslant \espcond{\int_\tau^\sigma \abs{Z^{t,x}_u}^2 \d u}{\mcF_\tau}.
	\end{align*}

The map $F$ is a Lipschitz function on the centred Euclidean ball of radius $\norme{v}_{L^{\infty}([0,T] \times \R^k)}\vee \norme{\tilde{v}}_{L^{\infty}([0,T] \times \R^k)}$. Denoting by $L$ its Lipschitz constant, we obtain for all $n \in \N$, 
\begin{align*}
\espcond{\int_\tau^\sigma \abs{Z^{t,x}_u}^2 \d u}{\mcF_\tau} &\leqslant L \times \espcond{\abs{v(\sigma, X^{t,x}_\sigma) - v(\tau, X^{t,x}_\tau)}}{\mcF_\tau}\\&\leqslant L \times \espcond{\abs{\sigma-\tau}^{\kappa'/2}+\abs{X^{t,x}_\sigma-X^{t,x}_\tau}^{\kappa'}}{\mcF_\tau}\\
& \leqslant C\abs{\sigma-\tau}^{\kappa'/2},
\end{align*}

where we have used in the last inequality a classical estimate for SDEs when $b$ and $\sigma$ are bounded.   
For $N \in \N^*$ we set $T_i = \frac{iT}{N}$. Then, for all $i \in \{0,...,N-1\}$ and stopping time $T_{i} \leqslant \tau \leqslant T_{i+1}$ we get
$$\espcond{\int_\tau^{T_{i+1}} \abs{Z^{t,x}_u}^2 \d u}{\mcF_\tau} \leqslant C\parenth{\frac{T}{N}}^{\kappa'/2},$$
and finally, for $N$ large enough, we have that 
$$\norme{\startstop{T_i}{\abs{Z^{(t,x)}}\star W}{T_{i+1}} }_{BMO} \leqslant \mbbK, \quad \forall i \in \{0,...,N-1\}$$
with $2 L_y \mbbK^2+2\sqrt{2} L_z\mbbK C'_2<1$. Obviously, this estimate is also true for $\tilde{Z}^{t,x}$ which means that $Z^{t,x} \star W$ and $\tilde{Z}^{t,x}\star W$ are in $\mcZ_{BMO}^{slic,2}$. By using \cref{remarque stabilite sliceable} we get that $v(t,x)=\tilde{v}(t,x)$ and $\esp{\int_t^T |w(s,X_s^{t,x})-w'(s,X_s^{t,x})|^2 \d s}=0$. Since this is true for all $(t,x) \in [0,T] \times \R^k$ and, due to \ref{HX}, $X_s^{0,x}$ has positive density on $\R^k$ for $s \in (0,T]$, $w=\tilde{w}$ a.s. with respect to the Lebesgue measure on $[0,T] \times \R^k$. Then, $(v,w)$ and $(v',w')$ are equal.

\paragraph{Regularity}
Now we prove \cref{th regularite markov}.  We consider $(v,w)$ a continuous Markovian solution of \eqref{XZquadra} such that $v$ is bounded. 
We set $t,t' \in [0,T]$ and $x,x' \in \R^k$. Without restriction, we can assume that $t \leqslant t'$. Then, we can write
$$\abs{v(t,x)-v(t',x')} = \abs{Y_t^{t,x} - Y_{t'}^{t',x'}} \leqslant \abs{Y_t^{t,x} - Y_{t}^{t',x'}}+\abs{Y_t^{t',x'} - Y_{t'}^{t',x'}},$$
where we have extend the definition of $(Y_u^{t',x'},Z_u^{t',x'})_{t' \leqslant u \leqslant T}$ to $[0,t']$ by setting
$$Y_s^{t',x'} = Y_{t'}^{t',x'} + \int^{t'}_s \tbf{f}\parenth{u,x',Y_u^{t',x'},0} \d u, \quad Z_s^{t',x'}=0,\quad 0 \leqslant s \leqslant t',$$
which is equivalent to take $\sigma(s,.)=0$ and $b(s,.)=0$ for $s \in [0,t']$. By the same token we can extend $(Y_u^{t,x},Z_u^{t,x})_{t \leqslant u \leqslant T}$ to $[0,t]$.
Then, a standard estimate gives us 
$$\abs{Y_t^{t',x'} - Y_{t'}^{t',x'}} \leqslant C\abs{t-t'},$$
where $C$ does not depend on $(t,x)$ and $(t',x')$.
To conclude we just have to study the remaining term $\abs{Y_t^{t,x} - Y_{t}^{t',x'}}$. Thanks to calculations done in the uniqueness part of the proof, we know that there exist some deterministic times $0=T_0 \leqslant T_1 \leqslant ... \leqslant T_N=T$ such that $ \startstop{T_i}{Z^{t,x}\star W}{T_{i+1}} \in \mcZ_{BMO}^{2}$ and $ \startstop{T_i}{Z^{t',x'}\star W}{T_{i+1}} \in \mcZ_{BMO}^{2}$ for all $i\in \{0,...,N-1\}$. Let us emphasize that times $0=T_0 \leqslant T_1 \leqslant ... \leqslant T_N=T$ can be chosen independently from $(t,x)$ and $(t',x')$. Then, by using once again \cref{remarque stabilite sliceable} and assumptions on $\mcG$ and $\tbf{f}$, we obtain
\begin{align*}
\abs{Y_t^{t,x} - Y_{t}^{t',x'}} &\leqslant \norme{Y^{t,x} - Y^{t',x'}}_{\mcS^{4}}\\
 & \leqslant C\left( \esp{\abs{\mcG \parenth{X^{t,x}_T}-\mcG \parenth{X^{t',x'}_T}}^4}^{1/4} + \esp{\abs{ \int_0^T \abs{ \tbf{f}\parenth{s,X_s^{t,x},Y_s^{t,x},Z_s^{t,x}} - \tbf{f}\parenth{s,X_s^{t',x'},Y_s^{t,x},Z_s^{t,x}} } \d s }^4}^{1/4}\right)\\
 & \leqslant C\left( \esp{\abs{X^{t,x}_T-X^{t',x'}_T}^{4\kappa}}^{1/4} + \esp{ \sup_{s \in [0,T]} \abs{X_s^{t,x} - X_s^{t',x'}}^{4\kappa} \left(1+ \int_0^T \abs{Z_s^{t,x}}^2 \d s \right)^4 }^{1/4}\right)\\
  & \leqslant C \esp{\abs{X^{t,x}_T-X^{t',x'}_T}^{4\kappa}}^{1/4} + C\esp{ \sup_{s \in [0,T]} \abs{X_s^{t,x} - X_s^{t',x'}}^{8\kappa}}^{1/8} \left(1+ \esp{\left(\int_0^T \abs{Z_s^{t,x}}^2 \d s \right)^8 }^{1/8}\right).
\end{align*}
where $C$ depends on $N$ and the BMO norms of $\startstop{T_i}{Z^{t,x}\star W}{T_{i+1}}$  (see an explicit value in \cref{remarque stabilite sliceable}). However $N$ stays finite in the sequel and, as already mentioned in the Proof of \cref{existuniq1} and \cref{existuniq1diag}, an uniform estimation on $\norme{\startstop{T_i}{Z^{t,x}\star W}{T_{i+1}}}_{BMO}$ and $\norme{\startstop{T_i}{Z^{t',x'}\star W}{T_{i+1}}}_{BMO}$ gives an uniform upper bound on $C$ (particularly $C$ does not depend on $(t,x)$ and $(t',x')$). Then an energy inequality gives us 
$$|Y_t^{t,x} - Y_{t}^{t',x'}| \leqslant C\esp{ \sup_{s \in [0,T]} \abs{X_s^{t,x} - X_s^{t',x'}}^{8\kappa}}^{1/8}.$$
Thus, we just have to use the classical estimate on SDEs given by 
$$ \norme{X^{t',x'}-X^{t,x}}_{\mcS^8} \leqslant C \parenth{\abs{x-x'} + \abs{t-t'}^{1/2}},$$
to obtain that 
$$|Y_t^{t,x} - Y_{t}^{t',x'}| \leqslant C \parenth{\abs{x-x'}^{\kappa} + \abs{t-t'}^{\kappa/2}}.$$
\end{Proof}

\paragraph{Acknowledgements.} The authors would like to thank Hao Xing for useful remarks and comments on a preliminary version of this article.

\section{Appendix – Technical proofs}

\begin{Proof}[of \cref{BMOstartstop}]
	Let us show the $BMO$ property for a started and stopped process. Let us consider a stopping time $\tau'$ such that $0 \leqslant \tau' \leqslant T \ps.$ we have

\begin{align*} \espcond{\quadra{\startstop{\tau}{\abs{Z} \star W}{\sigma}}_T - \quadra{\startstop{\tau}{\abs{Z} \star W}{\sigma}}_{\tau'}}{\mcF_{\tau'}}&=\espcond{\parenth{\quadra{{}^\tau \abs{Z} \star W}_\sigma - \quadra{{}^\tau \abs{Z} \star W}_{\min(\tau', \sigma)}}}{\mcF_{\tau'}}\\
&=\espcond{\parenth{ \quadra{{}^\tau \abs{Z} \star W}_\sigma - \quadra{{}^\tau \abs{Z} \star W}_{\tau'}}\indicat_{(0 \leqslant \tau' \leqslant \sigma)}}{\mcF_{\tau'}}.\end{align*}

Since ${}^\tau \abs{Z} \star W$ vanishes before $\tau$ and $ \mcF_{\tau'} \subset \mcF_{\max(\tau',\tau)}$, we get:
\begin{align*} \espcond{\parenth{ \quadra{{}^\tau \abs{Z} \star W}_\sigma - \quadra{{}^\tau \abs{Z} \star W}_{\tau'}}\indicat_{(0 \leqslant \tau' \leqslant \sigma)}}{\mcF_{\tau'}} =&\espcond{\parenth{ \quadra{\abs{Z} \star W}_\sigma - \quadra{\abs{Z} \star W}_{\max(\tau',\tau)}}\indicat_{(0 \leqslant \tau' \leqslant \sigma)}}{\mcF_{\tau'}} \\=&\espcond{\espcond{\parenth{\quadra{\abs{Z} \star W}_\sigma - \quadra{\abs{Z} \star W}_{\max(\tau',\tau)}}\indicat_{(0 \leqslant \tau' \leqslant \sigma)}}{\mcF_{\max(\tau,\tau')}}}{\mcF_{\tau'}} 	\\\leqslant & \supess \sup_{\tilde{\tau} \in \mcT^{\tau,\sigma}} \espcond{\quadra{\abs{Z} \star W}_{\sigma}-\quadra{\abs{Z} \star W}_{\tilde{\tau}}}{\mcF_{\tilde{\tau}}}.
	\end{align*}
 Finally we have shown that 
 $$\supess \sup_{0 \leqslant \tau' \leqslant T} \espcond{\quadra{\startstop{\tau}{\abs{Z} \star W}{\sigma}}_T - \quadra{\startstop{\tau}{\abs{Z} \star W}{\sigma}}_{\tau'}}{\mcF_{\tau'}} \leqslant \supess \sup_{\tilde{\tau} \in \mcT^{\tau,\sigma}} \espcond{\quadra{\abs{Z} \star W}_{\sigma}-\quadra{\abs{Z} \star W}_{\tilde{\tau}}}{\mcF_{\tilde{\tau}}},$$
	and the inequality is obviously an equality.\end{Proof}

\begin{Proof}[of \cref{SDEprop}]
	We are going to use inequalities given by Lemma \ref{BMOineq}. Let us suppose the existence of a solution $X$ for the equation \eqref{SDE}. We have for all $m \geqslant 1$,
	\begin{align*}
	\norme{X}_{\mcS^m} \leqslant & \norme{X_0}_{L^m}+\norme{\int_0^. F(s,X_s) \d s}_{\mcS^m}+\norme{\sum_{p=1}^k \int_0^. G^p(s,X_s) \d W^p_s}_{\mcS^m} 	
	\\\leqslant & \norme{X_0}_{L^m}+\esp{\sup_{0 \leqslant u \leqslant T}\left(\int_0^{u} \alpha_s \abs{X_s} \d s\right)^m}^{1/m}+C'_m\esp{\left(\sum_{p=1}^k\int_0^T \abs{G^p(s,X_s)}^2 \d s\right)^{m/2}}^{1/m}
	\\\leqslant & \norme{X_0}_{L^m}+\esp{\left(\int_0^{T} \alpha_s \abs{X_s} \d s\right)^m}^{1/m}+C'_m\esp{\left(\int_0^T \beta_s^2 \abs{X_s}^2 \d s\right)^{m/2}}^{1/m}.
	\end{align*}
	On the one hand, according to Lemma \ref{BMOineq} we have 
\begin{align*}
	\esp{\left(\int_0^{T} \alpha_s \abs{X_s} \d s\right)^m}^{1/m}&=\norme{\scalaire{\sqrt{\alpha} \star W}{(\sqrt{\alpha}\abs{X})\star W}_T}_{L^m} \\\leqslant & \sqrt{2} m \norme{\sqrt{\alpha}\star W}_{BMO}\norme{(\sqrt{\alpha}\abs{X}) \star W}_{\mcH^m}\\\leqslant & 2 m\norme{X}_{\mcS^m} \norme{\sqrt{\alpha}\star W}_{BMO}^2.
\end{align*}
	On the other hand, we get for the last term
	$$ \esp{\left( \int_0^T \beta_s^2 \abs{X_s}^2 \d s\right)^{m/2}}^{1/m}=\norme{(\beta |X|) \star W}_{\mcH^m}=\norme{|X| \star(\beta \star W)}_{\mcH^m}\leqslant \sqrt{2} \norme{X}_{\mcS^m} \norme{\beta \star W}_{BMO}.$$
	Hence we obtain the following inequality 
	\beq \norme{X}_{\mcS^m}\left( 1-2 m \norme{\sqrt{\alpha}\star W}_{BMO}^2-\sqrt{2}C'_m\norme{\beta \star W}_{BMO}\right)\leqslant \norme{X_0}_{L^m}.\label{XSm} \eeq
	The constant behind $\norme{X}_{\mcS^m}$ is not always positive, but we can use the sliceability assumption in order to construct piece by piece the process $X$, and on each piece the constant will be positive. 
	
	More precisely there exists a sequence of stopping times $0=T_0 \leqslant T_1 \leqslant ... \leqslant T_N=T \ps$ such that for all $i \in \ensemble{0,...,N-1}$:
	$$ \norme{\startstop{T_i}{\sqrt{\alpha}\star W}{T_{i+1}}}_{BMO} \leqslant \varepsilon_1, \quad \norme{{}^{T_i}\beta\star W^{\rfloor T_{i+1}}}_{BMO} \leqslant \varepsilon_2. $$
	 The process $X$ is equal to $$X_t=\sum_{i=1}^{N-1} \widetilde{X^i_t} \indicat_{[T_i,T_{i+1}[}(t)$$
	 where each $\widetilde{X^i}$ is the restriction of $X$ to the stochastic interval $[T_i,T_{i+1}]$. By convention we extend $\widetilde{X^i}$ to $[0,T]$ by zero outside $[T_{i},T_{i+1}]$. 
	 $\widetilde{X_i}$ satisfies the following SDE:
	 $$\widetilde{X^i_t} = \widetilde{X_{T_i}}^{i-1}+ \int_{T_i}^t F(s,\widetilde{X^i_s}) \mathrm{d}s + \sum_{p=1}^k \int_{T_i}^t G^p(s,\widetilde{X^i_s}) \mathrm{d} \startstop{T_i}{W_s}{T_{i+1}}, \quad t \in [T_i,T_{i+1}[, \quad \text{ and }\quad \widetilde{X}^{-1}=X_0.$$
	 For all $i \in \ensemble{0,...,N-1}$, by considering above computations on each $[T_i, T_{i+1}[$, \eqref{XSm} becomes
	$$\norme{\widetilde{X^i}}_{\mcS^m}\left( 1-2m\varepsilon^2_1 -\varepsilon_2\sqrt{2}C'_m \right)\leqslant \norme{\widetilde{X_{T_i}}^{i-1}}_{L^m}.$$
	Denoting by $\KK_{\varepsilon_1,\varepsilon_2}$ the constant
	$$\KK_{\varepsilon_1,\varepsilon_2}:=\frac{1}{1-2m\varepsilon_1^2-\varepsilon_2\sqrt{2}C'_m}>0,$$
	we have
	$$ \norme{\widetilde{X^i}}_{\mcS^m} \leqslant \KK_{\varepsilon_1,\varepsilon_2}^i \norme{X^0}_{L^m},$$
	and finally we obtain
	$$\norme{X}_{\mcS^m}\leqslant \sum_{i=0}^{N-1} \norme{\widetilde{X^i}}_{\mcS^m}\leqslant  \parenth{\sum_{i=0}^{N-1} \KK_{\varepsilon_1,\varepsilon_2}^{i}}\norme{X^0}_{L^m}.$$
	The result follows by setting $K_{m,\varepsilon_1,\varepsilon_2}=\displaystyle\sum_{i=0}^{N-1} \KK_{\varepsilon_1,\varepsilon_2}^{i}$.
\end{Proof}

\pagebreak

\pagebreak

\bibliographystyle{alpha}

\end{document}